\documentstyle[11pt,amstex,amsfonts,amssymb]{amsart}
\ifx\shlhetal\undefinedcontrolsequence\let\shlhetal\relax\fi

\newcount\skewfactor
\def\mathunderaccent#1#2 {\let\theaccent#1\skewfactor#2
\mathpalette\putaccentunder}
\def\putaccentunder#1#2{\oalign{$#1#2$\crcr\hidewidth
\vbox to.2ex{\hbox{$#1\skew\skewfactor\theaccent{}$}\vss}\hidewidth}}
\def\name{\mathunderaccent\tilde-3 }

\newcommand{\rest}{{\,|\grave{}\,}}
\newcommand{\forces}{\mathrel {{\vrule height 6.9pt depth -0.1pt}\! \vdash }}

\newcommand{\V}{{\bf V}} 
\newcommand{\lesdot}{\mathrel{\mathord{<}\!\!\raise 
0.8 pt\hbox{$\scriptstyle\circ$}}} 
   
\newcommand{\conc}{{}^\frown\!}
\newcommand{\lh}{\ell g\/} 
\newcommand{\dom}{{\rm Dom}} 
\newcommand{\rng}{{\rm Rang}}

\newcommand{\bP}{{\Bbb P}}
\newcommand{\bQ}{{\Bbb Q}}
\newcommand{\bR}{{\Bbb R}}
\newcommand{\nbQ}{{\name{\Bbb Q}}}

\newcommand{\QED}{\hfill\vrule width 6pt height 6pt depth 0pt 
\vspace{0.1in}} 
\newcommand{\Proof}{\noindent{\sc Proof} \hspace{0.2in}} 

\newcommand{\can}{{}^{\textstyle \omega}2} 
 
\newcommand{\baire}{{}^{\textstyle \omega}\omega}

\newcommand{\rk}{{\rm rk}}
\newcommand{\Lim}{{\rm Lim}}
\newcommand{\pd}{{\rm pd}}
\newcommand{\cl}{{\rm cl}}
\newcommand{\ZFC}{{\rm ZFC}}
\newcommand{\bp}{{\bf p}}
\newcommand{\proj}{{\rm proj}}
\newcommand{\code}{{\rm code}}
\newcommand{\Levy}{{\rm Levy}}
\newcommand{\truth}{{\frak t}}
\newcommand{\hc}{{\rm hc}}
\newcommand{\pr}{{\rm pr}}
\newcommand{\cK}{{\cal K}}
\newcommand{\ex}{{\rm ex}}
\newcommand{\Borel}{{\rm Borel}}
\newcommand{\Cohen}{{\rm Cohen}}
\newcommand{\BA}{{\rm BA}}
\newcommand{\dx}{{\rm dx}}
\newtheorem{theorem}{Theorem}[section] 
\newtheorem{claim}{Claim}[theorem]
\newtheorem{fact}[theorem]{Fact}
\newtheorem{lemma}[theorem]{Lemma} 
\newtheorem{proposition}[theorem]{Proposition} 
 
\newtheorem{thesis}[theorem]{Thesis} 

\theoremstyle{definition}
\newtheorem{definition}[theorem]{Definition}
\newtheorem{example}[theorem]{Example}
\newtheorem{question}[theorem]{Question}
\newtheorem{defthe}[theorem]{Definition/Theorem}
\newtheorem{problem}[theorem]{Problem} 

\theoremstyle{remark}
\newtheorem{notation}[theorem]{Notation}
\newtheorem{conclusion}[theorem]{Conclusion}
\newtheorem{remark}[theorem]{Remark}

\newtheorem{discussion}[theorem]{Discussion}
\newtheorem{convention}[theorem]{Convention}
\newtheorem{hypothesis}[theorem]{Hypothesis}
\newtheorem{choice}[theorem]{Choice}
\newtheorem{comment}[theorem]{Comment}

\title{Properness Without Elementaricity}  
\author[S. Shelah]{\uppercase {\bf S. Shelah}}
\address{Institute of Mathematics\\
The Hebrew University\\
Jerusalem 91904, Israel\\
and Rutgers University\\
Mathematics Department\\
New Brunswick, NJ 08854, USA
}
\email{shelah@@math.huji.ac.il}

\date{\today} 

\thanks{The research was partially supported by ``Basic Research Foundation''
of the Israel Academy of Sciences and Humanities. Publication 630.}


\begin{document} 
\maketitle 

\bigskip
\bigskip
\bigskip
{\it Abstract.}
We present reasons for developing a theory of forcing notions
which satisfy the properness demand for countable models which are not
necessarily elementary submodels of some $({\cal H}(\chi),\in)$.  This leads
to forcing notions which are ``reasonably" definable.  We present two specific
properties materializing this intuition: nep (non-elementary properness) and
snep (Souslin non-elementary properness).  For this we consider candidates 
(countable models to which the definition applies), and the older Souslin 
proper. A major theme here is ``preservation by iteration'', but we also show
a dichotomy: if such forcing notions preserve the positiveness of the set of
old reals for some naturally define c.c.c.\ ideals, then they preserve the
positiveness of any old positive set. We also prove that (among such forcing
notions) the only one commuting with Cohen is Cohen itself. 
\bigskip
\bigskip
\bigskip
\setcounter{section}{-1}

\stepcounter{section}
\subsection*{Annotated Content}\ 

\noindent {\bf Section 0:\quad Introduction}\quad
We present reasons for developing the theory of forcing notions which satisfy
the properness demand for countable models which are not necessarily
elementary submodels of some $({\cal H}(\chi),\in)$.  This will lead us to
forcing notions which are ``reasonably" definable.  
\medskip

\noindent {\bf Section 1:\quad Basic definitions}\quad
We present two specific properties materializing this intuition: nep 
(non-elementary properness) and snep (Souslin non-elementary properness). For
this we consider candidates (countable models to which the definition
applies), and the older Souslin proper.
\medskip

\noindent {\bf Section 2:\quad Connections between the basic definitions}\quad
We point out various implications (snep implies nep, etc.). We also point out
how much the properties are absolute.
\medskip

\noindent {\bf Section 3:\quad There are examples}\quad
We point out that not just the reasonably definable forcing notions in use fit
our framework, but that all the general theorems of Ros{\l}anowski Shelah
\cite{RoSh:470}, which prove properness, actually prove the stronger
properties introduced earlier.
\medskip

\noindent {\bf Section 4:\quad Preservation under iteration: first round}\quad
First we address a point we ignored earlier (it was not needed, but is
certainly part of our expectations).  In the definition of ``$q$ is
$(N,\bQ)$-generic" predensity of each ${\cal I}\in\pd(N,\bQ)$ was originally
designed to enable us to say things on $N[\name{G}_{\bQ}]$, i.e.\ $N[G_{\bQ}]
\cap{\cal H}(\chi)^{\V}=N$, but we should be careful saying what we intend by
$N[G_{\bQ}]$ now, so we replace it by $N\langle\name{G}_{\bQ}\rangle$. The
preservation theorem \ref{2.5} says that CS iterations of nep forcing notions
have the main property of nep.  For this we define $p^{\langle\langle N\rangle
\rangle}$ if $N \models$`` $p \in\Lim(\bar\bQ)$ ''.  We also define and should
consider (\ref{2.3A}) the ``$K$-absolute nep". 
\medskip

\noindent {\bf Section 5:\quad True preservation theorems}\quad 
We consider two closure operations of nep forcing notions $(\cl_1,\cl_2)$,
investigate what is preserved and what is gained and prove a general
preservation theorem (\ref{4.9}). This is done for the ``straight" version of
nep. 
\medskip

\noindent {\bf Section 6:\quad When a real is $(\bQ,\name{\eta})$--generic
over $\V$}\quad
We define the class $\cK$ of pairs $(\bQ,\name{\eta})$, in particular when
$\name{\eta}$ is the generic real for $\bQ$, and how nice is the subforcing
$\bQ'$ of $\bQ$ generated by $\name{\eta}$. 
\medskip

\noindent {\bf Section 7:\quad Preserving a little implies preserving much}
\quad  
We are interested in the preservation of the property (of forcing notions)
``retaining positiveness modulo the ideal derived from a c.c.c. nep forcing
notion'', e.g.\ being non-null (by forcing notions which are not necessarily
c.c.c.). In \cite[Ch.VI,\S1,\S2,\S3, Ch.XVIII,\S3]{Sh:f} this is dealt with
but mainly in the limit case. Our main aim is to show that for ``nice'' enough
forcing notion we have a dichotomy (which implies preservation under e.g.\ CS
iterations (of proper forcing) of the property above) retaining the
positiveness of $\baire$ (or in general every positive Borel set) implies
retaining the positiveness of any $X \subseteq \baire$.  
\medskip

\noindent {\bf Section 8:\quad Non-symmetry}
\quad  
We start to investigate for c.c.c. nep forcing: when does ``if $\eta_0$ is
$(\bQ_0,\name{\eta}_0)$-generic over $N$ and $\eta_1$ is $(\bQ_1,
\name{\eta}_1)$-generic over $N[\eta_0]$ then $\name{\eta}_1$ is $(\nbQ_0,
\name{\eta}_0)$-generic over $N[\eta_1]$''?  This property is known for Cohen
reals and random reals above.
\medskip

\noindent {\bf Section 9:\quad Poor Cohen commute only with himself}
\quad  
We prove that commuting with Cohen is quite rare.  In fact, c.c.c. Souslin
forcing which adds $\name{\eta}$ which is (absolutely) nowhere essentially
Cohen does not commute with Cohen. So such forcing makes the set of old reals
meagre. 
\medskip

\noindent {\bf Section 10:\quad Some c.c.c.\ nep forcing notions are not nice} 
\quad  
We define such forcing notions which are not essentially Cohen as long as
$\aleph_1$ is not too large in ${\bf L}$. This shows that ``c.c.c.  Souslin''
cannot be outright replaced by ``absolutely c.c.c. nep''.
\medskip

\noindent {\bf Section 11:\quad Preservation of ``no dominating real''}  
\quad  
We would like to strengthen the main conclusion of \S7, (that retaining of
positiveness is preserved by composition) of nice forcing notions
(i.e. if each separately has it, then so does its composition) to additional
natural ideals, mainly the one mentioned in the title, which does not flatly
fall into the context of \S7. Though \ref{9.2} contains a counterexample, 
we prove it for ``nice'' enough forcing notions.
\medskip

\noindent {\bf Section 12:\quad Open problems}  
\quad  
We formulate several open questions.

\stepcounter{section}
\subsection*{\quad 0. Introduction} 

The thema of \cite{Sh:b}, \cite{Sh:f} is:
\begin{thesis}
\label{0.1}
It is good to have general theory of forcings, particularly for
iterated forcing. 
\end{thesis}

Some years ago, Judah asked me a question (on inequalities on cardinal
invariants of the continuum). Looking for a forcing proof we arrived to the
following question: 
\begin{question}
Will it not be nice to have a theory of forcing notions $Q$ such that:
\begin{enumerate}
\item[$(\oplus)$] {\em if} $\bQ\in N\subseteq ({\cal H}(\chi),\in)$, $N$
a countable model of $\ZFC^-$ and $p\in N\cap\bQ$,

{\em then} there is $q\in\bQ$ which is $(N,\bQ)$-generic. 
\end{enumerate}
\end{question}
Note the absence of $\prec$ (i.e.\ $N$ is just a submodel of $({\cal H}(\chi),
\in)$), which is the difference between this property and ``properness", and is
alluded to in the name of this paper. This evolved to ``Souslin proper
forcing'' (see \ref{0.7}) in Judah and Shelah \cite{JdSh:292}, which was
continued in Goldstern Judah \cite{GoJu}.

There are still some additional desirable properties (absent there):
\begin{enumerate}
\item[(a)]  many ``nicely defined'' forcing notions do not satisfy ``Souslin
proper'', in fact not so esoteric ones: the Sacks focing, the Laver forcing;
\item[(b)]  actual preservation by CS iteration was not proved, just the
desired conclusion $(\oplus)$ hold for $\bP_\alpha$ when $\langle \bP_i, 
\nbQ_j:i \leq\alpha,j<\alpha\rangle$ is a countable support iteration and $i<
\alpha\quad \Rightarrow\quad \forces_{\bP_i}$`` $\nbQ_i$ is a Souslin proper
forcing notion'';
\item[(c)]  to prove for such forcing notions better preservation theorems
when we add properties in addition to properness.
\end{enumerate}
Martin Goldstern asked me some years ago on the inadequacy of Souslin proper
from clause (a). I suggested a version of the definitions here, and this was
preliminarily announced in Goldstern \cite{Go}.

The intention here is to include forcing notions with ``nice definition'' (not
ones constructed by diagonalization like Baumgartner's ``every
$\aleph_1$-dense sets of reals are isomorphic" \cite{B4} or the forcing
notions constructed for the oracle c.c.c., see \cite[Ch.IV]{Sh:b}, or forcing
notions defined from an ultrafilter).

Note that our treatment (nep/snep) in a sense stands between \cite{Sh:f} and
Ros{\l}anowski Shelah \cite{RoSh:470}. In \cite{Sh:f} we like to have theorems
on iterations $\bar{\bQ}$, mainly CS, getting results on the whole
$\Lim(\bar{\bQ})$ from assumptions on each $\bQ_i$, but with no closer look at
$\bQ_i$ -- by intention, as we would like to cover as much as we can. In
Ros{\l}anowski Shelah \cite{RoSh:470} we deal with forcing notions which are
quite concrete, usually built from countably many finite ``creatures'' (still
relative to specific forcing this is quite general). 

Here, our forcing notions are definable but not in so specific way as in
\cite{RoSh:470}, which still provides examples (all are included), and the
theorems are quite parallel to \cite{Sh:f}. So we are solving the ``equations''
\medskip

$x/$``theory of manifolds'' = 

theory of proper forcing \cite{Sh:b},\cite{Sh:f}/general topology = 

theory of forcing based on creatures \cite{RoSh:470}/theory of manifolds in
${\Bbb R}^3$. 

\begin{thesis}
\label{0.1B}
``Nice'' forcing notions which are proved to be proper, normally satisfy (even
by same proof) the stronger demands defined in the next section. 
\end{thesis}

\noindent{\bf History:}\qquad The paper is based on the author's lectures in
Rutgers University in Fall 1996, which results probably in too many
explanations. Answering Goldstern's question was mentioned above. A version of
\S8 (on non-symmetry) was done in Spring of '95 aiming at the symmetry
question, and the rest in the Summer and Fall of '96. I thank the audience of
the lectures for their remarks and mainly Andrzej Ros{\l}anowski for
correcting the paper. 
\bigskip

\noindent{\bf Notation:}\qquad We try to keep our notation standard and
compatible with that of classical textbooks on Set Theory (like Bartoszy\'nski
Judah \cite{BaJu95} or Jech \cite{J}). However in forcing we keep the
older tradition that {\em a stronger condition is the larger one}. 

For a regular cardinal $\chi$, ${\cal H}(\chi)$ stands for the family of sets
which are hereditarily of size less than $\chi$. The collection of all sets
which are hereditarily countable relatively to $\kappa$ is denoted by ${\cal
H}_{{<}\aleph_1}(\kappa)$.\\
${\rm Tc}^{\rm ord}(x)$ is defined by induction on $\rk(x)=\gamma$ as follows: 

if $\gamma =0$ then  ${\rm Tc}^{\rm ord}(x)=x\cup \{x\}$,

if $\gamma > 0$ then ${\rm Tc}^{\rm ord}(x)=x \cup \bigcup \{{\rm Tc}^{\rm
ord}(y):y \in x,y \mbox{ not an ordinal}\} \cup \{x\}$.\\
So ${\cal H}_{{<}\aleph_1}(\kappa)=\{x\in{\cal H}(\kappa):{\rm Tc}^{\rm ord}
(x)$ is countable$\/\}$.\\
We say that a set $M\subseteq{\cal H}(\chi)$ is ${\rm ord}$--transitive if 
\[x\in M\ \&\ x\mbox{ is not an ordinal}\quad\Rightarrow\quad x\subseteq M.\]

\begin{notation}
We will keep the following rules for our notation:
\begin{enumerate}
\item $\alpha,\beta,\gamma,\delta,\xi,\zeta, i,j\ldots$ will denote ordinals,
\item $\theta,\kappa,\lambda,\mu,\chi\ldots$ will stand for cardinal numbers,
$\theta \le \kappa$ if not said otherwise,
\item a tilde indicates that we are dealing with a name for an object in
forcing extension (like $\name{x}$),
\item a bar above a name indicates that the object is a sequence, usually
$\bar{X}$ will be $\langle X_i: i<\lh(\bar{X})\rangle$, where $\lh(\bar{X})$
denotes the length of $\bar{X}$,
\item For two sequences $\eta,\nu$ we write $\nu\vartriangleleft\eta$ whenever
$\nu$ is a proper initial segment of $\eta$, and $\nu\trianglelefteq\eta$ when
either $\nu\vartriangleleft\eta$ or $\nu=\eta$. The length of a sequence
$\eta$ is denoted by $\lh(\eta)$.
\item A {\em tree} is a family of finite sequences closed under initial
segments. For a tree $T$ the family of all $\omega$--branches through $T$ is
denoted by $\lim(T)$. 
\item The Cantor space $\can$ and the Baire space $\baire$ are the spaces of
all functions from $\omega$ to $2$, $\omega$, respectively, equipped with
natural (Polish) topology. 
\item The fix ``version'' $\ZFC^-_*$ should be such that the forcing theorem
holds and for any large enough $\chi$, the set of $({\frak B},\bar{\varphi},
\theta)$--candidates (defined in \ref{0.1C}) is cofinal in $\{N:N\subseteq
({\cal H}(\chi),\in)\}$ and whatever we should use (fully see \ref{0.9}). 
\item $\frak C$, $\frak B\ldots$ will denote models (with some countable
vocabulary). For a model $\frak C$, its universe is denoted $|{\frak C}|$ and
its cardinality is $\|{\frak C}\|$. Usually $\frak C$'s universe is an ordinal
$\alpha({\frak C})$ and $\kappa({\frak B})\subseteq |{\frak B}|\subseteq {\cal
H}_{{<}\aleph_1}(\kappa({\frak B}))$, $\kappa({\frak B})$ a cardinal.
\item $K$ will denote a family of forcing notions including the trivial one
(so a $K$--forcing extension of $\V$ is $\V[G]$ when $G\subseteq \bP\in K$ is
generic over $\V$) and $\forces_{\bP}$`` $\nbQ\in K^{\V^{\bP}}$ ''$\quad
\Rightarrow \quad \bP*\nbQ\in K$. Usually $K$ is the class of (set) forcing
notions. 
\end{enumerate}
\end{notation}

\stepcounter{section}
\subsection*{\quad 1. Basic definitions}
Let us try to analyze the situation. Our intuition is that: looking at $\bQ$
inside $N$ we can construct a generic condition $q$ for $N$, but if $N\nprec
({\cal H}(\chi),\in)$, $\bQ\cap N$ might be arbitrary.  So let $\bQ$ be a
definition. What is the meaning of, say, $N\models$``$r\in \bQ$''?  It is 
$N\models$``$r$ satisfies $\varphi_0(-)$" for a suitable $\varphi_0$. It seems
quite compelling to demand that inside $N$ we can say in some sense ``$r\in
\bQ$'', and as we would like to have
\[q\forces\mbox{`` }\name{G}_{\bQ}\cap\bQ^N \mbox{ is a subset of }\bQ^N
\mbox{ generic over $N$ ''},\]
we should demand
\begin{enumerate}
\item[$(*)_1$] $N\models$`` $r\in \bQ$ '' implies $\V\models$`` $r\in \bQ$ ''.
\end{enumerate}
So $\varphi_0$ (the definition of the set of members of $\bQ$) should have
this amount of absoluteness.  Similarly we would like to have: 
\begin{enumerate}
\item[$(*)_2$] if $N\models$`` $p_1\le_{\bQ} p_2$ '' and $p_2\in G_{\bQ}$ then
$p_1 \in G_{\bQ}$. 
\end{enumerate}
So we would like to have a $\varphi_1$ (or $<^{\varphi_1}$) (the definition of
the partial order of $\bQ$) and to have this upward absoluteness for
$\varphi_1$. 

But before we define this notion of properness without elementaricity, we
should define the class of models to which it applies.

We may have put in this section the ``straight nep'' (see \ref{4.8}) and/or
``absolute nep'' (see \ref{2.3A}).  {\em Advice:\ } The reader may concentrate
on the case of local correct explicit simply good and nep forcing notions
which are normal (see Definitions \ref{0.1C}, \ref{0.2}(11), \ref{0.2}(2),
\ref{0.2}(5), \ref{0.2}(1),(4), \ref{0.9}(3), \ref{0.9}(4), respectively). 

When we consider ``preservation by iteration'', it is natural to define the
following: 

\begin{definition}
\label{0.1C}
\begin{enumerate}
\item Let $\ZFC^-_*$ be a fixed version of set theory e.g.\ $\ZFC^- +
``\beth_7$ exists'' which may speak on $\frak C$ (or see more in the end of
this section), and let ${\frak C}$ be a fixed model with countable vocabulary
(say $\subseteq {\cal H}(\aleph_0)$) and universe an ordinal $\alpha =
\alpha_*({\frak C})$ and let $\Delta$ be a fixed set of first order formulas
in the vocabulary of ${\frak C}$ (closed under subformulas normally). Let
${\frak B}$ denote another such model (not fixed) but we may allow the
universe to satisfy $\kappa({\frak B}) \subseteq |{\frak B}| \subseteq {\cal
H}_{<\aleph_1}(\kappa({\frak B}))$ for some cardinal\footnote{We do not fix
the order between $\alpha_*({\frak C})$ and $\kappa({\frak B})$, but there is
no loss if we assume that $\theta\geq\alpha_*({\frak C})$, $\theta\geq\kappa({
\frak B})$.}  $\kappa({\frak B})$. 
\item We say that $N$ is a class $({\frak B},\bp,\theta)$--candidate if:
\begin{enumerate}
\item[(a)] $N\subseteq ({\cal H}(\chi),\in)$ for some $\chi$,
\item[(b)] $N$ is countable,
\item[(c)] $N$ is a model of $\ZFC^-_*$,
\item[(d)] ${\frak C}\in N$, $\bp\in N$, ${\frak B} \in N$ (but see below),
\item[(e)] ${\frak B}\rest N \prec_\Delta {\frak B}$ {\em but}\footnote{so if
$\Delta = \Delta_1=\{\exists y \psi:\psi \mbox{ is q.f.}\}$, ${\frak C},{\frak
B}$ have Skolem functions, we have\\ 
{\bf (e)'} \quad ${\frak C} \rest N \prec {\frak C}$, ${\frak B} \rest N \prec
{\frak B}$;\\
can use {\bf (e)'} instead of {\bf (e)}} for transparency we treat ${\frak B}$
as relations of $N$ and ${\frak B}\rest N$ are their interpretations in $N$;
so we allow $|{\frak B}|\cap |N|\setminus |{\frak B} \rest N|\neq\emptyset$,
and $\tau({\frak B})$, the vocabulary of $\frak B$, belongs to $N$, but
$N\models$``$x\in {\frak B}$'' $\Rightarrow\ x\in {\frak B}$. Similarly for
${\frak C}$ (this is less essential), (see \ref{1.10}(3)),
\item[(f)] if $N\models$``$\alpha$ is an ordinal $<\theta$, or $\le |{\frak
C}|$ (that is $\alpha_*({\frak C})$ or $\le \kappa({\frak B})$'' then $\alpha$
is an ordinal, 
\item[(g)] if $N \models$``$x$ is countable'' then $x \subseteq N$,
\item[(h)] if $N\models$``$x$ is an ordinal'' then $x$ is an ordinal.
\end{enumerate}
\item We omit the ``class'' if additionally
\begin{enumerate}
\item[(i)] $\bp =\bar{\varphi}$ is a tuple of formulas, $\varphi_0=\varphi_0
(x)$ and in $N$, $\varphi_0(x)$ defines a set\footnote{This is normally the
forcing notion $\bQ$.}. 
\end{enumerate}
\item We add the adjective ``semi'' if we omit clause (b) (the countability
demand).
\item If $\bp$ is absent (or clear from context) we may omit it, similarly
$\theta$ when $\theta = \aleph_0$ or clear from the context. We tend to
``forget" to mention ${\frak C}$ (e.g.\ demand ${\frak B}$ expands it).
\item We say that a formula $\varphi$ is upward absolute for (or from) class
$({\frak B},\bold p,\theta)$--candidates when: if $N_1$ is a class $({\frak B},
\bp,\theta)$--candidate, $N_1 \models \varphi[\bar x]$, and $N_2$ is a class
$({\frak B},\bp,\theta)$--candidate or is $({\cal H}(\chi),\in)$ for $\chi$
large enough, and $N_1$ is a set or just a class of $N_2$, {\em then} $N_2
\models \varphi[\bar x]$.

We say above ``through (class) $({\frak B},\bold p,\theta)$--candidates'' if
$N_2$ is demanded to be a (class) $({\frak B},\bold p,\theta)$--candidate. 
Note that we can omit ${\cal H}(\chi)$ in the correct case (see
\ref{0.2}(11)). 

If ${\frak B},\bp,\theta$ are clear from the context, we may forget to say
``for class $({\frak B},\bp,\theta)$--candidates''.
\item We say that $\varphi$ defines $X$ absolutely through $({\frak B},\bold
p,\theta)$--candidates if
\begin{enumerate}
\item[$(\alpha)$] $\varphi=\varphi(x)$ is upward absolute through $({\frak B},
\bold p,\theta)$--candidates,
\item[$(\beta)$]  $X=\bigcup\{X^N: N$ is a $({\frak B},\bold p,\theta
)$--candidate$\/\}$,\quad where $X^N=\{x\in N:N\models \varphi(x)\}$.
\end{enumerate}
If only clause $(\alpha)$ holds then we add ``weakly''.
\end{enumerate}
\end{definition}

\begin{discussion}
\label{1.1A}
Should we prefer $|{\frak B}| = \alpha$ an ordinal or $|{\frak B}|\subseteq
{\cal H}_{< \aleph_1}(\alpha)$?  The former is more convenient when we
``collapse $N$ over $\kappa \cup \theta$'' (see \ref{1.10}). Also then we can
fix the universe whereas $|{\frak B}| = {\cal H}_{< \aleph_1}(\alpha)$ is less
reasonable as it is less absolute. On the other hand, when we would like to
prove preservation by iteration the second is more useful (see \S5). To have
the best of both we adopt the somewhat unnatural meaning of ${\frak B} \rest N
\prec_\Delta {\frak B}$ in clause (e) of Definition \ref{0.1C}.\\
We may have forgotten sometimes to write $\|{\frak B}\|$ instead of $\kappa =
\kappa({\frak B})$. 

In some cases, we may omit the demand (h) in the definition \ref{0.1C} of
$({\frak B},\bp,\theta)$--candidates (and then calling them ``impolite
candidates''), but still we should demand then that 
\[N\models\mbox{`` $x$ is an ordinal from $\frak B$ or $\frak C$''}\quad
\Rightarrow\quad\mbox{$x$ is an ordinal,}\]
and we should change ``ordinal collapse'' appropriately. However, there is no
reason to attend the ``impolite'' company here.
\end{discussion}

This motivates:

\begin{definition}
\label{0.2}
\begin{enumerate}
\item Let $\bar{\varphi}=\langle\varphi_0,\varphi_1\rangle$ and ${\frak B}$ be
a model as in \ref{0.1C}, $\kappa=\kappa({\frak B})$, and of countable
vocabulary, say $\subseteq {\cal H}(\aleph_0)$. We say that $\bar{\varphi}$ or
$(\bar{\varphi},{\frak B})$ is a temporary $(\kappa,\theta)$--definition, or
$({\frak B},\theta)$--definition, of a nep-forcing notion\footnote{so in the
normal case (see \ref{0.9}(4), \ref{0.11}), $\bar{\varphi}$ defines $\bQ$}
$\bQ$ if, in $\V$: 
\begin{enumerate}
\item[(a)] $\varphi_0$ defines the set of elements of $\bQ$ and $\varphi_0$ is
upward absolute from $({\frak B},\bar{\varphi},\theta)$--candidates,
\item[(b)] $\varphi_1$ defines the partial (or quasi) ordering of $\bQ$, also
in every $({\frak B},\bar{\varphi},\theta)$--candidate, and $\varphi_1$ is
upward absolute from $({\frak B},\bar{\varphi},\theta)$--candidates, 
\item[(c)] if $N$ is a $({\frak B},\bar{\varphi},\theta)$--candidate and $p
\in \bQ^N$, {\em then} there is $q \in \bQ$ such that $p \le^{\bQ} q$ and
\[q \forces\mbox{`` } \name{G}_{\bQ}\cap \bQ^N \mbox{ is a subset of }\bQ^N
\mbox{ generic over } N\mbox{ ''}\]
where, of course, $\bQ^N = \{p:N \models \varphi_0(p)\}$.
\end{enumerate}
\item We add the adjective ``explicitly'' if $\bar{\varphi}=\langle\varphi_0,
\varphi_1,\varphi_2\rangle$ and additionally
\begin{enumerate}
\item[(b)$^+$]  we add: $\varphi_2$ is an $(\omega+1)$-place relation, upward
absolute through $({\frak B},\bar{\varphi},\theta)$--candidates and $\varphi_2(
\langle p_i:i \leq\omega\rangle)\ \Rightarrow\ \mbox{``}\{p_i:i\le\omega\}
\subseteq \bQ$ and $\{p_i:i<\omega\}$ is predense above $p_\omega$'', not just
in $\V$ but in every $\bQ$--candidate (which, if $\bQ$ is correct, implies the
case in $\V$); in this situation we say: $\{p_i:i<\omega\}$ is explicitly
predense above $p_\omega$, 
\item[(c)$^+$]  we add: if $N \models\mbox{``}{\cal I} \subseteq\bQ$ is dense 
open'' (or just predense) (so ${\cal I} \in N$) {\em then} for some list
$\langle p_i:i<\omega\rangle$ of ${\cal I} \cap N$ we have $\varphi_2(\langle
p_i:i < \omega\rangle\conc\langle q \rangle)$. 
\end{enumerate}
\item For a class $({\frak B},\bar{\varphi},\theta)$--candidate $N$ we let 
$\pd(N,\bQ)=\pd_{\bQ}(N)=\{{\cal I}:{\cal I}$ is a class of $N$ (i.e. defined
in $N$ by a first order formula with parameters from $N$) and is a predense
subset of $\bQ^N\}$. If $N$ is a candidate, it is $\{ {\cal I} \in N:N
\models$``${\cal I}$ is predense''$\}$.
\item  We replace ``temporary'' by $K$ if the relevant proposition holds not
only in $\V$ but in any forcing extension of $\V$ by a forcing notion $\bP\in
K$. If $K$ is understood from the context (normally: all forcing notions we
will use in that application) we may omit it.
\item We say that $(\bar{\varphi},{\frak B})$ is simply [explicitly]
$K$--$(\kappa,\theta)$--definition of a nep--forcing notion $\bQ$, {\em if}:
\begin{enumerate}
\item[$(\alpha)$] $(\bar{\varphi},{\frak B})$ is [explicitly] $K$--definition
of a nep-forcing notion $\bQ$,
\item[$(\beta)$]  $\bQ\subseteq {\cal H}_{<_{\aleph_1}}(\theta)$; i.e.\ $\bP
\in K$ implies $\forces_{\bP}$ ``if $\varphi_0(x)$ then $x \in {\cal H}_{<_{
\aleph_1}}(\theta)$'',
\item[$(\gamma)$] ${\frak B},\kappa,\theta$ are the only parameters of
$\bar{\varphi}$ (meaning there are no others, but even ${\frak B},\kappa,
\theta$ do not necessarily appear). 
\end{enumerate}
\item We add ``very simply'' {\em if} in addition:
\begin{enumerate}
\item[$(\delta)$] $\bQ \subseteq {}^\omega\theta$.
\end{enumerate}
\item We may say ``$\bQ$ is a nep-forcing notion", ``$N$ is a
$\bQ$-candidate'' abusing notion. If not clear, we write $\bQ^{\bar{\varphi}}$
or $(\bQ^{\bar{\varphi}})^{\V}$. If not said otherwise, $\Delta$ is the set of
first order formulas. Inversely, we write $({\frak B},\bar{\varphi},\theta)=
({\frak B}^{\bQ},\bar{\varphi}^{\bQ},\theta^{\bQ})$ and $\ZFC^{\bQ}$ for the
relevant $\ZFC^-_*$. 
\item We say ``${\cal I} \subseteq \bQ^N$ is explicitly predense over
$p_\omega$'' if $\varphi_2(\langle p_i:i\le\omega\rangle)$ for some list
$\{p_i:i<\omega\}$ of ${\cal I}$.
\item We add the adjective ``class'' if we allow ourselves (in clauses (b),
(c) of part (1) and (c)$^+$ of part (2)) class $({\frak B},\bar{\varphi},
\theta)$--candidates $N$; so in clauses (c), (c)$^+$, ${\cal I}$ is a class of
$N$; i.e.\ first order definable with parameters from $N$, and use the weak
version of absoluteness. 

If we use $({\frak B},\bp,\theta)$ we mean $\bar{\varphi}$ is an initial
segment of $\bp$. 
\item We say $({\frak B},\bar{\varphi},\theta)$ (or abusing notation, $\bQ$) is
class=set if every class $({\frak B},\bar{\varphi},\theta)$--candidate is a
$({\frak B},\bar{\varphi},\theta)$--candidate.
\item In \ref{0.2}(1) we add the adjective ``correctly'' (and we say that
$({\frak B},\bar{\varphi},\theta)$ is {\em correct}) if, for a large enough
regular cardinal $\chi$:
\begin{enumerate}
\item[(a)] the formula $\varphi_0$ defines the set of members of $\bQ$
absolutely through $({\frak B},\bar{\varphi},\theta)$--candidates, that is
\[\bQ=\bigcup\{\bQ^N:N\mbox{ is a }({\frak B},\bar{\varphi},\theta)
\mbox{--candidate }\},\]
$\bQ^N=\{x: N\models\varphi_0(x)\}$,
\item[(b)] the formula $\varphi_1$ defines the quasi order of $\bQ$ absolutely
through $({\frak B},\bar{\varphi},\theta)$--candidates, that is $\leq_{\bQ}=
\bigcup\{(\leq_\bQ)^N: N$ is a $({\frak B},\bar{\varphi},\theta
)$--candidate$\/\}$, $(\leq_\bQ)^N=\{(p,q): N \models\varphi_1(p,q)\}$.
\end{enumerate}
Similarly when we add ``explicitly''.

So in those cases we can ignore ${\cal H}(\chi)\models\varphi_\ell(x)$ and
just ask for satisfaction in suitable candidates. (Note: correct is less
relevant to snep.) 
\end{enumerate}
\end{definition}

\noindent{\bf Convention:}\quad We may say ``$\bQ$ is $\dots$'' when we mean
``$({\frak B},\bar{\phi},\theta)$ is $\ldots$'' or ``$({\frak B},\bp,\theta)$
is $\ldots$''. 
\medskip

\noindent{\bf Remark:}\quad The main case for us is candidates (not class
ones), etc; still mostly we can use the class version of nep. Also we can play
with various free choices.

\begin{discussion}
\label{0.3}
1)\quad Note: if $x\in{\cal I}\in N$, $N \models$`` ${\cal I}\subseteq\bQ$'',
possibly $x\notin\bQ$ so those $x$ are not relevant (e.g.\ though $\alpha<
\kappa({\frak B})$ have a special role). 

\noindent 2)\quad We think of using CS iteration $\bar{\bQ}=\langle\bP_i,
\nbQ_i:i<\delta\rangle$, each $\nbQ_i$ has a definition $\bar{\varphi}^i$ and
we would like to prove things on $\bP_\alpha$ for $\alpha\le\delta$. So the
relevant family $K_i$ of forcing notions we really should consider for
$\bar{\varphi}^i$ is $\{\bP_\beta/\bP_i:\beta\in [i,\delta)\}$, at least this
holds almost always (maybe we can look as help in other extensions). 

\noindent 3)\quad Note that a significant fraction of iterated forcing of
proper forcing related to reals are forcing notions called ``nice'' above. The
proof that they are proper usually gives more and we think that they will be
included even by the same proof. 

\noindent 4)\quad If $K$ is trivial, (i.e.\ has only the trivial forcing
notion as a member) this means we can replace it by ``temporarily''. 

\noindent 5)\quad See also \ref{2.3A} for ``$K$-absolutely". 

\noindent 6)\quad Note a crucial point in Definition \ref{0.2}, the relation
``$\{p_n:n < \omega\}$ is predense above $p$'' is not demanded to be absolute;
only a ``dense'' family of cases of it is demanded (we also allow other basic
relations; e.g.\ $q\notin\bQ$ to be non-absolute but those are less crucial).
This change may seem technical, but is central being the difference between
including not few natural examples and including all those we have in mind.

\noindent 7)\quad Note that in clause (c) of \ref{0.2}(1) we mean: \quad
$G\cap \bQ^N$ is directed (by $\leq^{\bQ^N}$, not only by $\leq^{\bQ}$) and
$G\cap N\cap {\cal I}\neq\emptyset$ for ${\cal I}\in \pd(N,\bQ)$.

\noindent 8)\quad Note that the demand described in 7) above almost implies
``incompatibility is upward absolute from $N$'', but not quite.
\end{discussion}

Let us consider a more restrictive class, where the absoluteness holds because
of more concrete reasons, the usual ones for upward absoluteness, the
relations are $\Sigma^1_1$, or more generally, $\kappa$--Souslin.

\begin{definition}
\label{0.4}
\begin{enumerate}
\item We say that $\bar T$ is a temporary $(\kappa,\theta)$--definition of a
snep--forcing notion $\bQ$ if: 
\begin{enumerate}
\item[(a)] $\bar{T}=\langle T_0,T_1\rangle$ where $T_0 \subseteq {}^{\omega
>}(\kappa\times\kappa)$ and $T_1\subseteq {}^{\omega >}(\kappa\times\kappa
\times \kappa)$ are trees (i.e.\ closed under initial segments, non-empty) and
$\theta\le\kappa$,
\item[(b)] the set of elements of $\bQ$ is 
\[\begin{array}{ll}
\proj_0(T_0)\stackrel{\rm def}{=}\{\nu\in {}^\omega\theta:&\mbox{for some }
\eta \in {}^\omega \kappa \mbox{ we have}\\
\ &\nu * \eta \stackrel{\rm def}{=} \langle(\nu(n),\eta(n)):n <\omega\rangle\in
\lim(T_0)\}, 
  \end{array}\]
\item[(c)] the partial order of $\bQ$, $\{(p_0,p_1):\bQ\models p_0 \le p_1\}$
is 
\[\begin{array}{ll}
\proj_1(T_1)\stackrel{\rm def}{=}&\{(\nu_0,\nu_1):\nu_0,\nu_1\in \bQ \mbox{
and for some }\eta\eta\in{}^{\textstyle\omega}\kappa \mbox{ we have}\\
\ &\quad\nu_0 * \nu_1 * \eta \stackrel{\rm def}{=}\langle (\nu_0(n),\nu_1(n),
\eta(n)):n<\omega\rangle\in\lim(T_1)\},
  \end{array}\]
\item[(d)] for a large enough regular cardinal $\chi$, if $N\subseteq ({\cal
H}(\chi),\in)$ is a $({\frak B}_{\bar{T}},\bar{T},\theta)$--candidate and
$\kappa\in N$, $\bar{T}\in N$, $p\in \bQ^N$ {\em then} there is $q\in \bQ$
such that $p \le^{\bQ} q$ and 
\[q\forces\mbox{`` }\name{G}_{\bQ}\cap \bQ^N \mbox{ is a generic subset of }
\bQ^N \mbox{ over } N\mbox{ ''},\]
where ${\frak B}_{\bar{T}}$ is the model with universe $\kappa$ and the
sequence relations $T_0 \cap {}^n(\kappa\times\kappa)$, $T_1\cap {}^n(\kappa
\times\kappa\times\kappa)$ for $n<\omega$.
\end{enumerate}
\item We add ``explicitly'' if $\bar{T}=\langle T_0,T_1,T_2\rangle$ and we add
\begin{enumerate}
\item[(a)$^+$]  also $T_2 \subseteq {}^{\omega >}(\theta\times\theta\times
\kappa)$ and we let 
\end{enumerate}
\[\begin{array}{ll}
\proj_2(T_2)\stackrel{\rm def}{=}\big\{\langle\nu_i:i \le\omega\rangle:&
\mbox{for some }\eta\in {}^\omega \kappa \mbox{ we have }\nu*\nu_\omega*\eta
\in \lim(T_2)\\
\ &\mbox{where }\nu=\code(\langle\nu_\ell:\ell<\omega\rangle)\mbox{ is the
member}\\
\ &\mbox{of } {}^\omega \theta\mbox{ satisfying } \nu\big(\binom{\ell+k+1}{2}
+ \ell \big)=\nu_\ell(k)\big\}
  \end{array}\]
and $\langle \nu_i:i \le\omega\rangle\in\proj_2(T_2)$ implies $\{\nu_i:i \le
  \omega\} \subseteq \bQ$ (even in candidates; the natural case is that
witnesses are coded). 
\begin{enumerate}
\item[(d)$^+$]  we add: $q$ is $\bar{T}$--explicitly $(N,\bQ)$--generic,
which means that\\
{\em if} $N \models ``{\cal I}$ is a dense open subset of $\bQ$''\\
{\em then} for some list $\langle p_n:n<\omega\rangle$ of ${\cal I}\cap N$ we
have $\langle p_n:n<\omega\rangle\conc\langle q \rangle\in\proj_2 (T_2)$,
\item[(e)$^+$]  if $\nu_i\in\bQ$ for $i\le\omega$ and for some $\eta\in
{}^\omega \kappa$ we have $\code(\nu_0,\nu_1,\ldots)*\nu_\omega*\eta\in
\lim(T_2)$  {\em then} $\{\nu_0,\nu_1,\ldots\} \subseteq \bQ$ is predense
above $\nu_\omega$ (and this holds in candidates too). 
\end{enumerate}
\item We will also say ``$\bQ$ is a snep-forcing notion'', ``$N$ is a
$\bQ$--candidate", etc. 
\item  We say $\eta$ is a witness for $\nu\in\bQ$ if $\nu*\eta\in\lim(T_0)$;
similarly for $T_1,T_2$. We say that ${\cal I}$ is explicitly predense over
$p_\omega$ if $\code(\langle p_i:i \le \omega\rangle)\in\proj_2(T_2)$ for some
list $\{p_i:i<\omega\}$ of ${\cal I}$. 
\end{enumerate}
\end{definition}

\begin{remark}
\label{0.5}
In clause (a)$^+$ we would like the $\proj_2(T_2)$ to be an $(\omega +
1)$-place relation on $\bQ$, but we do not like the first coordinate to give
too much information so we use the above coding, but it is in no way special.
Note: we do not want to have one coordinate giving $\langle \varphi_\ell(0):
\ell<\omega\rangle$.

Another possible coding is $\code(\nu_0,\nu_1,\ldots)\cong\langle\langle
\nu_\ell \rest i:\ell \le i \rangle:i < \omega \rangle$, so $T \subseteq
{}^{\omega>}({}^{\omega>}({}^{\omega>}\theta)\times\theta\times\kappa)$.
\end{remark}

\begin{proposition}
\label{added14A}
Assume that $\bar{T}$ is in $\V$ a temporary $(\kappa,\theta)$--definition of
a snep forcing notion which we call $\bQ$. Let $\V'$ be a transitive class of
$\V$ containing $\bar{T}$. Then:
\begin{enumerate}
\item also in $\V'$, $\bQ$ is snep,
\item if $\V'\models$``$p\in\bQ$'' then $\V\models$``$p\in\bQ$'',
\item if $\V'\models$``$p\leq^{\bQ} q$'' then $\V'\models$``$p\leq^\bQ q$'',
\item if in $\V'$, the model $N$ is a $({\frak B}_{\bar{T}},
\bar{\varphi}_{\bar{T}},\kappa_{\bar{T}})$--candidate then also in $\V$, $N$
is a $({\frak B}_{\bar{T}},\bar{\varphi}_{\bar{T}},
\kappa_{\bar{T}})$--candidate. \QED
\end{enumerate}
\end{proposition}

\begin{definition}
\label{0.5A}
\begin{enumerate}
\item Let $\bQ$ be explicitly snep. We add the adjective ``local'' if in the
``properness clause i.e.\ \ref{0.4}(2)(d)$^+$'' we can add:
\begin{enumerate}
\item[$(\otimes)$] the witnesses for ``$q\in\bQ$'', ``$\langle p^{\cal I}_n:
n<\omega\rangle$ is $\bQ$--explicitly predense above $q$'' are from ${}^\omega
(N\cap\kappa)$. 
\end{enumerate}
\item Let $\bQ$ be explicitly nep. We add the adjective ``$K$-local'' if in the
``properness clause i.e.\ \ref{0.2}(2)(b)$^+$'' we can add:\quad for each
candidate $N$ which is ${\rm ord}$--transitive we have 
\begin{enumerate}
\item[$(\oplus)$] for some $K$--extension $N^+$ of $N$, it is a
$\bQ$--candidate (in particular a model of $\ZFC^-_*$) and $N^+\models$
``$\bQ^N$ is countable'' and $q\in N^+$, $N^+\models$ ``$p\le^{\bQ} q$ and for
each ${\cal I}\in\pd(N,\bQ)$, ${\cal I}^N$ is explicitly predense over $q$''. 
\end{enumerate}
(Note that ${\frak B}\rest N^+={\frak B}\rest N$.)

If $K$ is the family of set forcing notions, or constant understood from the
context, we may omit $K$.
\end{enumerate}
\end{definition}

\begin{discussion}
\label{0.6}
1)\quad Couldn't we fix $\theta=\omega$? Well, if we would like to have the
result of ``the limit of a CS iteration $\bar{\bQ}$ of such forcing notions is
such a forcing notion'', we normally need $\theta\ge\lh(\bar{\bQ})$. Also
$\kappa>\aleph_0$ is good for including $\Pi^1_2$--relations. 

\noindent 2)\quad In ``Souslin proper" (starting with \cite{JdSh:292}) the
demands were 
\end{discussion}

\begin{definition}
\label{0.7}
A forcing notion $\bQ$ is Souslin proper if it is proper and: the relations
``$x \in \bQ$'', ``$x \le^{\bQ} y$'' are $\Sigma^1_1$ and ``the notion of
incompatibility in $\bQ$'' is $\Sigma^1_1$ (where, of course, the
compatibility relation is $\Sigma^1_1$). 
\end{definition}
This makes ``$\{p_n:n < \omega\}$ is predense over $p_\omega$'' a
$\Pi^1_2$--property, hence an $\aleph_1$-Souslin one. So we can get the
``explicitly'' cheaply, {\em however} possibly increasing $\kappa$. Note that
for a Souslin proper forcing notion $\bQ$, also $p \in \bQ^N \Leftrightarrow p
\in \bQ\ \&\ p\in N$ and similarly for $p \le^{\bQ} q$.
\medskip

\centerline {$* \qquad * \qquad *$}
\medskip

If you like to be more pedant on the $\ZFC^-_*$, look at the following
definition.  Normally there is no problem in having $\ZFC^-_*$ as required. 

\begin{definition}
\label{0.9}
\begin{enumerate}
\item We say $\ZFC^-_*$ is a $K$--good version [with parameter ${\frak C}$,
possibly ``for $({\frak B},\bp,\theta)$'' for ${\frak B},\bp,\theta$ as in
\ref{0.2} from the relevant family] if: 
\begin{enumerate}
\item[(a)] it contains ${\rm ZC}^-$; i.e.\ Zermelo set theory without power
set,

[and the axioms may speak on ${\frak C}$]
\item[(b)] ${\frak C}$ is a model with countable vocabulary (given as a well
ordered sequence, so ${\frak C}$ is an individual constant in the theory
$\ZFC^-_*$) and universe $|{\frak C}|$ is an ordinal $\alpha({\frak C})$, 
\item[(c)] for every $\chi$ large enough, if $X \subseteq {\cal H}(\chi)$ is
countable then for some countable $N \subseteq ({\cal H}(\chi),\in)$, $N
\models \ZFC^-_*$, $X\subseteq N$ and 
\[x \in N\ \&\ N\models \mbox{``}|x|=\aleph_0\mbox{''}\quad \Rightarrow\quad 
x \subseteq N\]
and ${\frak C} \in N$ and ${\frak C}\rest (N\cap |{\frak C}|)\prec{\frak C}$
(can be weakened to a submodel or $\prec_{\Delta}$, we do not loose much as we
can expand by Skolem functions); in the ``for $({\frak B},\bp,\theta)$''
version we add ``$N$ is a $({\frak B},\bp,\theta)$--candidate'',
\item[(d)] $\ZFC^-_*$ satisfies the forcing theorem\footnote{for \ref{6.5} we
need:\quad if $\bP,\bQ$ are forcing notions, $\name{G}$ is a $\bP$--name for a
subset of $\bQ$ such that $\forces$`` $\name{G}$ is a generic subset of $\bQ$
'', and $q\in\bQ\ \Rightarrow\ \not\forces_{\bP}$``$q\notin\name{G}$'' then
for some $\bQ$--name $\name{{\Bbb R}}$ of a forcing notion, $\bQ*\name{{\Bbb
R}}$, $\bP$ are equivalent} (see e.g.\ \cite[Ch. I]{Sh:f}) at least for
forcing notions in $K$,  
\item[(e)] those properties are preserved by forcing notions in $K$ (if $\bP
\in K$, $G \subseteq \bP$ generic over $\V[G]$ then $K^{\V[G]}$ will be 
interpreted as $\{\nbQ[G]:\bP*\nbQ\in K\}$).
\end{enumerate}
\item If $K$ is the class of all (set) forcing notions, we may omit it. 
\item We say $\ZFC^-_*$ is normal if for $\chi$ large enough any countable
$N\prec ({\cal H}(\chi),\in)$ to which ${\frak C}$ belongs is O.K. (for clause
(1)(c) above).
\item We say $\ZFC^-_*$ is semi-normal for $({\frak B},\bp,\theta)$ if for
$\chi$ large enough, for any countable $N\prec({\cal H}(\chi),\in)$ (to which
appropriate $\bp,{\frak C},{\frak B},\theta(\in {\cal H}(\chi)$) belong), for
some $\bQ\in N$ such that $N \models$``$\bQ$ is a forcing notion'' we have:
\begin{enumerate}
\item[$(*)$] {\em if} $N'$ is countable $N\subseteq N'\subseteq({\cal H}(
\chi),\in)$, $N' \cap \chi = N \cap \chi$ and 
\[(\forall x)[N' \models\mbox{``}x\mbox{ is countable ''}\quad\Rightarrow\quad
x \subseteq N'],\]
and $N'$ is a generic extension of $N$ for $\bQ^N$ 

{\em then} $N'$ is $({\frak B},\bp,\theta)$--candidate and $\bQ^{N'} \rest N =
\bQ \rest N$, $\varphi^{N'}_2 \rest N = \varphi^N_2 \rest N$.
\end{enumerate}
We say ``$K$--semi-normal'' if we demand $N\models\bQ\in K$.
\item We say $\ZFC^-_*$ is weakly normal for $({\frak B},\bp,\theta)$ if
clause (c) of part (1) holds.
\item In parts (4), (5) we can replace $({\frak B},\bp,\theta)$ by a family of
such triples meaning $N$ is a candidate for all of them. 
\item In parts (4), (5), (6) if $({\frak B},\bp,\theta)=({\frak B}^{\bQ},
\bar{\varphi}^{\bQ},\theta^{\bQ})$ we may replace $({\frak B},\bp,\theta)$
by $\bQ$.
\end{enumerate}
\end{definition}

\begin{discussion}
\label{0.10}
1)\quad What are the points of parameters? E.g.\ we may have $\kappa^*$ an
Erd\"os cardinal, ${\frak C}$ codes every $A \in {\cal H}(\chi)$ for each
$\chi<\kappa^*$, $\ZFC^-_* = \ZFC^- +$ ``$\kappa^*$ is an Erd\"os cardinal
${\frak C}$ as above'', $K =$ the class of forcing notions of cardinality $<
\kappa^*$. Then we have stronger absoluteness results to play with. 

\noindent 2)\quad On the other hand, we may use $\ZFC^-_* = \ZFC^- + (\forall
r \in {}^\omega 2)(r^\#$ exists) + ``$\beth_7$ exists''.  This is a good
version if $\V \models (\forall r \in {}^\omega 2)(r^\#$ exists) so we can
e.g.\ weaken the definition snep (or Souslin-proper or Souslin-c.c.c). 

\noindent 3)\quad What is the point of semi-normal? E.g.\ if we would like
$\ZFC^-_*\vdash{\rm CH}$, whereas in $\V$ the Continuum Hypothesis fails. But
as we have said in the beginning, the normal case is usually enough. 
\end{discussion}

\begin{proposition}
\label{0.11}
\begin{enumerate}
\item Assume $\ZFC^-_*$ is $\{\emptyset\}$--good. Then the clause (c)$^+$ of
\ref{0.2}(2) is equivalent to clause (c) + $(*)$, where
\begin{enumerate}
\item[$(*)$] {\em if} $p\in\bQ$ and ${\cal I}_n$ is predense over $p$ (for
$n<\omega$), each ${\cal I}_n$ is countable,\\
{\em then} for some $q$, $p\leq q\in\bQ$, and for some $p^n_\ell\in{\cal I}_n$
for $n<\omega$, $\ell<\omega$ we have $\varphi_2(\langle p^n_\ell:\ell<\omega
\rangle\conc q)$
\end{enumerate}
(this is an obvious abusing of notation, we mean that this holds in some
candidate). 
\item If $\ZFC^-_*$ is normal for $({\frak B},\bp,\theta)$ {\em then} in
Definition \ref{0.2}(1),(2) there is no difference between ``absolutely
through'' and ``weakly absolutely''. \QED
\end{enumerate}
\end{proposition}

\begin{proposition}
\label{0.12}
\begin{enumerate}
\item Assume $\V_1\subseteq\V_2$ (so $\V_1$ is a transitive class of $\V_2$
containing the ordinals, ${\frak C},{\frak B},\theta,\bp,\in\V_1$). If
$\ZFC^-_*$ is temporarily good {\em then} also in $\V_1$ it is temporarily
good. 
\item If co-$(\kappa+\theta)$-Souslin relations are downward absolute (from
$\V_2$ to $\V_1$) then also inverse holds.
\end{enumerate}
\end{proposition}

\Proof By Shoenfield--Levy absoluteness. \QED

\stepcounter{section}
\subsection*{\quad 2. Connections between the basic definitions}
We first give the most transparent implications: we can omit ``explicitly''
and we can replace snep by nep (this is \ref{1.1}) and the model ${\frak B}$
can be expanded, $\kappa,\theta$ increased, (see \ref{1.2}). Then we note that
if $\kappa \ge \theta + \aleph_1$ and we are in the simple nep case, we can get
from nep to snep because saying ``there is a countable model $N \subseteq 
({\cal H}(\chi),\in)$ such that $\ldots$'' can be expressed as a
$\kappa$--Souslin relation (see \ref{1.3}) and comment on the non-simple
case. Then we discuss how the absoluteness lemmas help us to change the
universe (in \ref{1.5}), to get the case with a class $K$ from the case of
temporarily (\ref{1.6}) and to get explicit case from snep or from Souslin
proper (in \ref{1.7}).  

\begin{proposition}
\label{1.1}
\begin{enumerate}
\item If $(\bar{\varphi},{\frak B})$ is explicitly a $K$--definition of a
nep-forcing notion $\bQ$, {\em then} $\bar{\varphi}\rest 2$ is a
$K$--definition of a nep-forcing notion $\bQ$.
\item If $\bar{T}$ is explicitly a $K$-definition of a snep-forcing notion
$\bQ$, {\em then} $(\bar{T}\rest 2)$ is a $K$--definition of an snep-forcing
notion $\bQ$.
\item If $\bar{T}$ is [explicitly] a $K$--$(\kappa,\theta)$--definition of a 
snep-forcing notion $\bQ$, and ${\frak B}$ any model with universe $\kappa$
coding the $T_\ell$'s and $\varphi_\ell$ is defined as $\proj_\ell(T_\ell)$,
{\em then} $(\bar{\varphi},{\frak B})$ is very simply [explicitly]
$K$--$(\kappa,\theta)$--definition of a nep forcing notion $\bQ$ (and let
${\frak B}={\frak B}_{\bar{T}}$, $\bar{\varphi}=\bar{\varphi}_{\bar{T}}$).  
\end{enumerate}
\end{proposition}

\Proof Read the definitions. \QED

\begin{proposition}
\label{1.2}
\begin{enumerate}
\item If $(\bar \varphi,{\frak B})$ is [explicitly] a $K$--definition of a
nep-forcing notion and ${\frak B}$ is definable in ${\frak B}'$ (and $\Delta$
is $L_{\omega,\omega}$, or change $\Delta$ accordingly to the interpretation),
{\em then} $(\bar{\varphi},{\frak B}')$ is [explicitly] a $K$-definition of a
nep-forcing notion; moreover, if ${\frak B}$ is the only parameter of the
$\varphi_\ell$, we can replace it by ${\frak B}'$ (changing trivially the
$\varphi_\ell$'s).  
\item Similarly we can increase $\kappa$ and $\theta$ and add ``simply'' (to
the assumption and the conclusion); we may also add ``very simply''.
\end{enumerate}
\end{proposition}

\Proof Straight. \QED
\medskip

A converse to \ref{1.1}(1)+(2) is

\begin{proposition}
\label{1.3}
\begin{enumerate}
\item Assume that $\kappa'=\kappa+\theta+\aleph_1+\|{\frak B}\|$ and
\begin{enumerate}
\item[$(\oplus)$]  $(\bar{\varphi},{\frak B})$ is a correct very simple
[explicit] $K$--$(\kappa,\theta)$--definition of a nep forcing notion $\bQ$.
\end{enumerate}
{\em Then} some $\bar{T}$ is an [explicitly]
$K$--$(\kappa',\theta)$--definition of a snep forcing-notion $\bQ$ (the same
$\bQ$).
\item If $\kappa=\theta=\kappa'=\aleph_0$ we get a similar result with the
$\varphi_\ell$ being $\Pi^1_2$-sets. 
\item If in clause $(\oplus)$ of \ref{1.3}(1) we replace very simple by simple 
(so we weaken $\bQ\subseteq {}^\omega \theta$ to $\bQ\subseteq {\cal H}_{<
\aleph_1}(\theta)$), {\em then} part (1) still holds for some $\bQ'$
isomorphic to $\bQ$.
\end{enumerate}
\end{proposition}

\Proof 1)\quad This is, by now, totally straight; still we present the case of
$\varphi_0$ for part (1) for completeness. If in Definition \ref{0.1C}(2),
clause (e) we use $\prec$, let $\langle\psi^1_n(y,x_0,\dotsc,x_{n-1}):n < 
\omega\rangle$ list the first order formulas in the vocabulary of ${\frak B}$
in the variables $\{y,x_\ell:\ell<\omega\}$, (so in $\psi^1_n$ no $x_\ell$,
$\ell\ge n$ appears, but some $x_\ell$, $\ell<n$ may not appear); if we use
$\prec_\Delta$ let it list subformulas of members of $\Delta$. Similarly
$\langle \psi^2_n(y,x_0,\dotsc,x_{n-1}):4 \le n<\omega\rangle$ for the
vocabulary of set theory. Let us define $T_0$ by defining a set of
$\omega$-sequences $Y_0$, and then we will let $T_0=\{\rho\rest n:\rho\in Y_0
\mbox{ and }n<\omega\}$. For $\alpha<\omega_1$ let $\{\beta_{\alpha,\ell}:\ell
<\omega\}$ list $\{\beta:\beta \le \alpha\}$.

Now let $Y_0$ be the set of $\omega$-sequences $\rho \in {}^\omega(\theta
\times\kappa')$ such that for some $({\frak B},\bar{\varphi},\theta
)$--candidate $N \subseteq ({\cal H}(\chi),\in)$ (so ${\frak B},\theta,\kappa$
belong to $N$) and some list $\langle a_n:n < \omega \rangle$ of the member of
$N$ we have:\quad $\rho = \nu * \eta$; i.e.\ $\rho(n) = (\nu(n),\eta(n))$ and
\begin{enumerate}
\item[(i)] $a_0 = {\frak B}$, $a_1 = \theta$, $a_2 = \kappa$, $a_3 = \nu$,
\item[(ii)] $\{n:N \models a_n \in \kappa'\}=\{\eta(8n+1):0<n<\omega\}$,
\item[(iii)] every $\eta(8n+2)$ is a countable ordinal such that:
\[N \models\mbox{`` }\rk(a_n) < \rk(a_m)\mbox{''\quad iff\quad}\eta(8n+2)< \eta(8m+2) <\aleph_1\le\kappa',\]
\item[(iv)] if ${\frak B}\models (\exists y)\psi^1_n(y,a_0,\ldots,a_{n-1})$
then  ${\frak B} \models \psi^1_n[a_{\eta(8(n+1)+3)},a_0,\ldots,a_{n-1}]$,
\item[(v)]  $N\models \varphi_0[\nu]$; i.e. $N \models \varphi_0[a_3]$,
\item[(vi)] $N\models$`` $a_\ell\in a_m$ ''\quad iff\quad $\eta(8(\binom{\ell
+m+1}{2}+ \ell)+4)=0$, 
\item[(vii)] if $n \ge 4$ and $N \models (\exists y)\psi^2_n(y,a_0,\ldots,
a_{n-1})$ 

then $N \models \psi^2_n[a_{\eta(8n+5)},a_0,\ldots,a_{n-1}]$ and $\eta(8n+6)=
1$,
\item[(viii)] if $N \models$``$a_n$ is a countable ordinal'' and $a_k =
\beta_{a_n,\ell}$ 

then $\eta(8(\binom{\ell+n+1}{2}+\ell)+7)=k$. 
\end{enumerate}
Let $T_0 = \{\rho \rest n:\rho \in Y_0,n < \omega\}$.

\begin{claim}
\label{1.3A}
\begin{enumerate}
\item $Y_0$ is a closed subset of ${}^\omega(\theta \times \kappa)$. 
\item $\bQ=\{\nu\in {}^\omega \theta:(\exists \eta)(\eta\in {}^\omega(\kappa')
\ \&\ \nu * \eta\in Y_0 (=\lim(T_0))\}=\proj_0(T_0)$.
\end{enumerate}
\end{claim}

\noindent{\em Proof of the claim:}\qquad 1)\quad Given $\nu*\eta\in\lim(T_0)$
we can define a model $N'$ with set of elements say $\{a'_n:n < \omega\}$ by
clause (vi), it is a model of $\ZFC^-_*$ by clause (vii) (and the demand $N
\models \ZFC^-_*$), it is well founded by clause (iii) (and the earlier
information). 

We start to define an embedding $h$ of $N'$ into ${\cal H}(\chi)$ and we put
$h(a'_0) = {\frak B}$, $h(a'_1) = \theta$, $h(a'_2) = \kappa$ and $h(a'_n) =
\eta(8n+1)$ if $N' \models a'_n \in a'_2$, $n>0$. Then let $h(a'_3) \in
{}^\omega \theta$ be such that $h(a'_3)(\ell)=\gamma$ iff letting $n$ be such
that $\psi^2_n\equiv [y = x_3(\ell)]$, so necessarily $N'\models$``$a'_3(\ell)
= a'_{\eta(8n+5)}$'', we have $\eta(8(n+5)+1)=\gamma$ (see clause (vii)).

Lastly we define $h(a'_n)$ for the other $a'_n$ by induction of $\rk^{N'}(
a'_n)$, note that we can give then dummy elements to relation $\rng(h)\cap
\kappa = \{\eta(8(n+1)+1):n < \omega\}$. 

The model $h[N']$ above should be built in such a way that it is ${\rm
ord}$--transitive. This (and clause (viii)) will ensure that the clause (g) of
the demand \ref{0.1C}(2) is satisfied.

Note that, actually, the coding (of candidates) which we use above does not
change when passing to the ${\rm ord}$--collapse. 

\noindent 2)\quad Should be clear from the above noting: $p \in \bQ$ iff for
some $N$ as above, $N \models\varphi_0(p)$ [as $\Leftarrow$ holds by the
definition and $\Rightarrow$ holds as there are countable $N \prec ({\cal
H}(\chi),\in)$ to which $p,{\frak B},\theta,\kappa$ belong]. 

This finishes the proof of the claim and so the first part of the proposition.

\noindent (2), (3)\quad Easy. \QED$_{\ref{1.3}}$
\medskip

What if in \ref{1.3} we omit ``the only parameters of $\bar \varphi$ are
${\frak B},\theta,\kappa$'', so what do we do?  Well, the role of ${\frak B}$
is assumed by the transitive closure of $\langle\bar{\varphi},{\frak B},
\theta,\kappa\rangle$, which we can then map onto some $\kappa^*\ge\kappa$.

\begin{proposition}
\label{1.5}
\begin{enumerate}
\item Assume $\ZFC^-_*$ is $\emptyset$-normal for $({\frak B},\bar{\varphi},
\theta)$, and, in $\V$, $\bar{\varphi}$ is a $({\frak B},\theta)$--definition
of an [explicit] nep forcing notion. Then we get ``correctly''.
\item Assume $\bar{\varphi}$ is a $K$--$({\frak B},\theta)$--definition of a
nep-forcing notion $\bQ$ (the ``nep'' part is not really needed). Let $\V'$ be
a transitive class of $\V$ such that 
\begin{enumerate}
\item[(i)]  $\bar{\varphi}$ and ${\frak B}$ belong to $\V'$ (and of course
${\frak C}$),
\item[(ii)]  the family of $({\frak B},\bar{\varphi},\theta)$--candidates is
unbounded in $\V'$, moreover 
\item[(iii)] for $\chi$ large enough, in $\V$ (or just in $\V'$) the set
\[\{N\subseteq ({\cal H}(\chi),\in):N\mbox{ is a $({\frak B},\bar{\varphi},
\theta)$--candidate }\}\]
is stationary, or at least 
\item[(iii)$^-$] $\V'\models$``$\varphi_\ell(\bar x)$'' implies that

for unboundedly many $({\frak B},\bar{\varphi},\theta)$--candidates $N$ (in
$\V$), $N\models$``$\varphi_\ell(\bar{x})$'', 

which in other words says that in the universe $\V'$, $({\frak B},
\bar{\varphi},\theta)$ is a correct definition of [explicitly] nep forcing
notion (see part (1) above and Definition \ref{0.2}(11)).  
\end{enumerate}
{\em Then}:
\begin{enumerate}
\item[(a)]  if $\V'\models$``$p\in \bQ$'' (i.e.\ $\varphi_0(p)$) then $\V
\models$``$p \in \bQ$'',
\item[(b)]  if $\V'\models$``$p\le^Q q$'' (i.e.\ $\varphi_1(p,q))$ then $\V
\models$``$p\le^\bQ q$'',
\item[(c)]  if in $\V'$, $N$ is a $({\frak B},\bar{\varphi},\theta)$--candidate
{\em then} also in $\V$, $N$ is a $({\frak B},\bar{\varphi},\theta)$--candidate.
\end{enumerate}
\item If in (2) we add ``explicitly" {\em then} 
\begin{enumerate}
\item[(d)]  if $\V'\models \varphi_2(\langle p_i:i \le \omega \rangle)$
then $\V\models \varphi_2(\langle p_i:i \le \omega \rangle)$,
\item[(e)]  if in $\V'$, $N$ is a $(\bar{\varphi},{\frak B})$--candidate and
$q$ is explicitly $(N,\bQ)$-generic {\em then} this holds in $\V$.
\end{enumerate}
\item If in (2) we add ``$\bar{\varphi}$ is a temporary explicit correct
$({\frak B},\theta)$--definition of a nep forcing notion'' (in $\V$) {\em
then} also in $\V'$, $\bar{\varphi}$ is a temporary explicit correct
$({\frak B},\theta)$--definition of a nep-forcing notion, 
\begin{enumerate}
\item[$(*)_3$]  $\kappa=\theta = \aleph_0$ {\em or} $\kappa=\aleph_0$ and
$([\theta]^{\le \aleph_0})^{\V'}$ is cofinal in $([\theta]^{\le \aleph_0})^{
\V_1}$ {\em or} (there are large enough cardinals to guarantee) any
co--$(\kappa + \theta + \aleph_1)$--Souslin relation in $\V'$ is upward
absolute to $\V_1$. 
\end{enumerate}
\item If in (2) we add $(*)_4$ below and we add ``local'' to the assumption,
{\em then} also in $\V'$, $\bar \varphi$ is a temporary explicit $({\frak B},
\theta)$--definition of a local nep-forcing notion, where 
\begin{enumerate}
\item[$(*)_4$]  $([\kappa \cup \theta]^{\le \aleph_0})^{\V'}$ is cofinal
in $([\kappa \cup \theta]^{ \le \aleph_0},\subseteq)^{\V}$.
\end{enumerate}
\end{enumerate}
\end{proposition}

\Proof  1)\quad  Straight. 
\medskip

\noindent 2)\quad There are two implications implicit in \ref{1.5}(2)
concerning the versions of clause (iii). Let 
\[\begin{array}{lr}
S_\chi\stackrel{\rm def}{=}\{N:&N\in \V,\ N\mbox{ is a countable submodel of
}({\cal H}(\chi),\in)^{\V}\ \ \\
\ &\mbox{and $N$ is a }({\frak B},\bar{\varphi},\theta)\mbox{--candidate }\}
  \end{array}\]
and let 
\[\begin{array}{lr}
S^\prime_\chi\stackrel{\rm def}{=}\{N:&N\in \V',\ N\mbox{ is a countable
submodel of }({\cal H}(\chi),\in)^{\V'}\ \ \\
\ &\mbox{and $N$ is a }({\frak B},\bar{\varphi},\theta)\mbox{--candidate }\}.
  \end{array}\]
Let $<^*_\chi\in \V$ be a well ordering of ${\cal H}(\chi)$ and let 
\[\begin{array}{lr}
C\stackrel{\rm def}{=}\{N:&N\in \V,\ N\mbox{ is a countable elementary
submodel of }\ \ \\
\ &({\cal H}(\chi)^\V,\in,\V\cap{\cal H}(\chi),<^*_\chi)\mbox{ to which
$({\frak B},\bar{\varphi},\theta)$ belongs }\}. 
  \end{array}\]

\noindent{\em If clause (iii) for $\V$ then clause (iii) for $\V'$.}\qquad
Why? Just observe that
\begin{enumerate}
\item[$(\bigstar)_1$] in $\V$:\quad $C$ is a club of $[{\cal H}(\chi)]^{\leq
\aleph_0}$ and $\{N\cap\V': N\in C\}$ is a club of $[{\cal H}(\chi)^{\V'}]^{
\leq\aleph_0}$. 
\end{enumerate}
Now suppose that, in $\V'$, $C'$ is a club of ${\cal H}(\chi)^{\V'}$ and we
should prove $C'\cap S_\chi'\neq\emptyset$ (say for some model ${\frak B}\in
\V'$ with countable vocabulary, the universe ${\cal H}(\chi)^{\V'}$ and Skolem
functions, $C'=\{N:\ N\prec{\frak B}$ countable $\}$\/). As $S_\chi'$ is
stationary in $\V$, also $C_1=\{N\in C: N\cap\V'\in C'\}$ is club of $[{\cal
H}(\chi)^{\V'}]^{\leq\aleph_0}$ in $\V$. Hence there is $N\in S_\chi\cap
C_1$. Now, $N\cap \V'$ is almost a member of $C'\cap S_\chi'$, it satisfies
the requirements in the definitions of $C'$ and $S_\chi'$. But $N\cap\V'$ is a
countable subset of ${\cal H}(\chi)^{\V'}$, so by Shoenfield--Levy
absoluteness it exists.  

\noindent{\em If clause (iii) for $\V'$ then clause (iii) for $\V$.}\qquad
Work in $\V'$. So let $x\in\bQ$ and let $\chi$ be large enough such that
$({\cal H}(\chi),\in)\prec_{\Sigma_n}\V'$ for $n$ large enough. The set 
\[\begin{array}{lr}
C^*\stackrel{\rm def}{=}\{N:&N\in \V,\ N\mbox{ is a countable elementary
submodel of }\ \ \\
\ &{\cal H}(\chi)\mbox{ to which $x,{\frak B},\bar{\varphi},\theta,{\frak C}$
belong }\}.  
  \end{array}\]
is a club of $[{\cal H}]^{\leq\aleph_0}$, hence has non-empty intersection with
any stationary subset of $[{\cal H}]^{\leq\aleph_0}$. In particular, by the
assumption, there is a $({\frak B},\bar{\varphi},\theta)$--candidate $N\in
C^*$. So $N\prec({\cal H}(\chi),\in)$, $x,{\frak B},\bar{\varphi},\theta,
{\frak C}\in N$. So 
\[\V'\models\varphi_0(x)\quad\Rightarrow\quad({\cal H}(\chi),\in)\models
\varphi_0(x)\quad\Rightarrow\quad N\models\varphi_0(x)\quad\Rightarrow\quad
x\in\bQ^N.\]

\noindent 3)\quad  Straight.\medskip
 
\noindent 4)\quad Suppose that 
\[\V'\models\mbox{`` }N \mbox{ is a }(\bar{\varphi},{\frak
B})\mbox{--candidate and } p \in \bQ^N\mbox{ ''}.\]
In $\V'$, let $\langle {\cal I}_n:n<\omega\rangle$ list the ${\cal I}$ such
that $N\models$``${\cal I}$ is a predense subset of $\bQ$''.  We know (by
\ref{1.5}(2)(c)) that $N$ is a candidate in $\V_1$. Hence, in $\V_1$, there
are $q,\langle p^n_\ell:\ell<\omega,n<\omega\rangle$ such that:
\begin{enumerate}
\item[(i)]   $\langle p^n_\ell:\ell<\omega\rangle$ lists ${\cal I}_n \cap N$,
\item[(ii)]  $p\le^{\bQ} q \in \bQ$,
\item[(iii)] $\varphi_2(\langle p^n_\ell:\ell<\omega\rangle\conc\langle q
\rangle)$ for each $n<\omega$.
\end{enumerate}
So there is a $({\frak B},\bar{\varphi},\theta)$--candidate $N_1$ such that
$N\in N_1$, $\langle p^n_\ell:\ell<\omega\rangle: n<\omega\rangle$, $q$ and
$\langle {\cal I}_n: n<\omega\rangle$ belong to $N_1$, and $N_1\models$``$p
\leq^{\bQ} q$'', and $N_1\models\varphi_2(\langle p^n_\ell:\ell<\omega\rangle
\conc \langle q\rangle)$ for $n<\omega$ (by ``correct''). It is enough to find
such $N_1\in\V'$, which follows from \ref{0.11}.\\
(We use an amount of downward absoluteness which holds as $\V'$ is a
transitive class including enough ordinals). 
\medskip

\noindent 5)\quad  Similar proof.  \QED$_{\ref{1.5}}$

\begin{proposition}
\label{1.6}
\begin{enumerate}
\item Assume $\bar{T}$ is an explicit temporary $(\kappa,\theta)$--definition
of a snep--forcing notion $\bQ$. For any extension $\V_1$ of $\V$, this still
holds if $(*)_3$ of \ref{1.5}(4) above holds. So we can replace ``temporary''
by $K=$ class of all set forcing notions. 
\item Assume $(\bar \varphi,{\frak B})$ is a simple explicit temporary
$(\kappa,\theta)$--definition of a nep--forcing notion $\bQ$. For any
extension $\V_1$ of $\V$ this still holds in $\V_1$ if $(*)_3$ of \ref{1.5}(4)
holds. So we can add/replace ``temporary" by the class $K$ of all forcing
notions preserving ``$([\theta]^{\le \aleph_0})^{\V}$ is cofinal in
$([\theta]^{\le \aleph_0})^{\V_1}$''. 
\item Assume $(\bar{\varphi},{\frak B})$ is a local explicit temporary
$(\kappa,\theta)$--definition of a nep forcing notion $\bQ$. {\em Then} for
any extension $\V_1$ of $\V$ this still holds, provided that:
\begin{enumerate}
\item[$(*)_4$]\quad $([\kappa\cup\theta]^{\leq\aleph_0})^{\V}$ is cofinal in
$([\kappa\cup\theta]^{\leq\aleph_0},\subseteq)^{\V_1}$.
\end{enumerate}
\end{enumerate}
\end{proposition}

\Proof 1), 2)\quad Left to the reader (and similar to the proof of part (3)). 

\noindent 3)\quad Let $\theta_1=\kappa+\theta$, and let $a\in [\theta_1]^{
\aleph_0}$, and consider the statement
\begin{enumerate}
\item[$\boxtimes_a$]  if $N$ is a $({\frak B},\bar{\varphi},\theta
)$--candidate satisfying $N\cap\theta_1\subseteq a$ and $p\in\bQ^N$ (i.e.\ $N
\models\varphi_0[p]$),\\
{\em then} there are $N'$, a generic extension of $N$ (so have the same
ordinals and $N$ is a class of $N'$) which is a $({\frak B},\bar{\varphi},
\theta)$--candidate such that $N'\models$``${\cal P}(\theta)^N$ is countable''
and 
\[N'\models\mbox{`` }(\exists q)[q \in\bQ\ \&\ q \mbox{ is explicitly } (N
\cap {\cal P}(\bQ),\bQ)\mbox{--generic] ''}.\]
\end{enumerate}
Note:\quad for $[x\in N'\wedge N'\models$``$x$ is countable''\quad$\Rightarrow
\quad x\subseteq N']$ just use a suitable collapse.

Now, only $(N',c)_{c \in N}/\cong$ and $(N,\alpha)_{\alpha \in a}/\cong$ and
$(N',N,\alpha)_{\alpha \in a}$ are important and we can code $N$ as a subset
of $a$ (as all three are countable). Thus the statement is essentially
\[\begin{array}{r}
(\forall N)[(N \mbox{ is not well founded (or not } {\frak B} \rest
(N \cap a) \prec_\Delta {\frak B}, \mbox{ etc}.) \vee\\
\vee (\exists N')(N' \mbox{ as above})].
  \end{array}\]
So it is $\Pi^1_2$, hence it is absolute from $\V$ to $\V_1$. Now, both in
$\V'$ and in $\V_1$ the statement ``$\bQ$ is simply, locally, explicitly nep''
is equivalent to $(\forall a \in [\theta_1]^{\aleph_0})\boxtimes_a$, which is
equivalent to ${\cal S}=\{a\in [\theta_1]^{\aleph_0}:\boxtimes_a\}$ is cofinal
in $[\theta_1]^{\aleph_0}$.  But by the previous paragraph ${\cal S}[\V] 
\subseteq {\cal S}[\V_1]$.  Now $(*)_4$ gives the needed implication.
\QED$_{\ref{1.6}}$

\begin{proposition}
\label{1.7}
\begin{enumerate}
\item Assume $\bar{T}$ is a temporarily $(\kappa,\theta)$--definition of a
snep--forcing notion $\bQ$. If $(*)$ below holds, {\em then} we can find a tree
$T_2 \subseteq {}^{\omega >}(\theta \times\theta \times \kappa')$ such that
$\bar{T}\conc\langle T_2 \rangle$ is an explicit temporary
$(\kappa,\theta)$--definition of a snep-forcing notion $\bQ$, where
\begin{enumerate}
\item[$(*)$] $\kappa=\theta=\aleph_0$, $\kappa' = \aleph_1$ or enough
absoluteness. 
\end{enumerate}
\item If $\bQ$ (i.e.\ $\langle \varphi_0,\varphi_1 \rangle$) is a Souslin
proper forcing notion (see \ref{0.7}) and ${\frak B}$ codes the parameter (so
has universe $\kappa=\aleph_0$ and let $\theta = \aleph_0$), {\em then}
$(\bar{\varphi},{\frak B})$ is a simple explicit temporary $(\kappa,\theta)$--definition of the nep-forcing notion $\bQ$.
\end{enumerate}
\end{proposition}

\Proof 1)\quad The question is to express ``$\{p_n:n < \omega\}$ is predense
above $q$'' which is equivalent to 
\[(\forall \nu \in {}^\omega \theta)[\varphi_0(\nu)\ \&\ \varphi_1(q,\nu)\ 
\Rightarrow\ (\exists\nu'\in {}^\omega \theta)(\bigvee_n \varphi_1(p_n,\nu')\
\&\ \varphi_1(\nu,\nu'))].\]
So, as $\kappa = \theta = \aleph_0$, this is a $\Pi^1_2$-formula and hence it
is $\aleph_1$-Souslin. 

\noindent 2)\quad Similarly (for $\varphi_2$ being upward absolute note that
the relation is now $\Pi^1_1$ and $\Pi^1_1$ formulas are upward absolute).
\QED$_{\ref{1.7}}$

\begin{definition}
\label{1.9}
Assume that $(\bar{\varphi},{\frak B})$ is a temporary
$(\kappa,\theta)$--definition of a nep forcing notion $\bQ$, and $N$ is a
$\bQ$--candidate. We say that a condition $q'\in\bQ$ is essentially explicitly 
$(N,\bQ)$-generic if for some candidate $N'$, $N \subseteq N'$, $N \in N'$,
$q'$ is explicitly $(N',\bQ)$-generic and for some $q_0\in \bQ^{N'}$, $q_0
\le^{\bQ} q'$ and $N' \models$``$q_0$ is $(N,\bQ)$-generic''. 
\end{definition}
Note: if $\bQ$ is a snep-forcing for $\bar{T}$, this relation is $(\kappa +
\theta+\aleph_1)$--Souslin, too. 

\begin{proposition}
\label{1.8}
Assume $\bQ$ is a correct explicitly nep-forcing notion, say by $(\bar{\varphi},
{\frak B})$. If $q$ is $(N,\bQ)$-generic, {\em then} for some $q'$ we have
\[q \le q' \in \bQ\quad \mbox{and}\quad q' \mbox{ is essentially explicitly }
(N,\bQ)\mbox{-generic}.\]
\end{proposition}

\Proof Let $\varphi_2(\langle p_n^{\cal I}:n<\omega\rangle,q)$ hold for some
list $\langle p_n^{\cal I}: n<\omega\rangle$ of ${\cal I}\in\pd(N,\bQ)$. For
${\cal I}\in\pd(N,\bQ)$ let $N_{\cal I}$ be a $\bQ$--candidate such that
$N_{\cal I}\models\varphi_2(\langle p_n^{\cal I}:n<\omega\rangle,q)$.  

Let $N' \subseteq {\cal H}(\chi)$ be a countable $\bQ$--candidate satisfying
\[\{N,q\}\cup\{N_{\cal I}:{\cal I}\in\pd(N,\bQ)\}\in N'.\]
By our assumptions there is $q'$ such that: $q\le q'\in\bQ$ and $q'$ is
explicitly $(N',\bQ)$-generic. \QED 

\begin{proposition}
\label{1.10}
\begin{enumerate}
\item If $N$ is a ${\frak B}$--candidate, so in particular
\[[N \models\mbox{``}\alpha<\kappa\vee \alpha<\theta\mbox{''}]\ \Rightarrow\
\alpha\in\kappa \vee \alpha\in\theta,\]
and $|{\frak B}|$ is an ordinal, {\em then} there is a unique $N' = \mbox{
MosCol}_{\kappa,\theta}(N)$ and $f$ such that 
\begin{enumerate}
\item[(a)] $f$ is an isomorphism from $N$ onto $N'$,
\item[(b)] $f\rest (N \cap \kappa) = {\rm id}$, $f(\kappa)=\kappa$ and $f
\rest (N \cap \theta)={\rm id}$, $f(\theta)=\theta$,
\item[(c)] if $x\in N \backslash(\kappa+1)\backslash(\theta +1)$ then $f(x)=
\{f(y):N \models$``$y \in x$''$\}$, 
\item[(d)] $N'$ is a ${\frak B}$--candidate.
\end{enumerate}
\item Note that if $N \models$``$x \in {\cal H}_{< \aleph_1}(\kappa)\cup
{\cal H}_{< \aleph_1}(\theta)$'' then $f(x) = x$.
\item If $|{\frak B}|$ is not an ordinal (so $\kappa \subseteq |{\frak B}|
\subseteq {\cal H}_{< \aleph_1}(\kappa)$), then $N'$ is still a $({\frak B},
\bar{\varphi},\theta)$--candidate, using the ``but'' of clause (e) of
Definition \ref{0.1C}(2). 
\end{enumerate}
\end{proposition}

\begin{fact}
\label{1.11}
In the definition of nep (or snep) in the ``properness'' clause, it is enough
to restrict ourselves to a family $\bf I$ of predense subsets of $\bQ^N$ such
that: 
\[\begin{array}{l}
\mbox{if }{\cal I}\in\pd(N,\bQ)\\
\mbox{then for some }{\cal J}\in{\bf I}\mbox{ we have }(\forall p \in {\cal I}
\cap N)(\exists q \in {\cal J})(N \models p \le^{\bQ} q).
  \end{array}\]
\end{fact}

\begin{proposition}
\label{2.11}
\begin{enumerate}
\item Assume $\bar{T}$ defines an explicit $(\kappa,\theta)$--snep forcing
notion. Let $\bar{\varphi}=\bar{\varphi}_T$, ${\frak B}={\frak B}_{\bar{T}}$
(see \ref{1.1}(3)). If $\bQ^{\bar{T}}$ is local then $\bQ^{\bar{\varphi}}$ is
local, in fact in Definition \ref{0.5A}(2).
\item If ($\ZFC^-_*$ is $K$--good and) $\ZFC^-_*$ says that $({\frak
B},\bar{\varphi},\theta)$ is explicitly nep, and $\bar{\varphi}$ is correct
then $({\frak B}^{\bar{\varphi}},\bar{\varphi},\theta)$ is explicitly nep and
local. \QED
\end{enumerate}
\end{proposition}

Moving from nep to snep (and inversely) we may ask what occurs to ``local".
It is usually preserved.  

\stepcounter{section}
\subsection*{\quad 3. There are examples} 
In this section we show that a large family of natural forcing notions
satisfies our definition. Later we will deal with preservation theorems but to 
get nicer results we better ``doctor" the forcing notions, but this is delayed
to the next section.
\medskip

\noindent In fact all the theorems of Ros{\l}anowski Shelah \cite{RoSh:470},
which were designed to prove properness, actually give one notion or another
from \S1 here (confirming the thesis \ref{0.1B} of \S0).  We will state them
without giving the definitions from \cite{RoSh:470} and give a proof of
(hopefully) well known specific cases, indicating why it works.    

\begin{lemma} [Ros{\l}anowski Shelah \cite{RoSh:470}] 
\label{2.1}
\begin{enumerate}
\item Suppose that $\bQ$ is a forcing notion of one of the following types:
\begin{enumerate}
\item[(a)] $\bQ^{\rm tree}_e(K,\Sigma)$ for some finitary tree-creating pair
$(K,\Sigma)$, where $e=1$ and $(K,\Sigma)$ is 2-big {\em or} $e=0$ and
$(K,\Sigma)$ is t-omittory (see \cite[\S2.3]{RoSh:470}; so e.g.\ this covers
the Sacks forcing notion), 
\item[(b)] $\bQ^*_{{\rm s}\infty}(K,\Sigma)$ for some finitary creating pair
$(K,\Sigma)$ which is growing, condensed and of the AB--type or omittory, of
the AB$^+$--type and satisfies $\oplus_0,\oplus_3$ of \cite[4.3.8]{RoSh:470}
(see \cite[\S3.4]{RoSh:470}; this captures the Blass--Shelah forcing notion
of \cite{BsSh:242}),
\item[(c)] $\bQ^*_{{\rm w}\infty}(K,\Sigma)$ for some finitary creating pair
which captures singletons (see \cite[\S2.1]{RoSh:470}) 
\item[(d)] $\bQ^*_f(K,\Sigma)$ for some finitary, 2-big creating pair
$(K,\Sigma)$ with the Halving Property which is either simple or gluing 
and an $H$-fast function $f:\omega\times\omega\longrightarrow\omega$ (see
\cite[\S2.2]{RoSh:470}). 
\end{enumerate}
{\em Then} $\bQ$ is an explicit $\aleph_0$--snep forcing notion, moreover, it
is local.
\item Assume that $\bQ$ is a forcing notion of one of the following types:
\begin{enumerate}
\item[(a)] $\bQ^{\rm tree}_e(K,\Sigma)$ for $e<3$ and a tree-creating pair
$(K,\Sigma)$, which is bounded if $e=2$ (see \cite[\S2.3]{RoSh:470}; this
includes the Laver forcing notion),
\item[(b)] $\bQ^*_\infty(K,\Sigma)$ for a finitary growing creating pair
$(K,\Sigma)$ (see \cite[\S2.1]{RoSh:470}; this covers the Mathias forcing
notion).  
\end{enumerate}
{\em Then} $\bQ$ is an explicit $\aleph_0$--nep forcing notion, moreover, it
is local. 
\end{enumerate}
\end{lemma}

\Proof Let $N$ be a $\bQ$-candidate and $p\in\bQ^N$. Let $\langle {\cal J}_n:
n<\omega\rangle$ list $\{{\cal J}:N\models ``{\cal J}\subseteq\bQ$ is open
dense''$\}$. Then there is a sequence $\langle (p_n,{\cal I}_n):n<\omega
\rangle$ such that $p_n,{\cal I}_n\in N$, $N\models p_n\le p_{n+1}$, ${\cal
I}_n \subseteq {\cal J}_n$ is a countable set, $\langle p_n:n<\omega\rangle$
has an upper bound in $\bQ$ and ${\cal I}_n$ is predense above $p_{n+1}$,
moreover, in an explicit way as described below (see the respective
subsections in \cite{RoSh:470}). Moreover,
\begin{quotation}
\noindent in part (1) cases (a)+(c), ${\cal I}_n$ is finite and moreover, we
can say ``${\cal I}_n$ is predense above $p_{n+1}$" in a Borel way.
\end{quotation}
For the Sacks forcing notion:\quad for some $k<\omega$, ${\cal I}_n =
\{p^{[\eta]}_{n+1}:\eta\in p_{n+1},\lh(\eta)=k\}$,  so ${\cal I}_n$
corresponds to a front of $p_{n+1}$, which necessarily is finite. This
property serves as $\varphi_2$ (compare with more detailed description for the
Laver forcing below).

In part (1) case (b) (e.g. the Blass--Shelah forcing notion) ${\cal I}_n$ is
countable. We do not know which level will be activated, but if use $n$, then
we get into ${\cal I}_n$, so ${\cal I}_n$ countable but the property is Borel
not $\Pi^1_1$.

Now, in part (2), ${\cal I}_n$ is countable and again it corresponds to some
front $A$ of $p_{n+1}$ in an appropriate sense. So ${\cal I}_n=\{p^{[\eta]}_{
n+1}:\eta \in A\}$, but to say ``$A$ is a front'' is $\Pi^1_1$ (in some
instances of 2(a) we have $e$-thick antichains instead of fronts, but the
complexity is the same).  

Recall that for a subtree $T \subseteq {}^{\omega >}\omega,A \subseteq T$ is
a front of $T$ if 
\[(\forall\eta\in\lim(T))(\exists n)(\eta \rest n \in A)\] 
(usually members of $A$ are pairwise incomparable).
\medskip

\noindent Specifically, for the Laver forcing notion, we can guarantee ${\cal
I}_n = \{p^{[\eta]}_{n+1}:\eta\in A\}$, where $A$ is a front of $p_{n+1}$. Now
being a front is a $\Pi^1_1$--sentence (see the definition above) which is
upward absolute and this is our choice for $\varphi_2$.  Let us write this
formula in a more explicit way (for the case of the Laver forcing notion):
\[\begin{array}{l}
\varphi_2(\langle p_i:i \le \omega \rangle) \equiv \mbox{ each }p_i \mbox{ is
a Laver condition and}\\
\qquad\bigwedge\limits_{i\in \omega}(\exists! \eta)(\eta \in p_\omega\ \&\
p_{2i}=p^{[\eta]}_\omega)\\
\mbox{[call this unique $\eta$ by $\eta_i$]}\quad \mbox{ and}\\
\qquad\bigwedge\limits_{i\neq j} \eta_i \ntrianglelefteq \eta_j \mbox{
(incomparable) }\ \&\ (\forall\rho\in\lim(p_\omega))(\bigvee\limits_n
\bigvee\limits_m \rho \rest n = \eta_m)\\
\mbox{[this is: $\{p_i:i\in\omega\}$ is explicitly predense above
$p_\omega$]}.  
  \end{array}\]
So it is $\Pi^1_1$ (of course, $\Sigma^1_2$ is okay, too.)\QED$_{\ref{2.1}}$
\medskip

Note that even for the Sacks forcing notion, ``$p,q$ are incompatible'' is
complete $\Pi^1_1$. So ``$\{p_n:n \in \omega\}$ is predense above $p$'' will
be $\Pi^1_2$. For Laver forcing we cannot do better. Now, generally $\Pi^1_2$
is not upward absolute from countable submodels, whereas $\Pi^1_1$ is.
\bigskip

\begin{proposition}
\label{2.1B}
All the forcing notions $\bQ$ defined in \cite{RoSh:470}, \cite{RoSh:628}, are
correct, and we can use $\ZFC^-_*={\rm ZC}^-$ which is good and normal (see
\ref{0.9}). Also the relation ``$p,q$ are incompatible members of $\bQ$'' is
upward absolute from $\bQ$--candidates (as well as $p \in \bQ$, $p \notin
\bQ$, $p \le q$, and ``$p,q$ are compatible'').
\end{proposition}

\Proof Check. \QED

\stepcounter{section}
\subsection*{\quad 4. Preservation under iteration: first round}
We give here one variant of the preservation theorem, but for it we need some
preliminary clarification. We have said ``there is $q$ which is
$(N,\bQ)$-generic"; i.e.\ $q \forces_{\bQ}$`` $\name{G}_{\bQ}\cap \bQ^N$ is a
generic subset of $\bQ^N$ over $N$ ''. Note that we have said $\bQ^N$ and not
$\bQ\cap N$ as we intended to demand $N\models$``$r\in\bQ$''\quad
$\Rightarrow\quad \V \models$``$r \in \bQ$'' rather than $r\in N\ \
\Rightarrow\ \ [N \models$``$r \in \bQ$''$\ \Leftrightarrow\ \V\models$``$r
\in\bQ$''] (the version we use is, of course, weaker and so better). Now, to
use Definition \ref{0.2}(1) we usually use $N[G_Q]$ (e.g.\ when iterating).  

But what is $N[G]$ here? In fact, what is the connection between
$N\models$``$\name{\tau}$ is a $\bQ$--name'' and $\V\models$``$\name{\tau}$ is
a $\bQ$--name''? Because $[x \in Y \in N \nRightarrow x \in N]$, none of the
implications holds. 

For our purpose, the usual $N[G]=\{\name{\tau}[G]:\name{\tau}\in N \mbox{ is a
}\bQ\mbox{-name}\}$ is not appropriate as it is not clear where being a
$\bQ$-name is defined. We use $N\langle G \rangle$ which is $N[G\cap\bQ^N]$
when we disregard objects in $\V\backslash N$. Of course, if the models are 
$\subseteq {\cal H}_{<\aleph_1}(\kappa\cup\theta)$ life is easier; but we
may lose $N \models\ZFC^-_*$.

We then prove (in \ref{2.5}) the first version of preservation by CS
iteration. We aim at proving only that $\bP_\alpha=\Lim(\bar{\bQ})$ satisfies
the main clause, i.e.\ clause (c) of Definition \ref{0.2} (but did not say
that $\bP_\alpha$ is nep itself). For this we need again to define what is
$N\langle G\rangle$ for $N$ which is not necessarily a candidate. The second
treatment (in \S5) depends just on Definition \ref{2.3} from this section. 

\begin{definition}
\label{2.3}
\begin{enumerate}
\item Assume $N\models$``$\bQ$ is a nep-forcing notion'' and $G\subseteq\bQ^N$
is generic over $N$. We define $N\langle G\rangle=N\langle G\cap\bQ^N\rangle$
``ignoring $\V$" and letting ${\frak B}^{N \langle G\rangle}={\frak B}^N$ for
the relevant ${\frak B}$. In details, 
\[N \langle G \rangle\stackrel{\rm def}{=}\{\name{\tau}^N\langle G\rangle:
N \models\mbox{``}\name{\tau}\mbox{ is a $\bQ$--name''}\},\]
where $\name{\tau}^N\langle G\rangle$ is defined by induction on $\rk^N(
\name{\tau})$ (see e.g.\ \cite[Ch.I]{Sh:f}):
\begin{enumerate}
\item[(a)] if for some $p\in G\cap \bQ^N$ and $x \in N$ we have $N\models[p
\forces_{\bQ}$``$\name{\tau} = x$''] then $\name{\tau}^N\langle G\rangle=x$,
\item[(b)] if not (a) then necessarily $N\models$``$\name{\tau}$ has the form
$\{(p_i,\name{\tau}_i):i<i^*\}$, $p_i \in \bQ$, $\name{\tau}_i$ a $\bQ$--name
of rank $<\rk(\name{\tau})$''; now we let 
\end{enumerate}
\[\name{\tau}^N\langle G\rangle=\{(\name{\tau}')^N\langle G\rangle: \name{
\tau}'\in N\mbox{ and for some }p \in G \cap \bQ^N\mbox{ we have }(p,
\name{\tau}')\in\name{\tau}\}.\]
\item If $N \models$``$\name{\tau}$ is a $\bQ$--name'' we define a $\bQ$--name
$\name{\tau}^{\langle N\rangle}$ as follows: 
\begin{enumerate}
\item[(a)] if $N\models$``$\name{\tau}=\check{x}$'', $x \in N$, we let
$\name{\tau}=\check{x}$ (see e.g.\ \cite[Ch.I]{Sh:f}),
\item[(b)] if $N\models$``$\name{\tau}=\{(p_i,\name{\tau}_i):i<i^*\}$, where
$p_i\in\bQ$, $\name{\tau}_i$ a $\bQ$--name of rank $<\rk(\name{\tau})$'' then 
\[\name{\tau}^{\langle N\rangle}=\big\{(p,(\name{\tau}')^{\langle N\rangle}):
N\models\mbox{``}(p,\name{\tau}')\in \name{\tau}\mbox{''}\big\}.\]
\end{enumerate}
\item We say ``$q$ is $\langle N,\bQ\rangle$--generic" if $q\forces_{\bQ}$``$
\name{G}_{\bQ} \cap\bQ^N$ is a subset of $(\bQ^N,<^N_{\bQ})$ generic over
$N$''.  
\end{enumerate}
\end{definition}

\begin{definition}
\label{2.3A}
\begin{enumerate}
\item In Definition \ref{0.2}(1) replacing ``temporarily'' by
``$K$--absolutely'' means 
\begin{enumerate}
\item[(a)] if $\V_1$ is a $K$--extension of $\V$ (i.e.\ a generic extension of
$\V$ by a forcing notion from $K^{\V}$) {\em then} 
\begin{enumerate}
\item[(i)]   $\V\models$``$x\in\bQ^{\bar{\varphi}}\mbox{''}\quad\Rightarrow
\quad\V_1\models$``$x\in \bQ^{\bar{\varphi}}$'',
\item[(ii)]  $\V\models$``$x<^{\bQ^{\bar{\varphi}}} y\mbox{''}\quad
\Rightarrow\quad\V_1\models$``$x<^{\bQ^{\bar{\varphi}}} y$'',
\item[(iii)] in the explicit case we have a similar demand for $\varphi_2$;
otherwise, if $N$ is a $\bQ^{\bar{\varphi}}$--candidate in $\V$, $q\in
\bQ^{\bar{\varphi}}$ is $\langle N,\bQ \rangle$--generic (see \ref{2.3}(3)) in
$\V$ {\em then} $q$ is $\langle N,\bQ \rangle$--generic in $\V_1$,
\end{enumerate}
\item[(b)] if $\V_1$ is a $K$--extension, {\em then} the relevant part of
Definition \ref{0.2} and clause (a) here holds in $\V_1$,
\item[(c)] if $\V_{\ell +1}$ is a $K$--extension of $\V_\ell$ for $\ell\in
\{0,1,2\}$, $\V_0 = \V$ {\em then} $\V_3$ is a $K$--extension of $\V_1$.
\end{enumerate}
\item We omit $K$ when we mean: any set forcing.
\end{enumerate}
\end{definition}

Note that (a)(i) + (ii) is automatic for explicitly snep, also (a)(iii). One
can  make ``absolutely nep" to the main case.  

The following is natural to assume.

\begin{definition}
\label{2.3B}
\begin{enumerate}
\item We say $\ZFC^-_*$ is nice to $\chi_1$ if $\chi_1$ is a constant in
${\frak C}$, $\ZFC^-_*$ says $\chi_1$ is strong limit and $\ZFC^-_*$ is
preserved by forcing by forcing notions of cardinality $<\chi_1$. 
\item We say $\bQ$ is nice (or $\ZFC^-_*$ nice to $\bQ$) if for some $\chi_1$,
$\ZFC^-_*$ is nice to $\chi_1$ and it says $\bQ \in {\cal H}(\chi_1)$.
\end{enumerate}
\end{definition}

\begin{proposition}
\label{2.4}
If $N$ is a $\bQ$--candidate, $\bQ$ is a nep-forcing notion, $G_{\bQ}
\subseteq \bQ$ is generic over $\V$ and $G_{\bQ} \cap\bQ^N$ is generic over
$N$ then: 
\begin{enumerate}
\item[(a)] $N\models$``$\name{\tau}$ is a $\bQ$-name'' implies $\name{\tau}^N
\langle G\rangle=\name{\tau}^{\langle N\rangle}[G]$,
\item[(b)] $N\langle G\rangle$ is a model of $\ZFC^-_*$ and moreover it is a
$\bQ$--candidate and is a forcing extension of $N$, provided that the forcing 
theorem applies, i.e.\ $\ZFC^-_*$ is $K$--good, $\bQ\in K$ (see Definition
\ref{0.9}),
\item[(c)] $N\langle G\rangle\cap\kappa=N\cap\kappa$, $N\langle G\rangle \cap
\theta = N\cap\theta$. \QED
\end{enumerate}
\end{proposition}

\noindent{\bf Remark:}\quad It seems that usually (but not in general) we
have:
\[({\cal H}_{< \aleph_1}(\kappa))^{\V[G]} \cap N \langle G \rangle={\cal H}_{
<\aleph_1}(\kappa)^{\V}\cap N\quad\mbox{ and}\]
\[({\cal H}_{< \aleph_1}(\theta))^{\V[G]}\cap N\langle G\rangle=({\cal H}_{<
\aleph_1}(\theta))^{\V} \cap N.\]

\begin{proposition}
\label{2.5}
Assume
\begin{enumerate}
\item[(a)] $\bar{\bQ}=\langle \bP_i,\nbQ_j:i \le\alpha,j<\alpha\rangle$ is a
CS iteration, 
\item[(b)] for each $i<\alpha$ 
\[\forces_{\bP_i}\mbox{``}(\bar{\varphi}_i,{\frak B}_i)\mbox{ is a temporary
$(\kappa_i,\theta_i)$--definition of a nep-forcing notion $\nbQ_i$''}\]
and the only parameter of $\bar{\varphi}_i$ is ${\frak B}_i$, so we are
demanding $\langle(\bar \varphi_i,{\frak B}_i):i<\alpha\rangle \in \V$,
\item[(c)] ${\frak B}$ is a model with universe $\alpha^*$, or including
$\alpha^*$ and included in ${\cal H}_{<\aleph_1}(\alpha^*)$, where $\alpha^*
\ge\alpha$, $\alpha^*\ge\kappa_i =\kappa({\frak B}_i)$, ${\frak B}$ codes
$\langle({\frak B}_i,\bar{\varphi}_i):i<\alpha\rangle$ and the functions
$\alpha -1,\alpha+1$. 
\end{enumerate}
We can use a vocabulary $\subseteq \{P_{n,m}:n,m < \omega\}$ where $P_{n,m}$
is an $n$-place predicate to code $\langle {\frak B}^i:i<\alpha\rangle$: let
$P^{\frak B}_{n+1,2m}=\{\langle i,x_1,\ldots,x_n \rangle:\langle x_1,\ldots,
x_n \rangle \in P^{{\frak B}_i}_{n,m}\}$, $P_{2,1} = \{(\alpha,\alpha +1):
\alpha +1 < \alpha^*\}$ (and $\Delta$ is the set of first order formulas).

\noindent{\em Then}:\quad if $N\subseteq ({\cal H}(\chi),\in)$ is a ${\frak
B}$--candidate, $p\in\bP_\alpha\cap N$, then for some condition $q$, $p\le
q\in \bP_\alpha$ and $q$ is $\langle N,\bP_\alpha \rangle$--generic (in
particular $\bP_\alpha$ is defined from ${\frak B}$) which is defined below. 
\end{proposition}

\begin{definition}
\label{2.6}
Under the assumptions of \ref{2.5}, in $N$ we have a definition of the
countable support iteration $\bar{\bQ}=\langle \bP_i,\nbQ_j:i\le\alpha,j<
\alpha\rangle$. We define by induction on $j\in N\cap(\alpha+1)$ when $q\in
\bP_j$ is $\langle N,\bP_j \rangle$--generic: 
\begin{enumerate}
\item[$(\circledast)$]  if $q\in G_j\subseteq\bP_j$ and $G_j$ is generic over
$\V$ {\em then} $G^{\langle N\rangle}_j$ is a generic subset of $\bP^N_j$ over
$N$, where
\[G^{\langle N\rangle}_j\stackrel{\rm def}{=}\{p:N\models{``}p\in \bP_j
\mbox{''\ \ and\ \ } p^{\langle\langle N\rangle\rangle} \in G_j\},\]
where $p^{\langle\langle N\rangle\rangle}$ is a function with domain
$\dom(p)^N$, and $p(\gamma)$ is the following $\bP_\gamma$--name:\quad {\em
if}  $p(\gamma)^{\langle N \langle\name{G}_\gamma\cap N\rangle\rangle}\in
\nbQ_\gamma$, {\em then} it is $p(\gamma)$; {\em if} not, {\em then} it is 
$\emptyset_{\nbQ_\gamma}$.
\end{enumerate}
\end{definition}

\begin{remark}
\label{2.7}
The major weakness is that $\bP_\alpha$ is not proved to be in some of our
classes (nep or snep).  We get the ``original property'' without the ``support
team", i.e.\ the $\nbQ_i$ are nep, but on $\bP_\alpha$ we just say it
satisfies the main part of nep. A minor one is that ${\frak B}_i$ is not
allowed to be a $\bP_i$--name in any way. In the later theorems, we use
$\bP'_\alpha \subseteq \bP_\alpha$ consisting of ``hereditarily countable''
names. 

Note: inside $N$, if ``$N\models p\in\bP_\alpha$'' then $\dom(p_\alpha)\in
[\alpha]^{\le \aleph_0}$ in $N$'s sense hence (see Definition \ref{0.1C}(2)),
$\dom(p_\alpha)\subseteq N$ and similarly the names are actually from $N$,
members outside $N$ do not count, they may not be in $\bP_\alpha$ at all.
\end{remark}

{\noindent{\sc Proof of \ref{2.5}} \hspace{0.2in}} We imitate the proof of the
preservation of properness. So we prove by induction on $j\in (\alpha+1)\cap
N$ that: 
\begin{enumerate}
\item[$(*)_j$] if $i\in j\cap N$, $q$ is $(N,\bP_i)$--generic, and
$q\forces_{\bP_i}$``$(p \rest i)^{\langle N \rangle}\in\name{G}_{\bP_i}$''
{\em then} we can find a condition $r\in \bP_j$ such that $r\rest i=q$, and $r
\forces_{\bP_j}$``$(p \rest j)^{\langle N \rangle}\in \name{G}_{P_j}$'', $r$
is $(N,\bP_j)$--generic, and $\dom(r)\setminus i\subseteq N$''.
\end{enumerate}

\noindent{\sc Case 0}:\quad  $j=0$.\\
Left to the reader.
\medskip

\noindent{\sc Case 1}:\quad  $j = j_1+1$.\\
So $j_1 \in N$ (why? use $P_{2,1}$ and \ref{0.2}(2)(e)), and by the inductive
hypothesis and the form of the conclusion without loss of generality $i
=j_1$. Let $q \in G_i\subseteq \bP_i$, $G_i$ generic over $\V$. So $N\langle
G^{\langle N \rangle}_i \rangle\cap\alpha^*=N\cap\alpha^*$ (by \ref{2.4}),
and hence ${\frak B}\rest N\langle G^{\langle N \rangle}_i \rangle = {\frak B}
\rest N \prec {\frak B}$. But $i\in N$, so this applies to ${\frak B}_i$,
too. So $\V[G_i]\models$``$N \langle G^{\langle N \rangle}_i \rangle$ is a 
${\frak B}_i$--candidate''. Also $N\langle G^{\langle N\rangle}_i\rangle
\models$``$p(i)^{\langle N\langle G^{\langle N\rangle}_i\rangle\rangle}\in
\bQ_i$'' because $G^{\langle N\rangle}_i$ is a generic subset of $P^N_i=\{x:x
\in N$, $N\models$``$x\in \bP_i$''$\}$ over $N$ and use the property of
$\bQ_i$. 
\medskip

\noindent{\sc Case 3}:\quad $j$ is a limit ordinal.\\
As in the proof for properness (see \cite[Ch.III, 3.2]{Sh:f}). 
\QED$_{\ref{2.5}}$ 
\medskip

\noindent{\bf Remark:}\quad Note that if $N \models$``$w$ is a subset of
$\alpha$'' then we can deal with $\bP_w$, as in \S5. 

\stepcounter{section}
\subsection*{\quad 5. True preservation theorems} 
Let us recall that $\bQ$ is nep if ``$p\in\bQ$'', ``$p\le_{\bQ} q$'' are
defined by upward absolute formulas for models $N$ which are $({\frak
B}^\theta,\bar{\varphi}^\theta,\theta^{\bQ})$--candidates; i.e.\
$N\subseteq({\cal H}(\chi),\in)$ countable, ${\frak B}^{\bQ}\in N$ a model on
some $\kappa$, ${\frak B}^{\bQ}\rest N\prec_\Delta {\frak B}^{\bQ}$, $N$ model
of $\ZFC^-_*$ and for each such model we have the properness condition. 
Usually $\bQ\subseteq {}^\omega \theta$, or ${\cal H}_{< \aleph_1}(\theta)$ or
so. We would like to prove that CS iteration preserves ``being nep'', but CS
may give ``too large'' names of conditions (of $\nbQ_i$, $i>0$) depending say
on large maximal antichains (of $\bP_i$). Note: if $\bQ_0$ is not
c.c.c. normally it has maximal antichain which is not absolutely so; start
with a perfect set of pairwise incompatible elements and extend it to a maximal
antichain. Then whenever a real is added, the maximality is lost. Finally,
c.c.c. is normally lost in $\bP_\omega$. So we will revise our iteration so
that we consider only hereditarily countable names. 

But in the iteration, trying to prove a case of properness for a candidate $N$
and $p \in\bP^N_{\alpha+1}$, considering $q\in\bP_\alpha$ which is $\langle N,
\bP^N_\alpha\rangle$--generic, we know that in $\V[G_{\bP_\alpha}]$ (if $q\in
G_{\bP_\alpha})$, there is $q'\in \nbQ_\alpha[G_{\bP_\alpha}]$ which is
$\langle N[G_{\bP_\alpha}],\bQ_\alpha[G_\alpha]\rangle$--generic. But under
present circumstances, we have no idea where to look for $q'$, so no way to
make a name of it, $\name{q}'$, which is hereditarily countable, without
increasing $q \in\bP_\alpha$. Except when $\bQ$ is local (see \ref{0.5A}), of
course; it is not unreasonable to assume it but we prefer not to and even
then, we just have to look for it in, essentially, a copy of the set of
reals. The solution is to increase $\bQ_i$ insubstantially so that we will
exactly have the right element $q'$:  
\[p(\alpha)\ \&\ \bigwedge_{{\cal I}\in \pd_{\bQ}(N)}\,\,\bigvee_{p\in {\cal
I} \cap N} p,\]
as explained below. We give two variants.

\begin{notation}
\label{4.0A}
Let $\pd_{\bQ}(N)=\pd(N,\bQ)=\{{\cal I}:N\models$``${\cal I}$ is a predense
subset of $Q$''$\}$ and ${\cal I}[N] = {\cal I}^N={\cal I}\cap N$.
\end{notation}

\begin{definition}
\label{4.1}
Let $\bQ$ be an explicitly nep-forcing notion. Then we define $\bQ'=\cl(\bQ)$
as follows:
\begin{enumerate}
\item[(a)] the set of elements is
\[\bQ\cup\big\{p\ \&\ \bigwedge_{{\cal I} \in \pd_{\bQ}(N)}\,\,\bigvee_{r\in
{\cal I} \cap N} r:\ \ p \in {\bQ}^N \mbox{ and }N \mbox{ is a
}\bQ\mbox{--candidate} \big\}\]
[we are assuming no incidental identification] and, in any reasonable way,
code them, if $\bQ$ is simple, as members of ${\cal H}_{< \aleph_1}(\theta)$,
for snep (or very simple) $\bQ$ we work slightly more to code them as members
of ${}^\omega \theta$, pedantically easier in ${}^\omega(\theta + \omega)$),
\item[(b)] the order $\le^{\bQ'}$ is given by $q_1 \le^{\bQ'} q_2$ if and
only if one of the following occurs:
\begin{enumerate}
\item[$(\alpha)$] $q_1,q_2 \in \bQ$, $q_1 \le^{\bQ} q_2$,
\item[$(\beta)$]  $q_1 \in \bQ$, $q_2 = p\ \&\ \biggl(\bigwedge\limits_{{\cal
I} \in\pd_Q(N)}\, \bigvee\limits_{r\in {\cal I}\cap N} r \biggr)$ and $q_1
\le^{\bQ} p$,
\item[$(\gamma)$] $q_1 = p\ \&\ \biggl(\bigwedge\limits_{{\cal I}\in\pd(N,Q)}
\,\bigvee\limits_{r\in {\cal I}\cap N} r\biggr)$ and $q_2\in\bQ$, $p \le^{\bQ}
q_2$ and if ${\cal I}\in\pd(N,\bQ)$ then 
\[\begin{array}{r}
(\exists q' \in \bQ)(\exists\langle p_n:n\in\omega\rangle)(q'\le^{\bQ} q_2\
\&\ \varphi^\bQ_2(\ldots, p_n, \ldots, q')\ \&\quad \\
\{p_n:n<\omega\}\mbox{ lists }{\cal I}\cap N)
  \end{array}\]
\item[$(\delta)$] $q_\ell=p_\ell \ \&\ \biggl(\bigwedge\limits_{{\cal I}\in
\pd(N,Q)} \,\bigvee\limits_{r\in {\cal I}\cap N} r\biggr)$ (for $\ell =1,2$)
and: $q_1 = q_2$ {\em or} $q_1 \le p_2$ by clause $(\gamma)$.
\end{enumerate}
\end{enumerate}
\end{definition}

\noindent{\bf Remark:}\quad In \cite{Sh:f}, for a hereditarily countable name,
instead of 
\[p\ \&\ \bigwedge\limits_{{\cal I} \in \pd_{\bQ}(N)}\,\bigvee\limits_{r\in
  {\cal I}\cap N} r\]
we use the first member of $\bQ_i$ which forces this. Simpler, but when we ask
whether this guy is $\le q$ (for some $q \in\bQ$) we run into uncountable
antichains.  

\begin{proposition}
\label{4.2}
\begin{enumerate}
\item Assume $\bQ$ is explicitly nep. {\em Then}:
\begin{enumerate}
\item[(a)] in Definition \ref{4.1}, $\bQ'$ is a (quasi) order,
\item[(b)] $\le^{\bQ'}\rest\bQ=\le^{\bQ}$,
\item[(c)] $\bQ$ is a dense subset of $\bQ'$.
\end{enumerate}
\item Assume in addition:
\begin{enumerate}
\item[$(\boxtimes_2)$] $\bQ$ is explicitly nep in every $\bQ$--candidate. 
\end{enumerate}
{\em Then}:
\begin{enumerate}
\item[(d)] if $N$ is a $\bQ$--candidate, $N\models$``$p\in\bQ'$'', then for
some $q\in N$ we have $N\models$``$p\le^{\bQ'}q \ \&\ q \in \bQ$'',
\item[(e)] $\bQ'$ is explicitly nep (with the same ${\frak B}^{\bQ}$ and
parameters). 
\end{enumerate}
\item Assume in addition
\begin{enumerate}
\item[$(\boxtimes_3)$] for any $\bQ$--candidate $N$, if $N'$ is a generic 
extension of $N$ for the forcing notion $\Levy(\aleph_0,|{\cal P}(\bQ)|^N)$,
{\em then} $N'$ is a $\bQ$--candidate.
\end{enumerate}
{\em Then} we can add
\begin{enumerate}
\item[(e)$^+$]  $\bQ'$ is explicitly local nep (see Definition \ref{0.5A}).
\end{enumerate}
\item  We can replace above (in the assumption and conclusion) nep by snep, or
nep by simple nep.
\end{enumerate}
\end{proposition}

\begin{remark}
\label{4.2A}
The definition of ``local'' (in \ref{0.5A}) and the statement $(\boxtimes_3)$
in \ref{4.2}(3) can be handled a little differently. We can (in \ref{0.5A}(2))
demand less on $N'$ (it is not a $\bQ$--candidate), just have some of its
main properties and in $\boxtimes_3$ of \ref{4.2}(3), $\ZFC^-_*$ says that
${\cal H}(\theta)$ is a set (so has a cardinality) and is a $\bQ$--candidate. 
So we may consider having $\ZFC^-_\ell$ for several $\ell$'s, $\ZFC^*_\ell$
speaks on $\chi_0>\ldots>\chi_{\ell -1}$ and the generic extensions of a model
of $\ZFC^*_{\ell +1}$ for $\Levy(\aleph_0,\chi_\ell)$ is a model of
$\ZFC^-_\ell$.  Similar remarks hold for \S7. But, as we can deal with the
nice case (see Definition \ref{2.3B}), we may start with a countable $N\prec({
\cal H}(\beth_\omega),\in)$ (or even better $({\cal H}(\beth_{\omega_1}),\in)$
so that ``countable depth can be absorbed''), we ignore this in our main
presentation.  

Does $(\boxtimes_3)$ of \ref{4.2}(3) occur at all? Let $G$ be a subset of
$\Levy(\aleph_0,|{\cal P}(\bQ)|^N)$ generic over $N$. Then $N'\stackrel{\rm
def}{=}N\langle G\rangle$ is a $\bQ$--candidate.
\end{remark}

{\noindent{\sc Proof of \ref{4.2}} \hspace{0.2in}} 1) {\bf Clause (a)}:\quad
Assume $q_1 \le q_2 \le q_3$; we have $2^3=8$ cases according to truth values
of $q_i\in \bQ$: 
\medskip

\noindent{\sc Case (A)}:\qquad  $q_1,q_2,q_3\in \bQ$.\\
Trivial.
\medskip

\noindent{\sc Case (B)}:\qquad  $q_1,q_2\in \bQ$, $q_3\notin\bQ$.\\
Check.
\medskip

\noindent{\sc Case (C)}:\qquad  $q_1\notin\bQ$, $q_2,q_3\in\bQ$.\\
Check.
\medskip

\noindent{\sc Case (D)}:\qquad  $q_1\in\bQ$, $q_2\notin\bQ$, $q_3\in\bQ$.\\ 
Then $q_2=p_2 \ \&\ \bigwedge\limits_{{\cal I} \in\pd(N,\bQ)}\,
\bigvee\limits_{r\in {\cal I}\cap N} r$ and $q_1\le^{\bQ} p_2$ (by
\ref{4.1}(b)$(\beta)$) and $p_2\le^{\bQ} q_3$ (by \ref{4.1}(b)$(\gamma)$). 
Hence $q_1 \le^{\bQ} q_3$ follows.
\medskip

\noindent{\sc Case (E)}:\qquad  $q_1\in\bQ$, $q_2\notin\bQ$, $q_3\notin\bQ$.\\
Let $q_\ell=p_\ell\ \&\ \bigwedge\limits_{{\cal I}\in\pd(N,Q)}\,
\bigvee\limits_{r\in {\cal I}\cap N} r$ for $\ell = 2,3$. So $q_1\le^{\bQ}
p_2$ (see \ref{4.1}(b)$(\beta)$) and $p_2\le^{\bQ} p_3$ (see
\ref{4.1}(b)$(\gamma),(\delta)$). Hence $q_1\le^{\bQ} p_2$ (as $\le^{\bQ}$ is
transitive) and so $q_1\le q_3$ (see \ref{4.1}(b)$(\beta)$).
\medskip

\noindent{\sc Case (F)}:\qquad  $q_1\notin\bQ$, $q_2\notin\bQ$, $q_3\in\bQ$.\\
Let $q_\ell=p_\ell \ \&\ \bigwedge\limits_{{\cal I}\in\pd(N,\bQ)}\,
\bigvee\limits_{r\in {\cal I}\cap N} r$ for $\ell = 1,2$ and suppose
that $q_1\neq q_2$ (otherwise trivial). Then, by \ref{4.1}(b)$(\delta)$, $q_1
\le p_2$ and by \ref{4.1}(b)$(\gamma)$, $p_2\le q_3$ so by the previous case
(C), $q_1 \le q_3$ as required. 
\medskip

\noindent{\sc Case (G)}:\qquad  $q_1\notin\bQ$, $q_2\in\bQ$, $q_3\notin\bQ$.\\
Let $q_\ell=p_\ell\ \&\ \bigwedge\limits_{{\cal I}\in\pd(N,\bQ)}\,
\bigvee\limits_{r\in{\cal I}\cap N} r$ for $\ell=1,3$. Now, by \ref{4.1}(b)$(
\beta)$, $q_2\le p_3$ and by the previous case (C), $q_1 \le p_3$ and hence, by
\ref{4.1}(b)$(\delta)$, $q_1 \le q_3$ as required. 
\medskip

\noindent{\sc Case (H)}:\qquad  $\bigwedge\limits_\ell q_\ell\notin \bQ$.\\
Let $q_\ell=p_\ell\ \&\ \bigwedge\limits_{{\cal I}\in \pd(N,\bQ)}\,
\bigvee\limits_{r\in {\cal I}\cap N} r$. If $q_1=q_2$ or $q_2=q_3$ then the
conclusion is totally trivial. So assume not. Thus 
\[\begin{array}{ll}
q_1\le p_2&\quad\mbox{(by clause }(\delta)\mbox{ a case defined in }(\gamma))\\
q_2\le p_3&\quad\mbox{(by clause }(\delta)).
  \end{array}\]
Hence $p_2\le p_3$ (see clause $(\gamma)$), so ``a previous case'' applies. 
This finishes the proof of the clause (a).
\medskip

\noindent{\bf Clause (b)}:\qquad Totally trivial.
\medskip

\noindent{\bf Clause (c)}:\qquad Let $q\in\bQ'$; if $q\in\bQ$ then there is
nothing to do; otherwise for some $\bQ$--candidate $N$ we have $q=p\ \&\
\bigwedge\limits_{{\cal I}\in\pd(N,\bQ)}\,\bigvee\limits_{r\in {\cal I}\cap N}
r$ and use nep (i.e.\ clause (c) of \ref{0.2}(1)) on the $\bQ$--candidate $N$.
\medskip

\noindent 2) Assume $(\boxtimes_2)$.

\noindent{\bf Clause (d)}:\qquad Proved inside the proof of clause (e).
\medskip

\noindent{\bf Clause (e)}:\qquad More pedantically we have to define
\[\varphi^{\bQ'}_0,\varphi^{\bQ'}_1,\varphi^{\bQ'}_2,{\frak B}^{\bQ'},
\theta^{\bQ'}\]
and then prove the required demands for a $\bQ'$--candidates. We let ${\frak
B}^{\bQ'}={\frak B}^{\bQ}$, $\theta^{\bQ'}=\theta^{\bQ}$, the formulas will be
different, but with the same parameters. So the $\bQ'$--candidates are the
$\bQ$--candidates. What is $\varphi^{\bQ'}_0$? It is
\[\begin{array}{r}
\varphi^{\bQ}_0(x) \vee\mbox{`` }x \mbox{ has the form }\quad p \ \&\  
\bigwedge\limits_{{\cal I}\in\pd_{\bQ(M)}}\,\bigvee\limits_{r\in{\cal I}\cap M}
r,\qquad \mbox{where}\\
M \mbox{ is a }\bQ \mbox{-candidate (so countable) and }\varphi^{\bQ}_0(p)
\mbox{ ''}.
  \end{array}\]
Clearly $\varphi^{\bQ'}_0$ defines $\bQ'$ through $\bQ'$--candidates. Note
that if $N$ is a $\bQ'$--candidate and $N\models$``$M$ is a countable
$\bQ$--candidate'', then we have $M\subseteq N$, and if $M\models$``$x$ is
countable'', then $x\subseteq M\subseteq N$; so $M$ is really a
$\bQ$--candidate. Consequently, $\varphi^{\bQ'}_0$ is upward absolute for
$\bQ'$-candidates and it defines $\bQ'$. So clause (a) of Definition
\ref{0.2}(1) holds.

Now we pay our debt proving clause (d). Let $N$ be a $\bQ'$--candidate and $N
\models$``$p\in\bQ'$'', i.e.\ $N \models\varphi^{\bQ'}_0(p)$. By the
definition of $\bQ'$, either $N\models$``$p\in\bQ$'' and we are done, or for
some $p',M\in N$ we have 
\[N\models\mbox{``}M\mbox{ is a }\bQ' \mbox{--candidate, $p'\in\bQ^M$, and }
p= \bigl(p'\ \&\ \bigwedge\limits_{{\cal I}\in\pd(M,\bQ)}\,\bigvee\limits_{r
\in{\cal I}\cap M} r\bigr)\mbox{''}.\]
By the assumption $(\boxtimes_2)$, for some $q\in\bQ^N$ we have $N \models$``
$q$ is explicitly $\langle M,\bQ\rangle$--generic'' and $N\models$``$p'\le^{
\bQ} q$''. Then for some $\langle\langle r_{{\cal I},\ell}:\ell<\omega\rangle:
{\cal I}\in \pd(M,\bQ)\rangle\in N$ we have:\quad $N\models$``$\{r_{{\cal I},
\ell}:\ell<\omega\}$ enumerates ${\cal I}\cap M$'' and $N\models$``$\varphi^{
\bQ}_2(r_{{\cal I}_\ell,0},r_{{\cal I},1},\ldots,q)$''.  Now it follows from
the definition of $\bQ'$ that $N\models$``$p\le^{\bQ'}q$'', so $q$ is as
required. 

What is $\varphi^{\bQ'}_1$? Just write the definition of $p\le^{\bQ'}q$ from
clause (b) of \ref{4.1}. Clearly also $\varphi^{\bQ'}_1$ is upward absolute for
$\bQ'$--candidates and it defines the partial order of $\bQ'$ (even in
$\bQ'$--candidates). So clause (b) of Definition \ref{0.2}(1) holds.

What is $\varphi^{\bQ'}_2$?  Let it be:
\begin{enumerate}
\item[\ ]$\varphi^{\bQ'}_2(p_0,p_1,\ldots,p_\omega)\stackrel{\rm def}{=}$\\
``there are $M,p,q$ such that:\quad $M$ is a $\bQ'$--candidate and $p\in\bQ^M$
and $q = \bigl(p\ \&\ \bigwedge\limits_{{\cal I} \in\pd(M,\bQ)}\,
\bigvee\limits_{r\in {\cal I}\cap M} r\bigr)$ and $q\le^{\bQ'} p_\omega$ and
for some ${\cal J}\in\pd(M,\bQ)$, if $r\in {\cal J}\cap M$ then there is
$\ell$ such that $p_\ell\le^{\bQ'} r$''. 
\end{enumerate}
To show that $\varphi^{\bQ'}_2$ is upward absolute for $\bQ'$--candidates
suppose that $N$ is a $\bQ'$--candidate and $N\models\varphi^{\bQ'}_2(p_0,
p_1,\ldots,p_\omega)$ and let $M,p,q$ witness it. Then, in $N$, $M$ is a
$\bQ'$--candidate, so $p\in\bQ$, $q\in\bQ'$ and for some ${\cal J}\in\pd(M,
\bQ)$ we have:
\begin{quotation}
if $r\in{\cal J}\cap M$, then there is $\ell$ such that $p_\ell\le^{\bQ'}r$.
\end{quotation}
By the known upward absoluteness all those statements hold in $\V$, too. 
Assume now that $\varphi^{\bQ'}_2(p_0,p_1,\ldots,p_\omega)$ holds as witnessed
by $M,p,q$ and ${\cal J}\in \pd(M,\bQ)$. Suppose $q'\geq p_\omega$ and we may
assume that $q'\in\bQ$ (by (1)(c)). Then $q\leq q'$ and (by clause $(\gamma)$
of the definition of $\leq^{\bQ'}$) we have $q''\in\bQ$ such that $q''\leq q'$
and $\varphi^{\bQ}_2(r_0,r_1,\ldots,q')$ for some list $\{r_n:n<\omega\}$ of
${\cal J}\cap M$. Thus ${\cal J}\cap M$ is predense (in $\bQ$) above $q''$ and
we find $r\in {\cal J}\cap M$ such that $r,q'$ are compatible. But now, there
is $\ell<\omega$ such that $p_\ell\leq^{\bQ'} r$, so necessarily $p_\ell,q'$
are compatible (in $\bQ'$). This shows \ref{0.2}(2)(b)$^+$. Let us turn to
clause (c)$^+$ of Definition \ref{0.2}(2). So suppose that $N$ is a
$\bQ'$--candidate and $p\in\bQ'\cap N$. By clause (d), there is $p'$ such that
$N\models$``$p\le^{\bQ'}p' \ \&\  p'\in\bQ$''. Let $q=p'\ \&\
\bigwedge\limits_{{\cal I} \in\pd(N,\bQ)}\,\bigvee\limits_{r\in {\cal I}\cap
N} r$, clearly $q \in\bQ'$ and $\bQ'\models$``$p'\le q$''. Hence, by
\ref{4.2}(1)(a), we know $\bQ'\models$``$p \le q$''. For ${\cal J}\in\pd(N,
\bQ')$ let 
\[{\cal J}'=\{q\in\bQ^N:N\models\mbox{``} q \mbox{ is above some member of
${\cal J}$ in }\bQ'\mbox{ ''}\}.\]
Note that if ${\cal J}\in\pd(N,\bQ')$ then ${\cal J}'\in\pd(N,\bQ)$, and so
${\cal J}'\cap N$ is predense above $q$. Moreover, $(\forall r\in {\cal J}
\cap N)(\exists r'\in {\cal J}'\cap N)(r\le^{\bQ'} r')$. So let ${\cal
J}\in\pd(N,\bQ')$ and let $\langle p_n:n<\omega\rangle$ be an enumeration of
${\cal J}\cap N$. It should be clear that $\varphi^{\bQ'}_2(p_0,p_1,\ldots, q)$
holds as witnessed by $N,p',q$ and ${\cal I}'$.
\medskip

\noindent 3) Compared to (e) of \ref{4.2}(2) we have also to prove (e)$^+$,
i.e.\ strengthen the clause (c)$^+$ of Definition \ref{0.2}(1) by $(*)$ of
Definition \ref{0.5A}(2). 

Let $N^+$ be a generic extension of a $\bQ'$--candidate $N$ by the forcing
notion $\Levy(\aleph_0,|{\cal P}(\bQ)|^N)$. Clearly for every $p \in Q^N$, the
condition 
\[p\ \&\ \bigwedge\limits_{{\cal I}\in \pd(N,\bQ)}\,\bigvee\limits_{r\in {\cal
I}\cap N} r\]
belongs to $N^+$. So by the proof of clause (c)$^+$ of Definition \ref{0.2}(1)
in the proof of (e) above, we are done.  \medskip

\noindent 4)  Similar proof.  \hfill$\square_{\ref{4.2}}$

\begin{discussion}
If we would like not to use \ref{4.2}, we may like to try the following
Definition \ref{4.3}. Note that there: $\cl_1(\bQ)$ cannot serve as a forcing
notion as it contains ``false'', $\cl_2(\bQ)$ is the reasonable restriction,
and $\cl_3(\bQ)$ has the same elements but more ``explicit'' quasi order. We
do not define a quasi order on $\cl_1(\bQ)$, but it is natural to use the one
of $\cl_2(\bQ)$ adding:\quad $\psi\leq\varphi$ if $\varphi\in\cl_1(\bQ)
\setminus\cl_2(\bQ)$. No harm in allowing in the definition of $\cl_1(\bQ)$
also $\neg$ (the negation). The previous $\cl(\bQ)$ is close to $\cl_3(\bQ)$. 
\end{discussion}

\begin{definition}
\label{4.3}
Let $\bQ$ be a forcing notion. 
\begin{enumerate}
\item Let $\cl_1(\bQ)$ be the closure of the set $\bQ$ by conjunctions and
disjunctions over sequences of members of length $\le \omega$ [we may add: and 
$\neg$ (the negation)]; wlog there are no incidental identification and $\bQ
\subseteq\cl_1(\bQ)$.  
\item  For a generic $G\subseteq\bQ$ over $\V$ and $\psi\in\cl_1(\bQ)$ let
$\psi[G]$ be the truth value of $\psi$ under $G$ where for $\psi=p\in\bQ$,
$\psi[G]$ is the truth value of $p\in G$. (We will use $\truth$ for
``truth''.)  
\item $\hat{\bQ}=\cl_2(\bQ)=\{\psi\in\cl_1(\bQ):\mbox{ for some } p\in\bQ
\mbox{ we have } p\forces$``$\psi[\name{G}_{\bQ}]=\truth$''$\}$, ordered by:
\[\psi_1 \le^{\hat{\bQ}}\psi_2\ \Leftrightarrow\ (\forall p \in \bQ)
[p \forces_{\bQ}\mbox{``}\psi_2[\name{G}_{\bQ}]=\truth\mbox{''}\ \
\Rightarrow\ \ p\forces_{\bQ}\mbox ``\psi_1[\name{G}_{\bQ}]=\truth\mbox{''}].\]
\item Let $\bQ$ be explicitly nep. We let $\cl_3(\bQ)$ be the following
forcing notion:
\begin{enumerate}
\item[(a)] the set of elements is $\cl_2(\bQ)$,
\item[(b)] the order $\leq^{\hat{\bQ}}_3=\leq^\bQ_3=\leq^{\cl_3(\bQ)}_3=\leq_{
\cl_3(\bQ)}$ is the transitive closure of $\le^{\hat{\bQ}}_0$ which is defined
by 

$\psi_1\le^{\hat{\bQ}}_0\psi_2$\quad{\em iff}\quad one of the following occurs
\begin{enumerate}
\item[(i)]   $\psi_1,\psi_2\in\bQ$ and $\psi_1\le^{\bQ} \psi_2$,
\item[(ii)]  $\psi_1$ is a conjunct of $\psi_2$ (meaning: $\psi_1=\psi_2$ or
$\psi_2=\bigwedge\limits_{n<\alpha}\psi_{2,n}$, and $\psi_1\in\{\psi_{2,n}:
n<\alpha\}$),  
\item[(iii)] $\psi_2\in\bQ$ and there is a $\bQ$-candidate $M$ such that $p,
\psi_1\in M$, $p\in \bQ^M$, $p\le^{\bQ}\psi_2$, $\psi_2$ is explicitly
$\langle M,\bQ\rangle$--generic and $M\models$``$p\forces\psi_1[\name{G}_{\bQ}
]=\truth$'' and if $q\in\bQ$ is a conjunct of $\psi_1$ then $M
\models$``$q\le^{\bQ} p$''. 
\end{enumerate}
\end{enumerate}
\end{enumerate}
\end{definition}

\begin{proposition}
\label{4.4}
\begin{enumerate}
\item $\bQ\subseteq\hat{\bQ}$, $\leq^{\hat{\bQ}}$ is a quasi order, and
$\le^{\hat{\bQ}}\rest\bQ =\{(p,q):q\forces_{\bQ}$``$p\in\name{G}_{\bQ}$''$\}$,
so if $\bQ$ is separative then $\le^{\hat{\bQ}}\rest \bQ=\le^{\bQ}$; and $\bQ$
is a dense subset of $\hat{\bQ}$. 
\item Assume $\bQ$ is temporarily explicitly nep.  Then:
\begin{enumerate}
\item[(a)]  $\bQ\subseteq\cl_3(\bQ)$ and $\le^{\bQ}_3\rest\bQ\supseteq
\le^{\bQ}$ and $\leq^{\bQ}_3\subseteq\leq^{\hat{\bQ}}$,
\item[(b)]  $\bQ$ is a dense subset of $\cl_3(\bQ)$.
\end{enumerate}
\item Assume in addition 
\begin{enumerate}
\item[$(\circledast_3)$] $\bQ$ is correctly explicitly nep in $\V$ and in
every $\bQ$--candidate.
\end{enumerate}
Then
\begin{enumerate}
\item[(d)] if $N$ is a $\bQ$--candidate and $N\models$``$p\in\cl_3(\bQ)$'' 

{\em then} for some $q\in N$ we have $N\models$``$p\leq_{\cl_3(\bQ)} q\ \&\
q\in\bQ$'',
\item[(e)] $\bQ'$ is explicitly nep and correct.
\end{enumerate}
\item Assume in addition
\begin{enumerate}
\item[$(\circledast_4)$] for any $\bQ$--candidate $N$, if $N'$ is a generic 
extension of $N$ for the forcing notion $\Levy(\aleph_0,|{\cal P}(\bQ)|^N)$,
{\em then} $N'$ is a $\bQ$--candidate.
\end{enumerate}
{\em Then} we can add
\begin{enumerate}
\item[(e)$^+$]  $\cl_3(\bQ)$ is explicitly local nep (see Definition
\ref{0.5A}). 
\end{enumerate}
\end{enumerate}
\end{proposition}

\Proof Straight, e.g.

\noindent {\bf (2) Clause (b):}\quad Assume $\psi\in\cl_3(\bQ)$, so $\psi\in
\cl_2(\bQ)$ and for some $p\in\bQ$ we have $p\forces_{\bQ}$``$\psi[\name{G}_{
\bQ}]=\truth$''. There is a $\bQ$--candidate $M$ to which $p$ and $\psi$
belong (as $\ZFC^-_*$ is $\emptyset$--good). Let $q$ be explicitly
$(M,\bQ)$--generic, and $\bQ\models p\leq q$. So, by clause (iii) of
\ref{4.3}(4)(b), we have $\cl_3(\bQ)\models$``$\psi\leq q''$, as required.

\noindent {\bf (3) Clause (e):}\quad Let $\varphi^3_0(x)$ say that there is a
$\bQ$--candidate $M$ such that $M\models$``$x\in\cl_3(\bQ)$''. Let 
$\varphi^3_1(x,y)$ say the definition of $\le^\theta_3$. Lastly, $\varphi^3_2
(\langle x_i:i \le\omega\rangle)$ says that for some $\langle y_i:i \le\omega
\rangle$ we have 
\[\varphi^{\bQ}_2(\langle y_i:i \le\omega\rangle),\quad y_\omega\le^{\bQ}_3
x_\omega\quad\mbox{(i.e.\ $\varphi^3_1(y_\omega,x_\omega)$) and}\quad
\bigwedge_{i<\omega}\,\bigvee_{j<\omega}x_j \le^{\bQ}_3 y_i.\]
\QED

\begin{remark}
\label{4.5}
Instead of using $\cl(\bQ)$ from \ref{4.1} below we can have in
$\bar{\varphi}$, a function which from an $\omega$--list of the elements of $N$
and from $p$ computes an element of $\bQ$ having the role of $p\ \&\
\bigwedge\limits_{{\cal I}\in\pd(N,\bQ)}\,\bigvee\limits_{r\in{\cal I}\cap N}
r$. The choice does not seem to matter. 
\end{remark}

\begin{definition}
\label{4.6}
For a forcing notion $\bP$ and a cardinal (or ordinal) $\kappa$, we define
what is an $\hc$-$\kappa$-$\bP$--name (here $\hc$ stands for hereditarily
countable), and for this we define by induction on $\zeta<\omega_1$ what is
such a name of depth $\le \zeta$. 
\medskip

\noindent{$\zeta = 0$}:\qquad It is $\alpha$, that is $\check{\alpha}$, for
some $\alpha < \kappa$.

\noindent{$\zeta > 0$}:\qquad It has the form $\name{\tau}=\{\langle p_i,
\name{\tau}_i\rangle:i<i^*\}$,  where $i^* < \omega_1$, $p_i\in\cl_1(\bP)$
from Definition \ref{4.3}(1) and $\name{\tau}_i$ an $\hc$-$\kappa$-$\bP$--name
of some depth $<\zeta$; that is for $G\subseteq \bP$ generic over $\V$, we let
$\name{\tau}[G]=\{\name{\tau}_i[G]:p_i[G]=\truth\}$.
\medskip

\noindent An $\hc$-$\kappa$-$\bP$--name is an $\hc$-$\kappa$-$\bP$--name of
some depth $< \omega_1$. An $\hc$-$\kappa$-$\bP$--name $\name{\tau}$ has depth
$\zeta$ if it has depth $\le \zeta$, but not $\le \xi$ for $\xi < \zeta$.
\end{definition}

\noindent{\bf Remark:}\quad Why did we use $p\in\cl_1(\bQ)$ and not $p\in
\cl_3(\bQ)$? As the membership in $\cl_1(\bQ)$ is easier to define.

\begin{proposition}
\label{4.7}
\begin{enumerate}
\item If $\name\tau$ is an $hc$-$\kappa$-$\bP$--name and $G\subseteq\bP$ is
generic over $\V$ {\em then} $\name{\tau}[G]\in {\cal H}_{< \aleph_1}(
\kappa)$. If in addition $\bP\subseteq {\cal H}_{< \aleph_1}(\kappa)$ then
$\name{\tau}\in {\cal H}_{<\aleph_1}(\kappa)$. 
\item Let $\varphi(x_0,\ldots,x_{n-1})$ be a first order formula and
$\name{\tau}_0,\ldots,\name{\tau}_{n-1}$ be $\hc$-$\kappa$-$\bP$--names. Then
there is $p\in\cl_1(\bP)$ such that for every $G\subseteq\bP$ generic over
$\V$: 
\[\bigl(\bigcup_{\ell < n} {\rm Tc}^{\rm ord}(\name{\tau}_\ell[G]),\in\bigr)
\models\varphi(\name{\tau}_0[G],\ldots,\name{\tau}_{n-1}[G])\quad\mbox{\rm
iff }\quad p[G] =\truth.\]
\item The set of $\hc$-$\kappa$-$\bP$--names is closed under the following
operations: 
\begin{enumerate}
\item[(a)] difference,
\item[(b)] union and intersection of two, finitely many and even countably
many, 
\item[(c)] definition by cases:\quad for $p_n\in\cl_1(\bP)$ and
$\hc$-$\kappa$-$\bP$--names $\name{\tau}_n$ (for $n<\omega$) there is a
$\hc$-$\kappa$-$\bP$--name $\name{\tau}$ such that for a generic $G\subseteq
\bQ$ over $\V$ we have 
\[\name{\tau}[G]\mbox{ is}:\quad\begin{array}{ll}
\name{\tau_n}[G] &\mbox{if } p_n[G]=\truth\ \&\ \bigwedge\limits_{\ell<n} \neg p_\ell[G]=\truth\\
\emptyset &\mbox{if }\bigwedge\limits_{\ell<\omega}\neg p_\ell[G]=\truth.
                            \end{array}\]
\end{enumerate}
\end{enumerate}
\end{proposition}

\begin{definition}
\label{4.8}
\begin{enumerate}
\item A forcing notion $\bQ$ (or $\bar{\varphi}$) is temporarily, explicitly
straight $(\kappa,\theta)$--nep for ${\frak B}$ if: old conditions from
Definition \ref{0.2}(1),(2) (for explicitly $(\kappa,\theta)$--nep) but
possibly ${\frak B} \subseteq {\cal H}_{< \aleph_1}(\kappa)$; and
\begin{enumerate}
\item[(d)] $\bQ \subseteq {\cal H}_{<\aleph_1}(\theta)$ (i.e.\ $\bQ$ is
simple) and $\aleph_1+\theta\le \kappa$,
\item[(e)] for $\ell<3$ the formula $\varphi^{\bQ}_\ell(\bar{x})$ is of the
form 
\[(\exists t)[t\in {\cal H}_{<\aleph_1}(\kappa)\ \&\ t={\rm Tc}^{\rm ord}(t)\
\&\ t\cap\omega_1\mbox{ is an ordinal}\ \&\ \psi^{\bQ}_\ell(\bar{x},t)],\] 
where in the formula $\psi^{\bQ}_i$ the quantifiers are of the form $(\exists
s \in t)$ and the atomic formulas are ``$x \in y$'',``$x$ is an ordinal'' and
those of ${\frak B}^{\bQ}$.
\end{enumerate}
\item In clause (e) of part (1), we call such $t$ an explicit witness for
$\varphi^{\bQ}_i(\bar{x})$. We call $t$ a weak witness, if for every
$\bQ$--candidate $N$, $\bar{x}\in N$, if $t\in N$ then $N\models\varphi^{
\bQ}_\ell(\bar{x})$.  We call it a witness if:
\begin{enumerate}
\item[(i)]   $\ell=0$ and it is an explicit witness, or
\item[(ii)]  $\ell=1$ (so $\bar{x}=\langle x_0,x_1\rangle$) and $t$ gives
$k$, $y_0,\ldots,y_k$, $t_0,\ldots,t_{k-1}$, $s_0,\ldots,s_k$ such that:\quad
$s_\ell$ explicitly witnesses $\varphi_0(y_\ell)$, $t_\ell$ explicitly
witnesses $y_\ell\le^{\bQ}y_{\ell+1}$ and $y_0=x_0$, $y_k=x_1$ (so $y_\ell\in
t$, $s_\ell \in t$, $x_\ell \in t$),
\item[(iii)] $\ell=2$ (so $\bar{x}=\langle x_i:i\le\omega\rangle$) and $t$
gives $\langle y_i:i \le\omega\rangle$, $\langle k_i:i<\omega\rangle$,
$\langle s_i:i\le\omega+1\rangle$ such that $s_\omega$ is a witness to
$y_\omega \le x_\omega$, $s_{\omega +1}$ is an explicit witness to
$\varphi^{\bQ}_2(\langle y_i:i \le\omega\rangle)$, $s_i$ is a witness to
$x_i \le^{\bQ} y_{k_i}$ (so also they all belong to $t$, as well as witnesses
to $x_i,y_j \in \bQ$).
\end{enumerate}
\end{enumerate}
\end{definition}

\begin{proposition}
\label{4.8A}
\begin{enumerate}
\item Assume $\bQ$ is temporarily explicitly straight $(\kappa,\theta)$--nep
for ${\frak B}$. {\em Then} $\bQ$ is temporarily simply explicitly
$(\kappa,\theta)$--nep for ${\frak B}$. Sufficient conditions for
``$K$--absolutely'' are as in \S2.
\item Assume $\bQ$ is temporarily correctly simply explicitly $(\kappa,
\theta)$--nep for ${\frak B}$ and $\theta+\aleph_1\leq\kappa$. {\em Then}
$\bQ$ is temporarily straight explicitly $(\kappa,\theta)$--nep for ${\frak
B}$ and is correct.
\end{enumerate}
\noindent {\em [Nevertheless, ``simple'' and ``straight'' are distinct as
properties of $({\frak B},\bar{\varphi},\theta)$, i.e.\ the point is changing
$\bar{\varphi}$.]} 
\end{proposition}

\Proof Straight. \QED

\begin{defthe}
\label{4.9}
By induction on $\alpha$ we define and prove the following situations:
\begin{enumerate}
\item[(A)]  {\em [Definition]}\quad $\bar{\bQ}=\langle(\bP_i,\nbQ_i,\bar{
\varphi}_i,\name{{\frak B}}_i,r_i,\theta_i):i<\alpha\rangle$ is
nep--CS--iteration.
\item[(B)]  {\em [Definition]}\quad $\kappa^{\bar{\bQ}}=\kappa[\bar{\bQ}]$, in
short $\kappa^\alpha$ abusing notation.
\item[(C)]  {\em [Definition]}\quad We define ${\frak B}^\alpha={\frak B}^{
\bar{\bQ}}$.
\item[(D)]  {\em [Definition]}\quad $\Lim(\bar{\bQ})=\bP_\alpha$ and $\bP_{
\alpha,w}$ for any set $w$ of ordinals $<\alpha$ for $\bar{\bQ}$ as above.
\item[(E)] {\em [Claim]}\quad If $\bar{\bQ}$ is a nep--CS--iteration, and
$\alpha=\lh(\bar{\bQ})$, {\em then} $\bar{\bQ}\rest\beta$ is a
nep--CS--iteration (for $\beta<\alpha$), $\Lim(\bar{\bQ}\rest\beta)=\bP_\beta$
and $\bP_\beta\subseteq {\cal H}_{< \aleph_1}(\kappa_\beta)$.
\item[(F)]  {\em [Claim]}\quad For $\bar{\bQ}$ as in (A), a ${\frak
B}^\alpha$--candidate $N$, $\gamma\le\beta\le\alpha$ and $p,q\in\bP_\beta$: 
\begin{enumerate}
\item[(a)] $p$ is a function with domain a countable subset of $\beta$
(pedantically see clause (D)),
\item[(b)] $\bP_\beta$ is a forcing notion (i.e.\ a quasi order) satisfying
(d) + (e) of \ref{4.8} and (a), (b), (b)$^+$ of \ref{0.2}(1),(2), 
\item[(c)] $p\rest\gamma\in\bP_\gamma$ and $\bP_\beta\models$``$p\rest\gamma
\le p$'', 
\item[(d)] $\bP_\gamma\models$``$p\rest\gamma\le q$'' implies $\bP_\beta
\models$``$p\le (q\cup p \rest [\gamma,\beta))$'', 
\item[(e)] $\bP_\gamma\subseteq\bP_\beta$ and even $P_\gamma\lesdot\bP_\beta$,
\item[(f)] $p\in\bP_\beta$\quad iff\quad $p$ a function with domain $\in [\beta]^{\le \aleph_0}$ and 
\[\zeta\in\dom(p)\quad \Rightarrow\quad p\rest\zeta\in\bP_{\zeta+1}.\]
\end{enumerate}
\item[(G)] {\em [Definition]}\quad For a ${\frak B}^\alpha$--candidate $N$ and
$w,\beta,\gamma$ such that $N\models$``$w\subseteq\alpha$'', $\gamma<\beta \le
\alpha$ and $\beta,\gamma\in (w \cup \{\alpha\})\cap N$, and $q\in\bP_\beta$,
$p\in N$ such that $N\models$``$p \in\bP_\beta$'' and $q\rest\gamma$ is $(N,
\bP_\gamma)$--generic we define when $q$ is $[\beta,\gamma)$--canonically
$(N,P_\beta,w)$--generic above $p$. 
\item[(H)] {\em [Theorem]}\quad If $q\in N$ is a $[\beta,\gamma)$--canonically
$(N,\bP_\beta,w)$--generic above $p$, {\em then} $q$ is $(N,
\bP_\beta)$--generic and $p \le q$.
\item[(I)] {\em  [Theorem]}\quad $\bP_\alpha$ is explicitly straight correct
$\kappa^\alpha$--nep for $\bar{\varphi},{\frak B}^\alpha$. 
\item[(J)] {\em [Theorem]}\quad  For any $\kappa\ge\kappa^\alpha$,
\[\forces_{\bP_\alpha}\mbox{`` }({\cal H}_{<\aleph_1}(\kappa))^{\V[\bP_\alpha
]}=\{\name{\tau}[\name{G}_{\bP_\alpha}]:\name{\tau}\mbox{ is an
$\hc$--$\kappa$--$\bP_\alpha$--name }\})\mbox{''}.\] 
\end{enumerate}
\end{defthe}
Let us carry out the clauses one by one.
\medskip

\noindent{\sc Clause (A)}, Definition:\qquad  $\bar{\bQ}=\langle(\bP_i,
\nbQ_i,\bar{\varphi}_i,\name{{\frak B}}_i,\kappa_i,\theta_i):i<\alpha\rangle$
is a nep--CS--iteration if: 
\begin{enumerate}
\item[$(\alpha)$] $\beta<\alpha\quad\Rightarrow\quad\bar{\bQ}\rest\beta$ is a
nep--CS--iteration, 
\item[$(\beta)$] if $\alpha=\beta+1$ then
\begin{enumerate}
\item[(i)]   $\bP_\beta=\Lim(\bar{\bQ}\rest\beta)$ (use clause (D))
\item[(ii)]  $\bar{\varphi}_\beta=\langle\varphi_{\beta,\ell}:\ell<3\rangle$
is formally as in the definition of nep (the substantial demand in (v) below,
but the parameter ${\frak B}_\beta$ is a name!) 
\item[(iii)] $\kappa_\beta,\theta_\beta$ are infinite cardinals (or ordinals)
\item[(iv)]  $\name{{\frak B}}_\beta$ is a $\bP_\beta$--name of a model with
universe $\kappa_\beta$ or even ${\cal H}_{< \aleph_1}(\kappa_\beta)$, whose
vocabulary is a fix countable one $\tau_0\subseteq {\cal H}(\aleph_0)$, but
for each atomic formula $\psi(x_0,\ldots,x_{n-1})$ and $\alpha_0,\ldots,
\alpha_{n-1}<\kappa_\beta$ the name of the truth value $\name{{\frak B}}_\beta
\models$``$\psi(\alpha_0,\ldots,\alpha_{n-1})$'' is an
$\hc$--$\kappa$-$\bP_\beta$--name (i.e.\ is defined by one $p=p^\beta_{
\psi^*}\in\cl_1(\bP_\beta)$)
\item[(v)] $\forces_{\bP_\beta}$``$\nbQ_\beta$ defined by
$\bar{\varphi}_\beta$ is temporarily straight explicitly $(\kappa_\beta,
\theta_\beta)$--nep as witnessed by $\name{{\frak B}}_\alpha$ and $\ZFC^-_*$
is good (hence is correct, see \ref{4.8A}(1))''.
\end{enumerate}
\end{enumerate}
\medskip

\noindent{\sc Clause (B)}, Definition:\qquad  We define $\kappa^\alpha=\sup[
\{\kappa_i:i< \alpha\}\cup\{ \alpha\}]$ (of course, if the result is an
ordinal we can replace it by its cardinality, coding it assuming the
$\kappa_i$'s are cardinals; remember that $\kappa_i \ge \theta_i$). 
\medskip

\noindent{\sc Clause (C)}, Definition:\qquad  We define ${\frak
B}^\alpha={\frak B}^{\bar{\bQ}}$, a model with universe $\subseteq {\cal H}_{<
\aleph_1}(\kappa^\alpha)$ or write $\kappa^\alpha$ and the usual vocabulary
such that 
\begin{enumerate}
\item[$(*)$]  ${\frak B}^\alpha$ codes (by its relations) $\alpha,\{(\beta,
\bar{\varphi}_\beta,\kappa_\beta,\theta_\beta):\beta<\alpha\}$ and $\langle
\name{{\frak B}}_\beta:\beta<\alpha\rangle$; i.e.\ for every atomic formula in
the vocabulary $\tau_0$ (so is of $\name{{\frak B}}_\beta$), $\psi=\psi(x_0,
\ldots,x_{n-1})$ for some function symbol $F_\psi$ we have:\quad if
$\alpha_\ell<\kappa_\beta$ for $\ell<n$ then $F_{\psi(\bar x)}(\beta;
\alpha_0,\ldots,\alpha_{n-1})$ is $p^\beta_{\psi(\alpha_0,\ldots,\alpha_{n-1}
)}$ (see clause (A)(iv)) and
\begin{quotation}
if the ${\frak B}$'s are on $\kappa$, we have also $F_{\psi,\ell}$, functions
of ${\frak B}^\alpha$ such that:\\
if $\alpha_\ell<\kappa_\beta$ for $\ell < n$ then $\{F_\psi(\beta,\ell,
\alpha_0,\ldots,\alpha_{n-1}):\ell<\omega\}$ lists the ordinals in ${\rm Tc}^{
\rm ord}(p^\beta_{\psi(\alpha_0,\ldots,\alpha_{n-1})})$ (the condition in
$\bP_\beta$ saying$\ldots$) and $F'_\psi$ codes how $p$ was gotten from them
(so we need $\kappa^\alpha\ge\omega_1$). 
\end{quotation}
\end{enumerate}
So in any case 
\begin{enumerate}
\item[$(**)$] if $N$ is a ${\frak B}^\alpha$--candidate, and $\beta\in\alpha
\cap N$ then $N$ is a ${\frak B}^\beta$--candidate. 
\end{enumerate}
\medskip

\noindent{\sc Clause (D)}, Definition:
\medskip

\noindent{\bf Case 1}:\quad If $\alpha = 0$ then $\bP_\alpha=\{\emptyset\}$.

\noindent{\bf Case 2}:\quad If $\alpha=\beta+1$ then
\[\begin{array}{ll}
\bP_\alpha=\bigl\{p:&p\mbox{ is a function, }\dom(p)\subseteq\alpha,\ p\rest
\beta\subseteq P_\beta\mbox{ and if }\beta\in\dom(p)\\
\ &\mbox{then for some } r=r_{p,\beta}\in \cl_1(\bP_\beta)\mbox{ determined by
$p$ we have:}\\
\ &p(\beta)\mbox{ is defined by  cases:}\\
\ &\mbox{\bf if } r[\name{G}_{P_\beta}]=\truth,\mbox{ it is in
}\cl(\nbQ_\beta),\mbox{ and an explicit witness }\\ 
\ &\mbox{is provided (say $p[\beta]$ codes it and having } r[\name{G}_{
\bP_\beta}]=\truth\mbox{ says so)},\\
\ &\mbox{{\bf if} not }, p(\beta)\mbox{ is }\emptyset=\emptyset_{\nbQ_\beta}=
\min(\nbQ_\beta)\bigr\}. 
  \end{array}\]
Pedantically, $p \in P_\alpha$ if and only if $p$ has the form $p'\cup\{
\langle\beta,\ell,x_\ell\rangle:\ell<3\}$ where $p'\in\bP_\beta$, $x_0\in\cl_1
(\bP_\beta)$, $x_1,x_2$ are $\hc$--$\kappa$--$\bP_\beta$--names of members of
${\cal H}_{<\aleph_1}(\theta_\beta)$ and $x_0$ is the truth value of ``$x_2[
G_\beta]$ is a witness to $x_1[G_\beta]\in\theta_\beta$''.

\noindent{\bf Case 3}:\quad If $\alpha$ is limit, then 
\[\bP_\alpha=\{p:p\mbox{ is a function, }\dom(p)\in [\alpha]^{\le\aleph_0}
\mbox{ and }\beta\le\alpha\quad\Rightarrow\quad p\rest\beta\in\bP_\beta\}.\]
\medskip

\noindent{\em The order}:\\
For $\alpha = 0$ nothing to do.\\
For $\alpha$ limit: $p\le q$ if and only if $\bigwedge\limits_{\beta<\alpha}
\bP_\beta\models$``$p\rest\beta\le q\rest\beta$'' (equivalently:
$\bigwedge\limits_{\beta\in\dom(\beta)}\bP_{\beta+1}\models$``$p\rest(\beta +
1)\le q\rest(\beta+1)$''), (see (C)).\\
For $\alpha=\beta + 1$:\quad the order is the transitive closure of the
following cases: 
\begin{enumerate}
\item[$(\alpha)$] $p\in\bP_\beta$, $q\in\bP_\alpha$, $\bP_\beta\models$``$p \le
q \rest \beta$'', 
\item[$(\beta)$]  $p(\beta)=q(\beta)$ and $\bP_\beta\models$``$p\rest\beta \le
q \rest \beta$'', 
\item[$(\gamma)$] $p\rest\beta=q\rest\beta$ and there is a ${\frak
B}^\alpha$--candidate $N$ such that $q\rest\beta$ is a $[0,\beta)$--canonical
$(N,\bP_\beta)$--generic above $p'\rest\beta$, $\bP_\beta\models$``$p'\rest
\beta\le q\rest\beta$'', $p'\in\bP^N_\alpha$ and
\[N\models\mbox{`` }p'\rest\beta\forces N[\name{G}_\beta]\models[\cl(
\nbQ_\beta)\models p(\beta)\le p'(\beta)\mbox{ and }p'(\beta)\in\nbQ_\beta]
\mbox{ ''}\]
and $q(\beta)$ is canonically generic for $(\bQ_\beta,N[\name{G}_\beta])$
above $p$, i.e.\  is 
\[p'(\beta) \ \&\  \bigwedge_{{\scriptstyle {\cal I}\in\pd(N,\bP_\alpha)}
\atop {\scriptstyle (\forall r \in {\cal I})r(\beta)\in\nbQ_\beta}}\bigvee
\{r(\alpha):N\models\mbox{`` }r \in {\cal I}\mbox{ '' and } r\rest\beta\in 
\name{G}_\beta\}\]
if $q\rest\beta\in\name{G}_{\bP_\beta}$ and $p'(\beta)$ if $q\rest\beta\notin 
\name{G}_{\bP_\beta}$. 
\end{enumerate}
Lastly, $\bP_{\alpha,w}=\{p\in\bP_\alpha:\dom(p)\subseteq w\}$ for $w\subseteq
\alpha$ (note that if $w \subseteq\beta\le\alpha$ we get the same forcing
notion). 
\medskip

\noindent{\sc Clause (E)}, Claim:\qquad  Trivial.
\medskip

\noindent{\sc Clause (F)}, Claim:\qquad Subclauses (a) and (c)--(f) are
trivial. 

\noindent{\bf Subclause (b)}:\quad  Here we should be careful as we do not ask
just that the order is forced but there is a $\hc$ witness; as we ask for a
witness and not explicit witness (see Definition \ref{4.8}) this is okay.  See
more in the proof of clause (I).
\medskip

\noindent{\sc Clause (G)}, Definition:

\noindent{\bf Case 1}:\quad  For $\beta<\alpha$ note that $N$ is also a
${\frak B}^\beta$--candidate and use the definition for $\bar{\bQ}\rest\beta$.

\noindent{\bf Case 2}:\quad  If $\gamma=\beta=\alpha$ -- trivial.

\noindent{\bf Case 3}:\quad  For $\beta=\alpha$, $\alpha = 0$ -- trivial.

\noindent{\bf Case 4}:\quad  For $\gamma<\beta=\alpha$ and $\beta=\beta'+1$,
$\beta'\notin w$ -- trivial. 

\noindent{\bf Case 5}:\quad  Suppose $\gamma<\beta=\alpha$, $\alpha=\beta'+
1$, $\beta' \in w$.\\
Then:\quad $q \rest\beta'$ is $[\gamma,\beta')$--canonically $(N,\bP_\beta,
w)$--generic and for some $\name{\tau}$, 
\[\begin{array}{ll}
N\models&\mbox{`` }\name{\tau}\mbox{ is a $\hc$--$\kappa_\beta$--$\bP_{\beta
\cap w}$--name of a member of }\nbQ_\beta\\ 
\ &\mbox{which is above $p(\beta)$ (which is in $\cl(\nbQ_\beta)$!) and is in
}N\langle \name{G}_{\bP_{\beta \cap w}}\rangle\mbox{ ''} 
  \end{array}\]
and 
\[q(\beta)=\name{\tau} \ \&\ \bigwedge_{{\scriptstyle {\cal I} \in\pd(N,
  \bP_{\alpha,w})}\atop {\scriptstyle (\forall r \in {\cal I})r(\beta)\in
\nbQ_\beta}}\bigvee\{r(\beta):N\models\models{`` }r\in {\cal I}\mbox{ '' and
}r \rest\beta\in \name{G}_{\bP_\beta}\}.\]

\noindent{\bf Case 6}:\quad $\gamma<\beta=\alpha$, $\beta$ a limit.\\
Say that diagonalization was used.
\medskip

\noindent{\sc Clause (H)}, Theorem:\qquad  Prove by induction.
\medskip

\noindent{\sc Clause (I)}, Theorem:\qquad  We have defined ${\frak B}^\alpha$
and $\kappa^\alpha$ (so $\theta^\alpha=\kappa^\alpha$). The formulas $\varphi^{
\bP_\alpha}_\ell$ ($\ell<3$) are implicitly defined (in the induction).

Why $\varphi^{\bP_\alpha}_0$ is absolute enough?  As the demand on $p(\beta)$
above says that $r_{p\rest (\beta+1),\beta}$, the witness for $p(\beta)\in
\cl(\bQ)$, is such that $r[\name{G}_{\bP_\beta}]=\truth$ gives all the
required information.

Why $\varphi^{\bP_\alpha}_1$ is absolute enough?  Because the canonical
genericity is about $\varphi_2$ and the properness requirement, see clause
(G), fit. 
\medskip

Now one proves by induction on $\beta \le \alpha$:
\begin{enumerate}
\item[$(\otimes)$]  {\em if}\ $N$ is a ${\frak B}^\alpha$--candidate, $w\in
N$, $N\models$``$w \subseteq \alpha$'', $\gamma_0\le\gamma_1\le\beta$,
$\{\gamma_0,\gamma_1,\beta\}\subseteq (\alpha+1)\cap N\cap w$, $p\in
\bP^N_\beta$, $q\in\bP_\gamma$, $p\rest\gamma_1\le q$, $q$ is $[\gamma_0,
\gamma_1)$--canonically $(N,\bP_{\gamma_1},w)$--generic  

\noindent {\em then}\ we can find $q^+$ such that:
\begin{enumerate}
\item[$(\alpha)$] $q^+\in\bP_\beta$, $q^+\rest\gamma=q$,
\item[$(\beta)$]  $p\leq^+$,
\item[$(\gamma)$] $q^+$ is $[\gamma,\beta)$--canonically
$(N,\bP_\beta,w)$--generic. 
\end{enumerate}
\end{enumerate}

\noindent{\sc Clause (I)}, Theorem:\qquad Straight. \QED$_{\ref{4.9}}$

\begin{proposition}
\label{4.10}
The iteration in \ref{4.9} is equivalent to the CS iteration.  More formally,
assume 
\[\bar{\bQ}=\langle(\bP_i,\nbQ_i,\bar{\varphi}_i,\name{{\frak B}}_i,\kappa_i,
\theta_i):i<\alpha\rangle\mbox{ is an CS--nep iteration}.\]
We can define $\bar{\bQ}'=\langle\bP'_i,\nbQ'_i:i<\alpha\rangle$ and $\langle
F_i:i<\alpha\rangle$ such that 
\begin{enumerate}
\item[(a)] $\bar{\bQ}'$ is a CS iteration,
\item[(b)] $F_i$ is a mapping from $\bP_i$ into $\bP'_i$,
\item[(c)] $j<i\quad\Rightarrow\quad F_j = F_i\rest\bP_j$,
\item[(d)] $F_i$ is an embedding of $\bP_i$ into $\bP'_i$ with dense range,
\item[(e)] $\nbQ_i$ is mapped by $F_i$ to $\nbQ'_i$.
\end{enumerate}
\end{proposition}

\Proof Straight. \QED

\begin{proposition}
\label{4.11}
In the context of \ref{4.10}:
\begin{enumerate} 
\item Assume that each $\name{{\frak B}}_\beta$ is essentially a real; i.e.\
$\kappa_\beta=\omega$ and if $R$ is in the vocabulary of $\name{{\frak
B}}_\beta$ then $R^{\name{{\frak B}}_\beta}\subseteq {}^{n(R)}\omega$. If
$\alpha<\omega_1$ then so is the ${\frak B}_\alpha$.  (If $\alpha\ge\omega_1$
we get weaker results).  
\item Assume that $\forces_{\bP_\beta}$`` the universe of ${\frak B}_\beta$ is
$\kappa_\beta$ ''. Then we can make ``${\frak B}_\beta$ has universe
$\kappa^\beta$'', coding the $p^\beta_\psi$'s. 
\end{enumerate}
\end{proposition}

\Proof  Left to the reader. \QED

\begin{remark}
\label{4.12}
1)\quad Note that \ref{4.9}, \ref{4.10} (and \ref{4.11}(4)) say something even
for $\alpha = 1$ so it speaks on $\cl(\bQ_0)=\bP_1$ (or $\cl_3(\bQ_0)=\bP_1)$.

\noindent 2)\quad Concerning \ref{4.11} note that if $\kappa({\frak B})\geq
\omega_1$, the difference between nep and snep is not large, so the case
$\alpha<\omega_1$ has special interest.

\noindent 3)\quad In \ref{4.9}, \ref{4.10}, we can replace the use of $\bQ'=
\cl(\bQ)$ from \ref{4.1} by $\cl_3(\bQ)$ from Definition \ref{4.3} (using
\ref{4.4}).

\noindent 4)\quad We can derive a theorem on local in \ref{4.10}, but for
strong enough $\ZFC^-_*$, then any follows.
\end{remark}

Of course, we can get forcing axioms.

\begin{proposition}
\label{4.14}
\begin{enumerate}
\item Assume for simplicity that $\V\models 2^{\aleph_0}=\aleph_1\ \&\
2^{\aleph_1}=\aleph_2$. Then for some proper $\aleph_2$--c.c. forcing notion
$\bP$ of cardinality $\aleph_2$ we have in $\V^{\bP}$: 
\begin{enumerate}
\item[$(\oplus)$] ${\rm Ax}_{\omega_1}[(\aleph_1,\aleph_1)\mbox{--nep}]$:
\qquad {\em if}\ $\bQ$ is a $(\kappa,\theta)$--nep forcing notion,
$\kappa,\theta\le\aleph_1$ and ${\cal I}_i$ is a dense subset of $\bQ$ for
$i<\omega_1$ and $\name{S}_i$ as a $\bQ$--name of stationary subset of
$\omega_1$ for $i<i(*) \le\omega_1$,\\
{\em then}\ for some directed $G\subseteq\bQ$ we have:\quad $i<\omega_1
\Rightarrow G\cap {\cal I}_i \ne\emptyset$ and 
\[\name{S}_i[G]\stackrel{\rm def}{=}\{\zeta<\omega_1:\mbox{for some } q\in G
 \mbox{ we have } q \forces_{\bQ}\mbox{`` }\zeta\in\name{S}_i\mbox{ ''}\}\]
is a stationary subset of $\omega_1$.
\end{enumerate}
\item We can demand that $\bP$ is explicitly $(\aleph_2,\aleph_2)$--nep
provided that in $(\oplus)$ we add ``explicitly simply'' to the requirements
on $\bQ$. 
\item In parts 1) and 2), we can strengthen $(\oplus)$ to ${\rm AX}_{\omega_1}
[\mbox{nep}]$. 
\end{enumerate}
\end{proposition}

\Proof Straight (as failure of ``$\bQ$, i.e.\ $\bar{\varphi}$, is nep'' is
preserved when extending the universe). \QED

\begin{proposition}
\label{4.15}
We can generalize the definitions and claims so far by:
\begin{enumerate}
\item[(a)] a forcing notion $\bQ$ is $(\bQ,\le,\le_{\pr},\emptyset_{\bQ})$,
where $\le_{\pr}$ is a quasi order, $p \le_{\pr} q \Rightarrow p \le q$ and
$\emptyset_{\bQ}$ the minimal element;
\item[(b)] in the definition of nep in addition to $\varphi_1$ we have
$\varphi_{1,\pr}$ defining $\le_{\pr}$, which is upward absolute from
$\bQ$--candidates, and in Definition \ref{0.2}(2)(c) we strengthen $p \le q$
to $p\le_{\pr} q$;
\item[(c)] the definition of CS iteration $\langle\bP_i,\nbQ_i:i<\alpha
\rangle$ is modified in one of the following ways:
\begin{enumerate}
\item[$(\alpha)$] $\bP_i=\bigl\{p:p$ is a function, $\dom(p)$ is a countable
subset of $i$, $j\in\dom(p) \Rightarrow \forces_{\bP_j}$``$p(j)\in\nbQ_j$'' 
and the set $\{j\in\dom(p):\ \not\forces_{\bP_j}$``$\emptyset_{\nbQ_j}\le_{\pr}
p(j)$''$\}$ is finite $\bigr\}$,

with the order

$p \le q$ if and only if $j\in\dom(p)\Rightarrow q\rest j\forces p(j)
\le^{\nbQ_j} q(j)$ and the set $\{j\in\dom(p):q\rest j\not\forces_{\bP_j}
\le^{\nbQ_j}_{\pr} q(j)$''$\}$ is finite;

\noindent (if each $\le^{\nbQ_j}_{\pr}$ is equality, this is FS iteration)
\item[$(\beta)$] $\bP_i=\bigl\{p:p$ is a function, $\dom(p)$ is a countable
subset of $i$, $j\in\dom(p)\ \Rightarrow\ \forces_{\bP_j}$``$p(j)\in\nbQ_j$''
$\bigr\}$,

with the order

$p \le q$ if and only if $j\in\dom(p)\Rightarrow q\rest j\forces_{\bP_j}$``$p
(j)\le^{\nbQ_j}q(j)$'' and $\{j\in\dom(p):q\rest j\not\forces_{\bP_j}$``$p(j) 
\le^{\nbQ_j}_{\pr} q(j)$''$\}$ is finite;
\end{enumerate}
\item[(d)] similarly for the CS-nep iteration.
\end{enumerate}
\end{proposition}

\Proof Left to the reader. \QED

\stepcounter{section}
\subsection*{\quad 6. When a real is $(\bQ,\eta)$--generic over $\V$}

\begin{definition}
\label{5.1}
\begin{enumerate}
\item We say that $(\bQ,\bar{W})$ is a temporary $({\frak B},\theta,\sigma,
\tau)$--pair if for some $\bQ$--name $\name{\eta}$ the following conditions
are satisfied:
\begin{enumerate}
\item[(a)] $\bQ$ is a nep-forcing notion for $({\frak B},\bar{\varphi},
\theta)$; possibly ${\frak B}$ expands ${\frak B}^{\bQ}$,
\item[(b)] $\forces_{\bQ}$``$\name{\eta}\in{}^\sigma\tau$'',
\item[(c)] $\bar{W}=\langle W_n:n<\sigma\rangle$,
\item[(d)] for each $n < \sigma$, $W_n\subseteq\{(p,\alpha):p \in\bQ,\alpha 
<\tau\}$,
\item[(e)] if $(p_\ell,\alpha_\ell)\in W_n$ for $\ell =1,2$ and $\alpha_1,
\alpha_2$ are not equal, then $p_1,p_2$ are incompatible in $\bQ$,
\item[(f)] for each $n<\sigma$ the set ${\cal I}_n = {\cal I}_n[\bar W]
\stackrel{\rm def}{=}\{p:(\exists \alpha)[(p,\alpha)\in W_n]\}$ is a predense
subset of $\bQ$, 
\item[(g)] so $\tau=\tau[\bar{W}]=\tau[\bQ,\bar{W}]$ and (abusing notation)
let $\sigma= \sigma[\bar{W}]=\sigma[\bQ,\bar{W}]$. 
\end{enumerate}
\item For $(\bQ,\bar{W})$ as above, $\name{\eta}=\name{\eta}[\bar{W}]=
\name{\eta}[\bQ,\bar{W}]$ is the $\bQ$--name 
\[\bigcup\{(p,(n,\alpha)):(\exists p\in \name{G}_{\bQ})((p,(n,\alpha))\in
W_n), \mbox{ so }n<\sigma\}.\] 
\item We replace the temporary by $K$ if this (specifically the demand (f))
holds in any $K$--extension. 
\item We may write $(\bQ,\name{\eta}),\bar{W}=\bar{W}^{\name{\eta}}$ abusing
notation. If we omit ${\frak B}$ we mean ${\frak B} = {\frak B}^{\bQ}$. If
$\tau=\aleph_0$ we may omit it; if $\tau=\sigma=\aleph_0$ we may omit them, if
$\theta=\sigma=\tau=\aleph_0$, we may write $\kappa$.
\item We say that $\name{\eta}[\bQ,\bar{W}]$ is a temporarily generic real (or
function) for $\bQ$ if for no distinct $G_1,G_2 \subseteq \bQ$ generic over
$\V$ do we have $\name{\eta}[G_1]=\name{\eta}[G_2]$. 
\item Instead $(\bQ,\bar{W})$ we may write $(({\frak B}^{\bQ},
\bar{\varphi}^{\bQ},\theta^{\bQ}),\bar{W})$ (or with $\name{\eta}$ instead
$\bar{W}$). 
\end{enumerate}
\end{definition}

\begin{definition}
\label{5.2A}
\begin{enumerate}
\item Let $\cK_{\kappa,\theta,\sigma,\tau}$ be the class of all $(\bQ,
\name{\eta})$ which are temporary $({\frak B},\theta,\sigma,\tau)$--pairs for
some ${\frak B}$ with $\kappa({\frak B})\le\kappa,\|{\frak B}\|\le\kappa$. 
\item Let $(\bQ,\name{\eta})$ be a temporary $(\kappa,\theta)$--pair (actually
more accurately write $({\frak B},\bar{\varphi},\theta)$, $\bar{W}$); so
$\sigma=\tau=\aleph_0$. 

Let $N$ be a $\bQ$--candidate and $\eta\in\baire$. We say that $\eta$ is a
$(\bQ,\name{\eta})$--generic real over $N$\ {\em if}\ for some $G\subseteq
\bQ^N$ which is generic over $N$ we have $\eta=\name{\eta}[G]$. 
\item We say that $\name{\eta}=\name{\eta}[\bar{W}]$ is hereditarily countable
if each $W_n$ is countable (note: the generic reals of the forcing notions
from \cite{RoSh:470} are like that, but for our purpose just ``absolute
enough'' suffices). 
\end{enumerate}
\end{definition}

\begin{definition}
\label{5.2B}
\begin{enumerate}
\item $(\bQ,\bar{W})$ is a temporary explicitly $({\frak B},\theta,\sigma,
\tau)$--pair (or nep pair) if for some $\bQ$--name $\name{\eta}$ we have: 
\begin{enumerate}
\item[(a)] $\bQ$ is an explicit nep forcing notion for $({\frak B},
\bar{\varphi},\theta)$,
\item[(b)] $\forces_{\bQ}$``$\name{\eta}\in {}^\sigma \tau$'',
\item[(c)] $\bar{W}=\langle\psi_{\alpha,\zeta}:\alpha<\sigma,\zeta<\tau
\rangle$, 
\item[(d)] $\psi_{\alpha,\zeta}\in\cl_1(\bQ)$ for $\alpha<\sigma$, $\zeta<
\tau$,
\item[(e)] $\forces_{\bQ}$`` $\name{\eta}(\alpha)=\zeta$ iff $\psi_{\alpha,
\zeta}[\name{G}_{\bQ}]=\truth$ ''.
\end{enumerate}
\item In this case $\name{\eta}=\name{\eta}[\bar{W}]=\name{\eta}[\bQ,
\bar{W}]$ is the $\bQ$--name above (it is unique). Abusing notation we may
write $(\bQ,\name{\eta})$ instead $(\bQ,\bar{W})$ and then let $\bar{W}=
\bar{W}[\name{\eta}]=\bar{W}[\bQ,\name{\eta}]$.
\item We introduce the notions from \ref{5.1}(3)--(6) for the current case
with almost no changes.
\end{enumerate}
\end{definition}

\begin{definition}
\label{5.2C}
$\cK^{\ex}_{\kappa,\theta,\sigma,\tau}=\{(\bQ,\name{\eta})\in K_{\kappa,
\theta,\sigma,\tau}:(\bQ,\name{\eta})$ is temporarily explicitly $({\frak
B},\theta,\sigma,\tau)$--pair for some model ${\frak B}$ with $\kappa({\frak
B}) \le \kappa,\|{\frak B}\| \le \kappa\}$.
\end{definition}

\begin{proposition}
\label{5.3A}
Assume that:
\begin{enumerate}
\item[(a)] $\bQ$ is an explicitly nep forcing notion which satisfies the c.c.c.
\item[(b)] $\forces_{\bQ}$``$\name{\eta}\in {}^\sigma \omega$'' and (for
$\alpha<\sigma$ and $m<\omega$) $\psi_{\alpha,m}\in\cl_1(\bQ)$ are such that
\[\forces_{\bQ}\mbox{`` }\name{\eta}(\alpha)=m\ \mbox{ iff }\ \psi_{\alpha,m}
[\name{G}_{\bQ}]=\truth\mbox{ ''}\]
\item[(c)] $\bQ'\stackrel{\rm def}{=}B_2(\bQ,\name{\eta})$ is the following
suborder of $\cl_2(\bQ)$:
\[\begin{array}{ll}
\{p\in\cl_2(\bQ):&p\mbox{ is generated by the }\psi_{\alpha,m}\mbox{ i.e.\ it
belongs to the closure of }\\
\ &\{\psi_{\alpha,m}:\alpha<\sigma,m<\omega\}\mbox{ under $\neg,
\bigwedge\limits_{i<\gamma}$ for $\gamma<\omega_1$ in }\cl_2(\bQ)\}
  \end{array}\]
(i.e.\ it is the quasi order $\le^{\bQ}_2$ restricted to this set).
\end{enumerate}
{\em Then}:
\begin{enumerate}
\item $\bQ'\lesdot\cl_2(\bQ)$ and $\name{\eta}\in {}^\sigma \omega$ is a
generic function for $\bQ'$. 
\item Assume additionally that
\begin{enumerate}
\item[$(*)$] if $M$ is a $\bQ$--candidate, $M\models$``${\cal I}$ is a maximal
antichain of $\bQ$'',\\
then ${\cal I}^M$ is a maximal antichain of $\bQ$. 
\end{enumerate}
Then we also have
\begin{enumerate}
\item[$(\alpha)$]  $\bQ'$ is $(\kappa,\theta)$-nep c.c.c. forcing notion,
\item[$(\beta)$]    if $\bQ$ is simple, then $\bQ'$ is simple,
\item[$(\gamma)$]   if $\bQ$ is $K$--local, then $\bQ'$ is $K$--local.
\end{enumerate}
\end{enumerate}
\end{proposition}

\Proof Straight. \QED

Now the hypothesis $(*)$ in \ref{5.3A}(2) is undesirable, so we use $B_3(\bQ,
\name{\eta})$ (see \ref{5.3B}(c) below), which has a suitable quasi order.

\begin{proposition}
\label{5.3B}
Assume that:
\begin{enumerate}
\item[(a)] $\bQ$ is explicitly nep forcing notion which satisfies the  c.c.c.
\item[(b)] $\forces_{\bQ}$``$\name{\eta}\in {}^\sigma \omega$'' and
$\psi_{\alpha,m}\in\cl_2(\bQ)$ are such that 
\[\forces_{\bQ}\mbox{`` }\name{\eta}(\alpha)=m\ \mbox{ iff }\ \psi_{\alpha,m}
[\name{G}_{\bQ}]=\truth\mbox{ ''},\]
\item[(c)] $\bQ'\stackrel{\rm def}=B_3(\bQ,\name{\eta})$ is a forcing notion
defined as follows:\\
the set of elements is like $B_2(\bQ,\eta)$; i.e.\ it is the closure of
$\{\psi_{\alpha,m}:\alpha<\sigma,m<\omega\}$ under $\neg,\bigwedge\limits_{i <
\gamma}$ for $\gamma<\omega_1$ inside $\cl_2(\bQ)$;\\
the quasi order $\le_3=\le^{B_3(\bQ,\name{\eta})}_3$ is $\leq^{\cl_3(\bQ)}$
restricted to $B_3(\bQ,\name{\eta})$,
\item[(d)] the statements $(\circledast_3)$ and $(\circledast_4)$ of \ref{4.4}
hold. 
\end{enumerate}
Then:
\begin{enumerate}
\item[$(\alpha)$]  $\bQ'$ is essentially a suborder of $\cl_2(\bQ)$; i.e.\
$\psi \in\bQ'\ \Rightarrow\ \psi\in\cl_2(\bQ)$, and for $\psi_1,\psi_2\in
\bQ'$ we have:\quad $\psi_1\leq_3 \psi_2\ \Leftrightarrow\ \psi_1\leq^{\cl_3(
\bQ)}\psi_2$,
\item[$(\beta)$] $\name{\eta}$ is a $\bQ'$--name,
$\forces_{\bQ'}$``$\name{\eta}\in {}^\sigma\omega$'' and $\name{\eta}$ is a
generic function for $\bQ'$, 
\item[$(\gamma)$] $\bQ'$ is explicitly nep c.c.c.\ forcing notion with ${\frak
B}^{\bQ'}={\frak B}^{\bQ}$, $\bar{\varphi}^{\bQ'}=\bar{\varphi}^{B_3(\bQ,
\name{\eta})}$, $\theta^{\bQ'}=\theta^{\bQ}$,
\item[$(\gamma)^+$] each forcing extension of $\V$ which preserves the
assumption (a) (hence also (b)) preserves $(\gamma)$,
\item[$(\delta)$]  if $\bQ$ is simple (or straight) then $\bQ'$ is simple
(or straight),
\item[$(\varepsilon)$] if $\bQ$ is $K$--local, then $\bQ'$ is $K$-local.
\end{enumerate}
\end{proposition}

\Proof Straight. \QED

\begin{proposition}
\label{5.3C}
In \ref{5.1}--\ref{5.3B} above, we can replace $\bQ$ by $\bQ\rest\{p\in\bQ:p
\ge q\}$ preserving the properties of $(\bQ,\name{\eta})$. \QED
\end{proposition}

\begin{fact}
\label{5.5A}
If $\bQ$ is simply correctly nep for $K$, $\bQ$ is in $\V$, and $\V_1$ is a
$K$--extension of $\V$ {\em then}
\begin{enumerate}
\item[(i)] in $\V_1$, $\bQ^{\V}\le_{ic}\bQ^{\V_1}$ (see \cite[Ch.IV]{Sh:f}),
i.e.\ for $p,q \in V_0$, ``$p\in\bQ$'', ``$p \le q$'', ``$\neg(p\le q)$'',
``$p,q$ compatible'',``$p,q$ compatible'' are preserved from $\V$ to $\V_1$, 
\item[(ii)] for $p,p_n\in\V$ the statements ``$p\notin\bQ$'' and ``${\cal I}=
\{p_n:n<\omega\}$ is predense above $p$ in $\bQ$'' are preserved from $\V$ to
$\V_1$, 
\item[(iii)] if $\bQ$ satisfies the c.c.c.\ then in clause (ii) above we can
omit the countability of ${\cal I}$.
\end{enumerate}
\end{fact}

\Proof Straight, for example:

``$p,q$ are incompatible'' iff there is no $\bQ$--candidate $M$ such that 
\[M\models\mbox{`` }p,q\mbox{ have a common $\leq_{\bQ}$--upper bound ''}.\]
So by Shoenfield--Levy absoluteness, if this holds in $\V$, it holds in
$\V_1$. 

\noindent (ii) Similarly.

\noindent (iii) Follows (and repeated in \ref{6.8}).  \QED

\begin{proposition}
\label{5.12} 
Let $(\bQ,\name{\eta})$ be temporarily explicitly nep pair. Assume $N$ is a
$\bQ$--candidate. If $N\models$``$\eta^*$ is $(\bQ,\name{\eta})$--generic over
a $\bQ$--candidate $M$'', {\em then} $\eta^*$ is a
$(\bQ,\name{\eta})$--generic for $M$. 
\end{proposition}

\Proof Straight. \QED

\begin{proposition}
\label{5.13}
Assume that:
\begin{enumerate}
\item[(a)] $\bQ$ is explicitly nep,
\item[(b)] $\bQ$ is c.c.c.\ moreover it satisfies the c.c.c.\ in every
$\bQ$--candidate, 
\item[(c)] incompatibity in $\bQ$ is upward absolute from $\bQ$--candidates
(but see \ref{5.5A}),
\item[(d)] $\name{\eta}$ is a $\hc$--$\kappa({\frak B}^{\bQ})$--$\bQ$--name of
a member of $\baire$ defined from ${\frak B}^{\bQ}$ (so we demand this in
every $\bQ$--candidate). 
\end{enumerate}
Furthermore, suppose that 
\begin{enumerate}
\item[(A)] $N_1,N_2$ are $\bQ$--candidates, $N_2$ is a generic extension
of $N_1$ for a forcing notion $\bR$, (so ${\frak B}^{N_2}={\frak B}^{N_1}$ and
$N_1\models$``$\bR$ is a forcing notion''),
\item[(B)] $N_1\models$`` for every countable $X \subseteq\bQ$ and $n<\omega$
there is a $\bQ$-candidate $N_0 \prec_{\Sigma_n} N_1$ to which $X$ and $\bR$
belong'',
\item[(C)] $\eta^*\in\baire$ is a $(\bQ,\name{\eta})$--generic real over
$N_2$.
\end{enumerate}
{\em Then}\ $\eta^*\in\baire$ is a $(\bQ,\name{\eta})$--generic real over
$N_1$. 
\end{proposition}

\begin{remark}
\label{5.13A}
\begin{enumerate}
\item In (B), we can replace $X$ by ``a maximal antichain of $\bQ,
\name{\eta}$''. 
\item Clearly we can replace ``maximal antichain'' by ``predense set'' or
``predense set over $p$'' (note ${\cal I}^{N_2}={\cal I}^{N_1}$ as $N_2 =
N^{\bR}_1$).
\item We can weaken ``$N_0 \prec_{\Sigma_n} N_1$'' in clause (B).
\end{enumerate}
\end{remark}

{\noindent{\sc Proof of \ref{5.13}} \hspace{0.2in}} Clearly it suffices to
prove that (assuming (a)-(d),(A),(B) and (C)):
\begin{enumerate}
\item[$(*)$] {\em if}\quad $N_1\models$``${\cal I}$ is a maximal antichain of
$\bQ$'',\\
{\em then}\quad ${\cal I}^{N_1}={\cal I}^{N_2}$ and $N_2\models$``${\cal I}$
is a maximal antichain of $\bQ$''. 
\end{enumerate}
Assume that this fails for ${\cal I}$. Then some $r\in\bR$ forces this failure
(in $N_1$). By assumption (b), in $N_1$ the set ${\cal I}^{N_1}$ is countable
so let $N_1 \models$``${\cal I}=\{p_n:n<\alpha\}$'',  where $\alpha\le\omega$.
Let $n<\omega$ be large enough. By clause (B) in $N_1$ there is a
$\bQ$--candidate $N_0$ to which ${\cal I}$ and $r$ and $\bR$ belong and
$N_0\prec_{\Sigma_n}N_1$. Since 
\[\begin{array}{r}
N_1\models{``}(\exists r\in\bR)[r\forces_{\bR}\mbox{``}{\cal I}\mbox{
is not a maximal antichain of $\bQ$}\quad\\
\mbox{(and $N_1[\name{G}_{\bR}]$ is a $\bQ$--candidate)'']''},
  \end{array}\]
there is $r_0\in\bR\cap N_0$ such that 
\[\begin{array}{r}
N_0\models\mbox{``}[r_0\forces_{\bR}\mbox{``${\cal I}$ is not a maximal
antichain of $\bQ$}\quad\\
\mbox{(and $N_0[\name{G}_{\bR}]$ is a $\bQ$--candidate)]''}.
  \end{array}\]
Now, as $N_1$ satisfies enough set theory and $N_1$ ``thinks'' that $N_0$ is
countable and $\bR^{N_0}$ is a forcing notion in $N_0$, there is in $N_1$ a
subset $G'_{\bR}$ of $\bR\cap N_0=\bR^{N_0}$ generic over $N_0$ to which $r_0$
belongs. So in $N_0[G'_{\bR}]$ there is $p\in\bQ^{N_0[G'_{\bR}]}$ incompatible
(in $\bQ^{N_0[G'_{\bR}]}$) with each $p_n$. By the assumption (c) this holds
in $N_1$, contradiction to the choice of ${\cal I}$ (see $(*)$). 
\QED$_{\ref{5.13}}$ 

\begin{definition}
\label{5.14}
\begin{enumerate}
\item We say that $\bar{\varphi}$ or $(\bar{\varphi},{\frak B})$ is a
temporary $(\kappa,\theta)$-definition of a strong c.c.c.--nep forcing notion
$\bQ$ if: 
\begin{enumerate}
\item[(a)] $\varphi_0$ defines the set of elements of $\bQ$ and $\varphi_0$ is
upward absolute from $({\frak B},\bar{\varphi},\theta)$--candidates, 
\item[(b)] $\varphi_1$ defines the partial ordering of $\bQ$ (even in $({\frak
B},\bar{\varphi},\theta)$--candidates) and $\varphi_1$ is upward absolute from
$({\frak B},\bar{\varphi},\theta)$--candidates,  
\item[(c)] for any $({\frak B},\bar{\varphi},\theta)$--candidate $N$, if $N
\models$``${\cal I}\subseteq\bQ$ is predense'', {\em then} also in $\V$,
${\cal I}^N$ is a predense subset of $\bQ$. 
\end{enumerate}
\item We say that $\bar{\varphi}$ or $(\bar{\varphi},{\frak B})$ is a
temporarily [explicitly] $(\kappa,\theta)$--definition of a c.c.c.--nep
forcing notion $\bQ$ if
\begin{enumerate}
\item[$(\alpha)$] it is a temporary [explicitly] $(\kappa,\theta)$-definition
of a nep forcing notion,
\item[$(\beta)$] for every $\bQ$--candidate $N$ we have $N\models$``$\bQ$
satisfies the c.c.c.''. 
\end{enumerate}
\item The variants are defined as usual.
\end{enumerate}
\end{definition}

\begin{proposition}
\label{5.15}
\begin{enumerate}
\item If $\bQ$ is strongly c.c.c.--nep forcing notion and $N_1\subseteq N_2$
are $\bQ$--candidates, {\em then} every $\eta$ which is
$(\bQ,\name{\eta})$--generic over $N_2$ is also $(\bQ,\name{\eta})$--generic
over $N_1$.  
\item If $\ZFC^-_*$ is normal and $\bQ$ is temporarily c.c.c.--nep {\em then}
$\bQ$ satisfies the c.c.c.
\end{enumerate}
\end{proposition}

\begin{comment}
\label{6.14}
We can spell out various absoluteness, e.g.
\begin{enumerate}
\item  If $\bQ$ is simple nep, c.c.c.\ and ``$\langle p_n:n<\omega\rangle$ is
predense'' has the form $(\exists t\in {\cal H}_{<\aleph_1}((\kappa+\theta)))
[t\models\ldots]$ (e.g.\ $\kappa^{\bQ}=\omega$ and it is $\Pi^1_2$) then
predensity of countable sets is preserved in any forcing extension. 
\item Note that strong c.c.c.--nep (from \ref{5.15}(1)) does not imply
c.c.c.--nep (from \ref{5.15}(2)). But if $\ZFC^-_{**}\vdash\ZFC^-_*$ and
$\ZFC^-_{**}$ says that $\ZFC^-_*$ is normal and $\bQ$ is strong c.c.c.--nep
for $\ZFC^-_*$, then $\bQ$ is c.c.c.--nep for $\ZFC^-_*$.
\end{enumerate}
\end{comment}

\stepcounter{section}
\subsection*{\quad 7. Preserving a little implies preserving much} 
Our main intention is to show that, for example if a ``nice'' forcing notion
$\bP$ satisfies $\forces_{\bP}$``$(\can)^{\V}$ is not null'', {\em then} it
preserves ``$X\subseteq\can\ \ (X\in\V)$ is not null''.

By Goldstern Shelah (\cite[Ch.XVIII, 3.11]{Sh:f}) if a Souslin proper forcing
preserves ``$(\baire)^{\V}$ is non-meagre'' then it preserves ``$X\subseteq
\baire$ is non-meagre'' and more (in a way suitable for the preservation
theorems there). 

The main question not resolved there was: is it special for Cohen forcing
(which is a way to speak on non-meagre), or it holds for nice c.c.c.\ forcing
notions in general, in particular does a similar theorem hold for ``non-null''
instead of ``non-meagre''. Though there have been doubts about it, we succeed
to do it here. In fact, even for a wider family of forcing notions but we have
to work more in the proof. 

See \S11 on a generalization. The reader may concentrate on the case that
$\bQ$ is strongly c.c.c.\ nep and $\bP,\bQ$ are explicitly $\aleph_0$--nep and
simple. It is natural to assume that $\name{\eta}$ is a generic real for $\bQ$
but we do not ask for it when not used.
\bigskip

\begin{convention}
\label{6.0}
\begin{enumerate}
\item $\bQ$ is an explicitly nep forcing notion. 
\item $\name{\eta}\in \baire$ is a hereditarily countable $\bQ$-name which is
${\frak B}$--definable. 
\end{enumerate}
\end{convention}
We would like to preserve something like: ``$x$ is $\bQ$--generic over $N$''.

\begin{definition}
\label{6.1}
\begin{enumerate}
\item  $I_{(\bQ,\name{\eta})}\stackrel{\rm def}{=}\{A\in\Borel(\baire):\
\forces_{\bQ}$``$\name{\eta}\notin A$''$\}$ (it is an ideal on the Boolean
algebra of Borel subsets of $\baire$). 
\item $I^{\ex}_{(\bQ,\name{\eta})}$ is the ideal generated by $I_{(\bQ,
\name{\eta})}$ on ${\cal P}(\baire)$. (So for $A\in\Borel(\baire)$ we have: $A
\in I_{(\bQ,\name{\eta})}\ \Leftrightarrow\ A\in
I^{\ex}_{(\bQ,\name{\eta})}$). Let 
\[\begin{array}{ll}
\hspace{-0.5cm}I^\dx_{(\bQ,\name{\eta})}\stackrel{\rm def}{=}\{X\subseteq
\baire:&\mbox{for a dense set of $q\in \bQ$, for some Borel set $B\subseteq
\baire$,}\\
&\mbox{we have $X\subseteq B$ and $q\forces$``$\name{\eta}\notin B$''}\}.
  \end{array}\]
\item For an ideal $I$ (on Borel sets, respectively), the family of
$I$--positive (Borel, respectively) sets is denoted by $I^+$.\\
(Thus, for a Borel subset $A$ of $\baire$, $A\in I^+_{(\bQ,\name{\eta})}$\ 
{\em iff}\ there is $q\in\bQ$ such that $q\forces_{\bQ}$``$\name{\eta}\in
A$''.
\end{enumerate}
\end{definition}

\begin{definition}
\label{6.3}
\begin{enumerate}
\item A forcing notion $\bP$ is $I_{(\bQ,\name{\eta})}$--preserving if for
every Borel set $A$ 
\[A\in (I_{(\bQ,\name{\eta})})^+\quad \Rightarrow\quad \forces_{\bP}\mbox{``
}A^{\V}\in (I^{\ex}_{(\bQ,\name{\eta})})^+\mbox{''}\]
($A^{\V}$ means: the same set, i.e.\ $A \cap\V$).
\item $\bP$ is strongly $I_{(\bQ,\name{\eta})}$--preserving if for all $X
\subseteq\baire$ (i.e.\ not only Borel sets)
\[X\in (I^{\dx}_{(\bQ,\name{\eta})})^+\quad \Rightarrow\quad \forces_{\bP}
\mbox{``}X\in (I^{\ex}_{(\bQ,\name{\eta})})^+\mbox{''}.\]
[See \ref{6.4}(7) for $\bQ$ which is c.c.c.]
\item  We say that a forcing notion $\bP$ is weakly $I_{(\bQ,
\name{\eta})}$--preserving if $\forces_{\bP}$`` $(\baire)^{\V}\in
(I^{\ex}_{(\bQ,\name{\eta})})^+$ ''.
\item $\bP$ is super--$I_{(\bQ,\name{\eta})}$--preserving if for all
$X\subseteq\baire$ we have:
\[X\in(I^\dx_{(\bQ,\name{\eta})})^+\quad\Rightarrow\quad \forces_{\bP} X^\V\in 
(I^\dx_{(\bQ,\name{\eta})})^+.\]
\end{enumerate}
\end{definition}

\begin{proposition}
\label{6.4}
\begin{enumerate}
\item $I_{(\bQ,\name{\eta})}$ is an $\aleph_1$--complete ideal (in fact, if
$\langle A_i:i\le\alpha\rangle\in\V$, each $A_i\in\Borel(\baire)$ and
$\forces_{\bQ}$``$A^{\V[\name{G}]}_\alpha\subseteq\bigcup\limits_{i<\alpha} 
A^{\V[\name{G}]}_i$'' and $A_i\in I_{(\bQ,\name{\eta})}$ for $i<\alpha$ then
$A_\alpha\in I_{(\bQ,\name{\eta})}$). 
\item If $(\bQ,\name{\eta})$ is not trivial (i.e.\
$\forces_{\bQ}$``$\name{\eta}\notin (\baire)^{\V}$), {\em then} singletons
belong to $I_{(\bQ,\name{\eta})}$. 
\item $\baire\notin I_{(\bQ,\name{\eta})}$.
\item Assume ($\ZFC^-_*$ is $K$--good and) $\bQ$ is correct. If in $\V$, $X
\in I^{\ex}_{(\bQ,\name{\eta})}$ and $\bP\in K$, {\em then}\ in $\V^{\bP}$
still $X\in I^{\ex}_{(\bQ,\name{\eta})}$ (but see later).
\item Assume ($\ZFC^-_*$ is $K$--good, particularly (c) of \ref{0.9} and)
$\bQ$ is correct. If, in $\V$, $B$ is a Borel subset of $\baire$ from $I_{(
\bQ,\name{\eta})}$ and $\V_1=\V^{\bP}$ then also $\V_1\models$``$B\in I_{\bQ,
\name{\eta}}$''. 
\item $I^\ex_{(\bQ,\name{\eta})}$, $I^\dx_{(\bQ,\name{\eta})}$ are ideals of
${\cal P}(\baire)$ and
\[I^\ex_{(\bQ,\name{\eta})}\rest(\mbox{the family of Borel sets})= I^\dx_{(
\bQ,\name{\eta})}\rest(\mbox{the family of Borel sets}).\]
\item If $\bQ$ satisfies the c.c.c.\ then $I^\dx_{(\bQ,\name{\eta})}$ is
generated by $I_{(\bQ,\name{\eta})}$, so equal to $I^\ex_{(\bQ,\name{\eta})}$.
\item $I^\ex_{(\bQ,\name{\eta})}$ is $\aleph_1$--complete.
\item If for some stationary $S\subseteq [\chi]^{\aleph_0}$, $\bQ$ is
$S$--proper then $I^\dx_{(\bQ,\name{\eta})}$ is $\aleph_1$--complete.
\end{enumerate}
\end{proposition}

\Proof We will prove parts 5) and 4) only, the rest is left to the reader.

\noindent 5)\quad First work in $\V^\bP$. If the conclusion fails then for
some $q\in\bQ$ we have $q\forces$``$\name{\eta}\in B$''. So there is a
$\bQ$--candidate $M$ to which $q,B$ (i.e.\ the code of $B$) belong. There is
$q'$ such that $q\leq q'$ and $q'$ is $(M,\bQ)$--generic. Now for every
$G\subseteq\bQ$ generic oner $\V^\bP$, $\name{\eta}[G]\in B^{\V^\bP[G]}$. By 
absoluteness, also $M\langle G\rangle\models\name{\eta}\langle G\cap\bQ^M
\rangle\in B^{M\langle G\rangle}$ and hence (by the forcing theorem) for some
$p\in G\cap\bQ^M$ we have $M\models [p\forces_\bQ\mbox{``}\name{\eta}\in B
\mbox{''}]$. Now, returning to $\V$, by Shoenfield--Levy absoluteness there
are such $M',p'$ in $\V$. Let $p''$ be $(M',\bQ)$--generic, $p'\leq^\bQ p''$. 
So similarly to the above, $p''\forces_{\bQ}$``$\name{\eta}\in B$''.

\noindent 4)\quad As $X\in I^\ex_{(\bQ,\name{\eta})}$, clearly for some Borel
set $B\in I_{(\bQ,\name{\eta})}$ we have $X\subseteq B$. By part (5), also in
$\V^\bP$ we have $B\in I^\ex_{(\bQ,\name{\eta})}$ and trivially $X\subseteq
B^\V\subseteq B^{\V^\bP}$. \QED

\begin{proposition}
\label{6.4A}
\begin{enumerate}
\item If a forcing notion $\bP$ is $I_{(\bQ,\name{\eta})}$--preserving, {\em
then} $\bP$ is weakly $I_{(\bQ,\name{\eta})}$--preserving. 
\item If $\bP$ is strongly $I_{(\bQ,\name{\eta})}$--preserving, {\em then}
$\bP$ is $I_{(\bQ,\name{\eta})}$--preserving. 
\item Assume that $\bQ$ satisfies the c.c.c.\ and $(\bQ,\name{\eta})$ is
homogeneous (see $(\circledast)$ below). Then:\quad $\bP$ is $I_{(\bQ,\name{
\eta})}$--preserving iff $\bP$ is weakly $I_{(\bQ,\name{\eta})}$--preserving,
where
\begin{enumerate}
\item[$(\circledast)$] $(\bQ,\name{\eta})$ is homogeneous if:

for any (Borel) sets $B_1,B_2\in (I_{(\bQ,\name{\eta})})^+$ we can find a
Borel set $B_1'\subseteq B_1$, $B_1\in (I_{(\bQ,\name{\eta})})^+$ and a Borel
function $F$ from $B_1'$ into $B_1$ such that
\begin{enumerate}
\item[$(\alpha)$] for every Borel set $A\in I_{(\bQ,\name{\eta})}$, $F^{-1}[
A\cap B_2]\in I_{(\bQ,\name{\eta})}$,
\item[$(\beta)$]  this is absolute (or at least it holds also in $\V^\bP$).
\end{enumerate}
\end{enumerate}
\end{enumerate}
\end{proposition}

\Proof 3)\quad By part (1) it suffices to show ``non--preserving'' assuming
``not weakly preserving''. So there are $p,B^*,\name{A},\name{q}$ such that
$B\in (I_{(\bQ,\name{\eta})})^+$ is a Borel subset of $\baire$ and
\[\begin{array}{ll}
p\forces_{\bP}\mbox{``}&(a)\ \name{A}\mbox{ is a Borel set}\\
\ &(b)\ \name{q}\in\bQ\mbox{ (in $\V^\bP$!)}\\
\ &(c)\ \name{q}\mbox{ witnesses }\name{A}\in I_{(\bQ,\name{\eta})},\mbox{
that is }\name{q}\forces_{\bQ}\mbox{``}\name{\eta}\notin \name{A}\mbox{''}\\
\ &(d)\ \nu\in\name{A}\mbox{ for every }\nu\in (B^*)^\V.\ \mbox{ ''}
  \end{array}\]
Let 
\[\begin{array}{ll}
{\cal J}=\{B:&B\in (I_{(\bQ,\name{\eta})})^+,\mbox{ so a Borel subset of
}\baire, \mbox{ and}\\
\ &\mbox{for some Borel one-to-one function $F$ from $B$ to $B^*$ we have}\\
\ &\mbox{$F$ is absolutely $(I_{(\bQ,\name{\eta})})^+$--preserving }\}.
  \end{array}\]
Choose a maximal family $\{B_i:i<i^*\}\subseteq{\cal J}$ such that $i\neq j\
\Rightarrow\ B_i\cap B_j\in I_{(\bQ,\name{\eta})}$. As $\bQ$ satisfies the
c.c.c.\ necessarily $i^*<\omega_1$, so wlog $i^*\leq\omega$. By the
assumption, $\baire\setminus\bigcup\limits_{i<i^*} B_i\in
I_{(\bQ,\name{\eta})}$. Let $F_i$ witness that $B_i\in {\cal J}$. Let 
\[\name{A}_i=\{\eta\in\baire:\eta\in B_i\ \mbox{ and }\
F_i(\eta)\in\name{A}\}.\] 
Then $\name{A}_i$ is a Borel subset of $\baire$ and
$p\forces_{\bP}$``$\name{A}_i\in I_{(\bQ,\name{\eta})}$'' as
$p\forces_{\bP}$``$\name{A}\in I_{(\bQ,\name{\eta})}$''. Hence 
\[p\forces\mbox{`` }\bigcup_{i<i^*}\name{A}_i\cup (\baire\setminus\bigcup_{i<
i^*})\in I_{(\bQ,\name{\eta})}\mbox{ ''}\]
(call this set $\name{A}^*$). Now
\[p\forces_{\bP}\mbox{`` }(\baire)^{\V}=(\baire\setminus\bigcup_{i<i^*}B_i)^\V
\cup\bigcup_{i<i^*} B_i^\V\subseteq (\baire\setminus \bigcup_{i<i^*} B_i)\cup
\bigcup_{i<i^*} A_i\in I_{(\bQ,\name{\eta})}\mbox{ ''}\]
so we are done. \QED
\medskip

\noindent{\bf Comment}:\qquad  1)\quad It is easy to find a forcing notion
$\bP$ which is $I_{(\bQ,\name{\eta})}$--preserving, but not strongly $I_{(\bQ,
\name{\eta})}$--preserving, e.g.\ for $\bQ =$ Cohen (see \ref{6.6} below). 
However, for sufficiently nice forcing notion $\bP$, ``$I_{(\bQ,\name{\eta})
}$--preserving'' and ``strongly $I_{(\bQ,\name{\eta})}$-preserving'' coincide,
as we will see in \ref{6.5}. (Parallel to the phenomenon that for ``nice''
sets, CH holds). 

\noindent 2)\quad It is even easier to find a weakly $I_{(\bQ,\name{\eta})
}$--preserving forcing notion $\bP$ which is not $I_{(\bQ,\name{\eta})
}$--preserving.\\
Assume that for $\ell<2$ we have $(\bQ_\ell,\name{\eta}_\ell)$ as in
\ref{6.0}, e.g.\ $\bQ_0$ is Cohen forcing, $\bQ_1$ is random real forcing. 
Let $\bQ =\{\emptyset\}\cup\bigcup\limits_{\ell<2}\{\ell\}\times\bQ_\ell$,
$\emptyset$ minimal, $(\ell_1,q_1) \le (\ell_2,q_2)$ iff $\ell_1 = \ell_2$ and
$\bQ_\ell\models q_1 \le q_2$. We define a $\bQ$--name $\name{\eta}$ by
defining for a generic $G\subseteq\bQ$ over $\V$:  
\[\name{\eta}[G]\mbox{ is }\begin{array}{ll}
\langle 0\rangle\conc(\eta_0[G_0])&\mbox{if }\{0\}\times\bQ_0\cap G\ne
\emptyset,\mbox{ and }G_0=\{q\in\bQ_0:(0,q)\in G\}\\
\langle 1\rangle\conc(\eta_1[G_1])&\mbox{if }\{1\}\times\bQ_1\cap G\ne
\emptyset,\mbox{ and }G_1=\{q\in\bQ_1:(1,q)\in G\}.
                                \end{array}\] 
Then usually (and certainly for our choice) we get a counterexample.

\begin{proposition}
\label{6.6x}
Assume that $A$ is a Borel subset (better: a definition of a Borel subset) of
$\baire$, $M$ is a $\bQ$--candidate (so $\name{\eta}\in M$, i.e.\ $\langle
\psi_{\alpha,m}:\alpha<\omega,m<\omega\rangle\in M$) and $A\in M$ (i.e.\ the
definition). Further, suppose that $q\in \bQ^M$ is such that
$q\forces_{\bQ}$``$\name{\eta}\in A$''. Then
\begin{enumerate}
\item[$(\alpha)$] $M\models$``$q\forces_{\bQ}\name{\eta}\in A$'',
\item[$(\beta)$]  there is $\eta\in A$ which is a $(\bQ,\name{\eta})$--generic
real over $M$. 
\end{enumerate}
\end{proposition}

\Proof As for $(\alpha)$, if it fails then for some $q'\in\bQ^M$, we have
\[M\models\mbox{`` }q \le^{\bQ} q'\mbox{ and }q'\forces_{\bQ}\name{\eta}
\notin A\mbox{ ''},\]
and let $r\in\bQ$ be $\langle M,\bQ\rangle$--generic above $q'$. So if $G$ is
a subset of $\bQ$ generic over $\V$ to which $r$ belongs then $q'\in G$ and
$G\cap\bQ^M$ is a subset of $\bQ^M$ generic over $M$ to which $q'$
belongs. Hence $M\langle G\rangle\models$``$\name{\eta}[G\cap\bQ^M]\notin A$''
and $\name{\eta}[G\cap Q^M]\in\baire$. By absoluteness also $\V[G]\models
\name{\eta}[G\cap\bQ^M]\notin A$ and $\name{\eta}[G\cap\bQ^M]\in\baire$. But
as $\name{\eta}\in M$ clearly $\name{\eta}[G\cap\bQ^M]=\name{\eta}[G]$ and as
$q'\in G$ also $q\in G$, so we get contradiction to
$q\forces_{\bQ}$``$\name{\eta}\in A$''.  

By clause $(\alpha)$ clause $(\beta)$ is easy: we can find a subset $G\in\V$
of $\bQ^N$ to which $q$ belongs which is generic over $M$. So $\name{\eta}[G]
\in\baire$ and it belongs to $A$ as $M\models$``$q\forces_{\bQ}\name{\eta}\in
A$''. \QED$_{\ref{6.6x}}$

\begin{proposition}
\label{6.7}
Assume $\bQ$ is correct and satisfies the c.c.c. The following conditions are
equivalent for a set $X\subseteq\baire$: 
\begin{enumerate}
\item[(A)] $X\in I^{\ex}_{(\bQ,\name{\eta})}$,
\item[(B)] for some $\rho\in\can$, for every $\bQ$--candidate $N$ to which
$\rho$ belongs there is no $\eta\in X$ which is $(\bQ,\name{\eta})$--generic
over $N$,
\item[(C)] for every $p\in\bQ$ for some $\bQ$--candidate $N$ such that $p\in
\bQ^N$, there is no $\eta\in X$ which is $(\bQ,\name{\eta})$--generic over
$N$. 
\end{enumerate}
\end{proposition}

\Proof (A)$\ \Rightarrow\ $(B):\quad So assume (A), i.e.\ $X\in I^{\ex}_{(\bQ,
\name{\eta})}$. Then for some Borel set $A\in I_{(\bQ,\name{\eta})}$ we have
$X\subseteq A$. Let $\rho\in\can$ code $A$. Since
$\forces_{\bQ}$``$\name{\eta}\notin A^{\V[\name{G}_{\bQ}]}$'', it follows
from \ref{6.6x} that
\begin{enumerate}
\item[(*)] for any $\bQ$--candidate $N$ to which $\rho$ belongs there is no
$(\bQ,\name{\eta})$--generic real $\eta$ over $N$ which belongs to $X$ (or
even just to $A$).
\end{enumerate}

\noindent (B)$\ \Rightarrow\ $(C):\quad Easy as $Q$ is correct.

\noindent (C)$\ \Rightarrow\ $(A):\quad Assume (C). Let 
\[\begin{array}{ll}
{\cal I}=\{p\in\bQ:&\mbox{for some Borel subset } A=A_p\mbox{ of }\baire\\
\ &\mbox{we have }p\forces\mbox{`` }\name{\eta}\notin A_p\mbox{ '' and }
X\subseteq A_p\}.
  \end{array}\]
Suppose first that ${\cal I}$ is predense in $\bQ$. Clearly it is open, so we
can find a maximal antichain ${\cal J}$ of $\bQ$ such that ${\cal J}\subseteq
{\cal I}$. As $\bQ$ satisfies the c.c.c., necessarily ${\cal J}$ is countable.
So $A\stackrel{\rm def}{=}\bigcap\limits_{p\in {\cal J}} A_p$ is a Borel
subset of $\baire$ (as ${\cal J}$ is countable) and it includes $X$ (as each
$A_p$ does). Moreover, since ${\cal J}$ is a maximal antichain of $\bQ$ (and
$p\in {\cal J}\ \Rightarrow\ p\in {\cal I}\ \Rightarrow\ p\forces_{\bQ}
$``$\name{\eta}\notin A_p\mbox{''}\ \Rightarrow\ p\forces_{\bQ}$``$\name{\eta}
\notin A$'') we have $\forces_{\bQ}$``$\eta\notin A$''. Consequently (A) holds. 

Suppose now that ${\cal I}$ is not predense in $\bQ$ and let $p^*\in\bQ$
exemplifies it, i.e.\ it is incompatible with every member of ${\cal I}$. Let
$N$ be a $\bQ$--candidate to which belongs some $\rho$ given by the assumption
(C) for $p^*$. Thus $p^*\in\bQ^N$ and no $\eta\in X$ is $(\bQ,\name{\eta
})$--generic over $N$. Let $q$ be a member of $\bQ$ which is above $p^*$ and
is $\langle N,\bQ^N\rangle$--generic (i.e.\ $q\forces$``$G^P\cap\bQ^N$ is
generic over $N$''). Let $A\stackrel{\rm def}{=}\{\eta\in\baire:\eta$ is not
$(\bQ,\name{\eta})$--generic over $N\}$. Now
\begin{enumerate}
\item[(a)] $A$ is a Borel subset of $\baire$ and $X\subseteq A$

\noindent (why? as $N$ is countable),
\item[(b)] $q\forces_{\bQ}$``$\name{\eta}\notin A^{\V[G_{\bQ}]}$''

\noindent (why? by the definition of $A$),
\item[(c)] $q\in {\cal I}$

\noindent (why? by (a)+(b)).
\end{enumerate}
Thus $p^*\le q\in {\cal I}$ and we get contradiction to the choice of
$p^*$. \QED$_{\ref{6.7}}$ 

\begin{theorem}
\label{6.5}
Assume that:
\begin{enumerate}
\item[(a)] $\bQ$, $\name{\eta}$ are as above (see\ref{6.0}), and $\bQ$ is
correct, 
\item[(b)] $\bP$ is nep-forcing notion with respect to {\em our fixed}
version $\ZFC^-_*$, 
\item[(c)] $\bP$ is $I_{(\bQ,\name{\eta})}$--preserving,
\item[(d)] $\ZFC^-_{**}$ is a stronger version of set theory including
clauses (i)--(v) below for some $\chi_1 <\chi_2$,
\begin{enumerate}
\item[(i)]   $({\cal H}(\chi_2),\in)$ is a (well defined) model of $\ZFC^-_*$,
\item[(ii)]  (a), (b) and (c) (with ${\frak B}^{\bP}$, ${\frak B}^{\bQ}$,
$\name{\eta}$ as individual constants),
\item[(iii)] $\bQ,\bP\in {\cal H}(\chi_1)$ and $({\cal H}(\chi_2),\in)$ is a
semi $\bP$--candidate and a semi $\bQ$-candidate with $({\frak B}^{\bP})$
interpreted as $({\frak B}^{\bP})^N\rest {\cal H} (\chi_2)^N$ and similarly
for $\bQ$, so (natural to assume) ${\frak B}^{\bP}$, ${\frak B}^{\bQ}\in {\cal
H}(\chi_2)$,

(remember, ``semi'' means omitting the countability demand)
\item[(iv)] forcing of cardinality $<\chi_1$ preserves the properties (i),
(ii), (iii), and $\chi_1$ is a strong limit cardinal,
\item[(v)] forcing by $\bP$ preserves ``${\cal I}$ is a predense subset of
$\bQ$'' (follows if $\bQ$ satisfies the c.c.c.\ by \ref{5.5A}(ii)). 
\end{enumerate}
\end{enumerate}
{\em Then}:
\begin{enumerate}
\item[$(\alpha)$] {\em if}, additionally, 
\begin{enumerate}
\item[(e)] $\ZFC^-_{**}$ is normal (see Definition\ref{0.9}(3))
\end{enumerate}
{\em then} $\bP$ is strongly $I_{(\bQ,\name{\eta})}$--preserving,
\item[$(\beta)$]  {\em if}\ $N$ is a $\bP$--candidate (and $\bQ$--candidate)
and moreover it is a model of $\ZFC^-_{**}$ and $N\models$``$p\in\bP$'' and
$\eta^*$ is $(\bQ,\name{\eta})$--generic over $N$,\\
{\em then}\ for some $q$ we have:
\begin{enumerate}
\item[(i)]   $p\le q$ and $q\in\bP$,
\item[(ii)]  $q$ is $\langle N,\bP\rangle$--generic; i.e.\ $q\forces_{
\bP}$``$\name{G}_{\bP}\cap\bP^N$ is generic over $N$'' (see \ref{2.3}), 
\item[(iii)] $q\forces_{\bP}$`` $\eta^*$ is $(\bQ,\name{\eta})$--generic over
$N[\bP^N\cap\name{G}_{\bP}]$ ''.
\end{enumerate}
\item[$(\alpha)^+$] We can strengthen the conclusion of $(\alpha)$ to 

``$\bP$ is super--$I_{(\bQ,\name{\eta})}$--preserving''.
\end{enumerate}
\end{theorem}

\begin{remark}
\label{7.8A}
1)\quad We consider, for a nep forcing notion $\bQ$
\begin{enumerate}
\item[$(*)_1$] $\bQ$ satisfies the c.c.c.
\end{enumerate}
We also consider  
\begin{enumerate}
\item[$(*)_2$] being a predense subset (or just a maximal antichain) of $\bQ$
is $K$--absolute.
\end{enumerate}
By results of the previous section, $(*_1)\ \Rightarrow\ (*_2)$ under
reasonable conditions. You may wonder whether $(*_2)\ \Rightarrow\ (*_1)$, but
by the examples in section 11 the answer is not.

\noindent 2)\quad Note that in $(\alpha), (\alpha)^+$ we can use the weak
normality if $\bQ$ satisfies the c.c.c., see \ref{6.8}. We do not use ``$\bP$
is explicitly nep'' so we do not demand it.
\end{remark}

Before we prove the theorem, let us give an example for a forcing notion
failing the conclusion and see why many times we can simplify assumptions.

\begin{example}
\label{6.6}
Start with $\V_0$.  Let $\bar{s}=\langle s_i:i<\omega_1\rangle$ be a sequence
of random reals, forced by the measure algebra on ${}^{\omega_1}(\can)$. Let 
$\V_1=\V_0[\bar{s}]$, $\V_2=\V_1[r]$, $r$ a Cohen over $\V_1$ and 
\[\begin{array}{ll}
\V_3=\V_2[\bar{t}]&\mbox{ where }\bar{t}=\langle t_i:i<\omega_1\rangle 
\mbox{ is a sequence of random reals}\\
\ &\mbox{forced by the measure algebra}.
 \end{array}\]
Then in $\V_3$ (in fact, already in $\V_2$), $\{s_i:i<\omega_1\}$ is a null
set, whereas $\{t_i:i<\omega_1\}$ is not null. But $\bar{t}$ is also generic
for the measure algebra over $\V_1$. So $\V'_2=\V_1[\bar{t}]$ is a generic 
extension of $\V_1$. We have $\V_3=\V'_2[r]$, where $r$ is generic for some
algebra, more specifically for 
\[\bR\stackrel{\rm def}{=}(\mbox{Cohen $*$ measure algebra adding }\bar{t})/
\bar{t}.\]
So in $\V'_2$ the sets $\bar{t}$ and $\bar{s}$ are not null and $\bR$ makes
$\bar{s}$ null, but not $\bar{t}$. 

How can $\bR$ do that?  $\bR$ uses $\langle t_i:i<\omega_1\rangle$ in its 
definition, so it is not ``nice'' enough.\QED$_{\ref{6.6}}$
\end{example}

\noindent{\bf Remark}\qquad In the proof of \ref{6.5}, of course, we may
assume $N\prec ({\cal H}(\chi,\in))$ if $({\cal H}(\chi,\in))\models
\ZFC^-_{**}$, as this normally holds. In $(\alpha)$ the use of such $N$ does
not matter. In $(\beta)$ it slightly weakens the conclusion. Now, $(\alpha)$ is
our original aim. But $(\beta)$ both is needed for $(\alpha)$ and is a step
towards preserving them (as in \cite{Sh:f}). So typically $N$ is an elementary
submodel of appropriate ${\cal H}(\chi)$. 
\medskip

{\noindent{\sc Proof of \ref{6.5}}\hspace{0.2in}} {\bf Clause $(\alpha)$}:
\qquad To prove $(\alpha)$ we will use $(\beta)$. So let $X\subseteq\baire$, 
$X\in (I^{\dx}_{(\bQ,\name{\eta})})^+$. Then there is a condition $q^*\in\bQ$
such that 
\begin{enumerate}
\item[$(*)_1$] for no Borel subset $B$ of $\baire$ do we have:\quad
$X\subseteq B$ and $q^*\forces_\bQ$``$\name{\eta}\notin B$''.
\end{enumerate}
Let $\chi$ be large enough. We can find $N\subseteq ({\cal H}(\chi),\in)$ as
in $(\beta)$, moreover $N\prec ({\cal H}(\chi),\in)$ a model of $\ZFC^-_{**}$
(and so a $\bP$--candidate and a $\bQ$-candidate) [it exists because by clause
(e) of the assumptions, $\ZFC^-_{**}$ is normal so for $\chi$ large enough any
countable $N\prec ({\cal H}(\chi),\in)$ to which ${\frak C},{\frak B}^{\bQ},
{\frak B}^{\bP}$ belong is a model of $\ZFC^-_{**}$ and is a $\bP$--candidate
and a $\bQ$--candidate, so as required].  

Towards a contradiction, assume $p^*\in\bP$ and $p^*\forces_{\bP}$``$X\in
I^{\dx}_{(\bQ,\name{\eta})}$''. So for some $\bP$--name $\name{A}$ we have 
\[p^*\forces_{\bP}\mbox{`` }\name{A}\mbox{ is a Borel subset of }\baire,\ X
\subseteq\name{A}\mbox{ and }\name{A}\in I_{(\bQ,\name{\eta})},\mbox{
i.e. }\forces_{\bQ}\name{\eta}\notin\name{A}\mbox{ ''.}\]  
Without loss of generality the name $\name{A}$ is hereditarily countable and
$\name{A},p^*,q^*$ belong to $N$. In $\V$, let 
\[\begin{array}{ll}
B=\{\eta\in\baire:&\eta\mbox{ is a $(\bQ^{\geq q^*},\name{\eta})$--generic
real over $N$, which means:}\\
\ &\eta=\name{\eta}[G]\mbox{ for some $G\subseteq\bQ^N$ generic over }N\\
\ &\mbox{such that }q^*\in G\}
  \end{array}\]
Clearly, it is an analytic set (if $\name{\eta}$ was generic real then Borel;
both holds as ``$\name{\eta}$ is a generic real for $\bQ$'' follows from
$\ZFC^-_{**}$). So $B=\bigcup\limits_{i<\omega_1}B_i$, each $B_i$ is Borel. 
Let $q\in\bQ$ be $\langle N,\bQ\rangle$--generic and $q^*\leq q$. Then
$q\forces_{\bQ}$``$\name{\eta}\in B$'' and hence wlog for some $i<\omega_1$ we
have $q\forces_\bQ$``$\name{\eta}\in B_i$''. Since $q\forces_\bQ$``$\name{
\eta}\notin (\baire\setminus B_i)$'' (as $q^*\leq q$, $q^*\forces$``$\name{
\eta}\in B_i$''), we may apply $(*)_1$ to the set $\baire\setminus B_i$ to
conclude that $X\not\subseteq\baire\setminus B_i$. Take $\eta^*\in X\cap
B_i$ (so it is $(\bQ,\name{\eta})$--generic over $N$). So by clause $(\beta)$
(proved below), there is a condition $p\in\bP$, $p\ge p^*$ which is $\langle
N,\bP\rangle$--generic (i.e.\ it forces that $\name{G}_{\bP}\cap\bP^N$ is
generic over $N$, not necessarily $\name{G}_{\bP}\cap N$) and such that 
\[p\forces_\bP\mbox{`` }\eta^*\mbox{ is }(\bQ,\name{\eta})\mbox{--generic over }
N[\name{G}_{\bP}\cap \bP^N]\mbox{ ''.}\]
Choose $G_{\bP}\subseteq\bP$, generic over $\V$, such that $p\in G_{\bP}$. In
$\V[G_{\bP}]$, $N[G_{\bP}\cap\bP^N]$ is a generic extension of $N$ (for
$\bP^N$!), a $\bQ$ candidate, and $\eta^*$ is $(\bQ,\name{\eta})$--generic
over it. As $p^*\leq p\in G_{\bP}$, clearly if $G_{\bQ}\subseteq\bQ^{\V[G_{
\bP}]}$ is generic over $\V[G_{\bP}]$ then $\name{\eta}[G_{\bQ}]\notin\name{A}
[G_{\bP}]$. But $N[G_{\bP}\cap \bP^N]\prec({\cal H}(\chi)^{\V[G_{\bP}]},\in)$,
so $N[G_{\bP}\cap \bP^N]$ satisfies the parallel statement. Since $\eta^*$ is
$(\bQ,\name{\eta})$--generic over $N[G_{\bP}\cap\bP^N]$, it cannot belong to
$\name{A}[G_{\bP}\cap \bP^N]$. But easily $\name{A}[G_{\bP}\cap\bP^N]=
\name{A}[G_{\bP}]$ and hence, by absoluteness, $\eta^*\in X\subseteq
\name{A}[G_{\bP}]$, a contradiction. This ends the proof of \ref{6.5}, clause
$(\alpha)$. 
\medskip

\noindent {\bf Clause $(\alpha)^+$}:\qquad Like the proof of clause
$(\alpha)$. We start like there but now we choose functions $r^*,A^*, {\cal
I}$ such that 
\begin{enumerate}
\item[$(*)_2$] $\dom(r^*)=\dom({\cal I})$ is the set of all hereditarily
countable canonical $\bP$--names for elements of $\bQ$ (so it is a member of
${\cal H}_{<\aleph_1}(\kappa(\bP)+\kappa(\bQ))$), and
 
$\dom(A^*)=\{(p,\name{q}): p\in {\cal I}(\name{q}),\ \name{q}\in\dom(r^*)\}$,
\item[$(*)_3$] for each $\name{q}\in\dom(r^*)=\dom({\cal I})$, ${\cal
I}(\name{q})$ is a predense subset of $\bP$ such that for each $p\in{\cal
I}(\name{q})$ we have:
\[\begin{array}{ll}
p\forces_{\bP}\mbox{``}&A^*(\name{q})\mbox{ is a Borel subset of }\baire\mbox{
''},\\ 
p\forces_{\bP}&[r^*(\name{q})\forces_{\bQ}\mbox{`` }\name{\eta}\notin A^*
(\name{q})\mbox{ ''}],\\
p\forces_{\bP}\mbox{``}&X\subseteq A^*(\name{q})\mbox{ ''}.
  \end{array}\]
\end{enumerate}
Without loss of generality, the set $X$,  and the functions $r^*,A^*,{\cal I}$
belong to $N$. We choose conditions $q\in\bQ$, $p\in\bP$ and a real $\eta^*\in
X$ and a generic filter $G_\bP\subseteq\bP$ over $\V$ in a similar manner as
in clause $(\alpha)$. We note that  
\[\name{q}\in\dom(r^*)\cap N\quad \Rightarrow\quad N\cap{\cal I}(\name{q})\cap
G_\bP\neq\emptyset,\]
so say $p[\name{q}]\in G_\bP\cap N$. Since $\eta^*$ is
$(\bQ,\name{\eta})$--generic over $N[G_\bP\cap\bP^N]$, there is $G^*\subseteq
\bQ^{N[G_\bP\cap \bP^N]}$ generic over $N$ such that $\eta^*=\name{\eta}[G^*]$. 
By the choice of $r^*,A^*$ there is $\name{q}\in N\cap\dom(r^*)$ such that
$r^*[\name{q}][G_\bP\cap \bP^N]\in G^*$. Now, $A=A^*(p[\name{q}],\name{q})\in
N[G_\bP\cap \bP^N]$ is a Borel subset of $\baire$ and $N[G_\bP\cap\bP^N]
\models$``$\name{\eta}\notin A$'', hence $N[G_\bP\cap \bP^N]\models$``$\name{
\eta}[G^*]\notin A$. But
\[N[G_\bP]=N[G_\bP\cap\bP^N]\models\mbox{`` }X\setminus A=\emptyset\mbox{ ''},\]
contradicting $\eta^*=\name{\eta}[G^*]\in X\setminus A$. 
\medskip

\noindent {\bf Clause $(\beta)$}:\qquad So $N,\eta^*,\bQ,\bP,p$ are given. Let
$N_1=N[G^*]$ be a generic extension of $N$ by a subset $G^*$ of $\bQ^N$
generic over $N$ and such that $\eta^*=\name{\eta}[G]$ (see\ref{5.2A}). Now
choose (in $N$) a model $M\prec ({\cal H}(\chi_2),\in)^N$ such that 
\begin{enumerate}
\item[(i)]   $\bP,\bQ,\name{\eta},p\in M$,
\item[(ii)]  $\bQ^N\subseteq M$ and $\bP^N\subseteq M$,
\item[(iii)] the family of maximal antichains of $\bP$ and of $\bQ$ from $N$
are included in $M$,
\item[(iv)]  $M\in N$, moreover $M\in {\cal H}(\chi_2)^N$,
\item[(v)]   $M\models$``forcing by $\bP$ preserves predensity of subsets
of $\bQ$''
\end{enumerate}
[Why is clause (v) possible? As $N\rest {\cal H}(\chi_1)$ inherits clause (v)
of (d) of the assumptions].  

Hence, by assumption (d),
\begin{enumerate}
\item[$(\bigotimes)$] $M$ is a $\bP$-candidate and a $\bQ$-candidate and 
\[N\models\mbox{`` }M\mbox{ is a semi $\bP$-candidate and semi
$\bQ$-candidate ''}.\]
\end{enumerate}
Let $\bR=\Levy(\aleph_0,|M|)$. In $\V$ let $G_{\bR}\subseteq \bR$ be generic
over $N_1=N[G^*]$ (note that as $N_1$ is countable, clearly $G_{\bR}$ exists)
and let $N_2=N_1[G_{\bR}]$ (note that it too is a $\bP$-candidate and a
$\bQ$-candidate).  

Note: $\eta^*$ is $(\bQ,\name{\eta})$--generic over $M$ too and $G^*$ is a
subset of $\bQ^M$ generic over $M$ (by clauses (ii) + (iii)) and
$\bQ^N=\bQ^M$, $\bP^N=\bP^M$ (note that in $N_2$ the model $M$ is countable).

Now we ask the following question: 
\begin{quotation}
Is there $q\in\bP^{N_2}$ such that 

\noindent $N_2\models$`` $p\le^{\bP} q$, $q$ is $(M,\bP^M)$--generic and
$q\forces_{\bP}$``$\eta^*$ is $(\bQ,\name{\eta})$--generic over $M[G_{\bP}\cap
\bP^M]$'' ''?
\end{quotation}
Depending on the answer, we consider two cases.
\smallskip 

\noindent{\bf Case 1}:\quad The answer is ``yes''. \\
Choose $q'\in \bP$, $q'\ge q$, $q'$ is $(N_2,P^{N_2})$-generic. Then we have
\[\begin{array}{ll}
q'\forces_{\bP}\mbox{``}&\mbox{in }\V[\name{G}_{\bP}],\ \name{G}_{\bP}\cap
\bP^{N_2}\mbox{ is generic over }N_2,\ p,q\in \name{G}_{\bP},\mbox{ and}\\
\ &\mbox{in }N_2[\name{G}_{\bP}\cap \bP^{N_2}],\ \eta^*\mbox{ is $(\bQ,
\name{\eta})$--generic over }M[\name{G}_{\bP}\cap\bP^M],\\
\ &\mbox{hence also over }N[\name{G}_{\bP}\cap \bP^N]\mbox{ ''.}
  \end{array}\]
[Why does $q'$ force this? As:
\begin{enumerate}
\item[(A)] ``$\name{G}_{\bP}\cap \bP^{N_2}$ is generic over $N_2$'' holds
because $q'$ is $(N_2,\bP^{N_2})$--generic; 
\item[(B)] ``$p,q\in\name{G}_{\bP}$'' holds as $p\le q\le q'\in G_{\bP}$
(forced by $q'$!); 
\item[(C)] ``in $N_2[\name{G}_{\bP}\cap \bP^{N_2}]$, $\eta^*$ is $(\bQ,
\name{\eta})$--generic over $M[G_{\bP}\cap \bP^M]$'' holds because of the
choice of $q$ (i.e.\ the assumption of the case ad as $q\in G_\bP$); 
\item[(D)] ``$\eta^*$ is $(\bQ,\name{\eta})$--generic over $N[G_{\bP}\cap
\bP^N]$ for $\bQ$'' holds by clause (C) above and clause (iii) of the choice
of $M$.] 
\end{enumerate}
By absoluteness we can omit the ``in $N_2[\name{G}_{\bP}\cap\bP^M]$'', i.e.\
$q'\forces_{\bP}$``$\eta^*$ is $(\bQ,\name{\eta})$--generic over
$N[\name{G}_{\bP}\cap \bP^N]$''. So $q'$ is as required.
\smallskip

\noindent{\bf Case 2}:\quad The answer is ``no''. \\
Let $\psi(x)$ be the following statement:
\begin{quotation}
there is no $q$ such that:\\
$q\in \bP$, $\bP\models$``$p\le q$'', $q$ is $(M,\bP^M)$--generic and $q
\forces_{\bP}$``$x$ is a $(\bQ,\name{\eta})$--generic real over
$M[\name{G}_{\bP}\cap \bP^M]$''.
\end{quotation}
So $\psi$ is a first order formula in set theory, all parameters are in $N_1=
N[G^*]\subseteq N_2=N[G^*][G_{\bR}]$, and by the assumption of the case
\[N[G^*][G_{\bR}]\models\psi[\eta^*].\]
Since $\bR$ is homogeneous we may assume that $r=\emptyset$. As $G_{\bR}
\subseteq\bR$ is generic over $N[G^*]$ for $\bR$, necessarily (by the forcing
theorem), for some $r\in G_{\bR}$  
\[N[G^*]\models\mbox{`` }r\forces_{\bR}\psi[\eta^*]\mbox{ ''.}\]
So necessarily, for some $q\in G^*\subseteq \bQ^N=\bQ^M$ we have
\[N\models\bigl(q\forces_{\bQ} [r\forces_{\bR}\psi(\name{\eta}[\name{G}_{\bQ}
])]\bigr).\]
Now $\bR\in N$ (as it members are finite sets of pairs of ordinals) so
\begin{enumerate}
\item[$(\otimes)$]  $N\models\bigl((q,r)\forces_{\bQ\times \bR}\psi(
\name{\eta}[\name{G}_{\bQ}])\bigr)$.
\end{enumerate}
Next, $N[G_{\bR}]$ is a generic extension by a ``small'' forcing of $N$ which
is a model of $\ZFC^-_{**}$, so $N[G_{\bR}]$ satisfies (i), (ii) and (iii) of
the clause (d) of the assumptions. Note that $N\models$``$M$ is a semi
$\bQ$-candidate and a semi $\bP$--candidate'', see clause (d)(iii) of the
assumptions and the choice of $M$, so also $N[G_{\bR}]$ satisfies this. 
Moreover, $N[G_\bR]\models$``$M$ is countable'', so $N[G_\bR]\models$``$M$ is
a $\bQ$-candidate and a $\bP$--candidate''. Hence by assumption (d)(ii) there are $p_1,\eta^\otimes,G^\otimes_\bQ\in N[G_{\bR}]$ such that:
\[\begin{array}{ll}
N[G_{\bR}]\models\mbox{``}&p_1\in\bP,\ p\le^{\bP} p_1,\ p_1\mbox{ is $(M,
\bP^M)$--generic and}\\
\ &p_1\forces_{\bP}[\eta^\otimes\mbox{ is a $(\bQ,\name{\eta})$--real over
$M[\name{G}_{\bP}\cap \bP^M]$ satisfying }q]\\
\ &\mbox{moreover, }p_1\forces_\bP\mbox{`` }\eta^\otimes=\name{\eta}[
G^\otimes_\bQ]\mbox{ '', and}\\
\ &G^\otimes_\bQ\subseteq \bQ^{M[G_\bP\cap\bP^M]}\mbox{ is a generic set over
$M$ such that $q\in G^\otimes_\bQ$  ''}.
  \end{array}\]
[Here we use the following:\quad if $G\subseteq\bP^{N[G_\bR]}$ is generic over
$N[G_\bR]$ then $N[G_\bR]\langle G\rangle$ is a $\bQ$--candidate (apply clause
(iv) of the assumption (d) to $\bR*\name{\bP}$).] It follows from clause
(d)(v) of the choice of $M$ that 
\begin{quotation}
$G^\otimes_{\bQ}\cap \bQ^M$ is generic over $M$.
\end{quotation}
Let $\name{p}^1,\name{\eta}^\otimes,\name{G}^\otimes_\bQ\in N$ be $\bR$--names
such that $\name{\eta}^\otimes[G_{\bR}]=\eta^\otimes$, $G^\otimes_\bQ=
\name{G}^\otimes_\bQ[G_\bR]$ and $\name{p}^1[G_{\bR}]=p_1$, and without loss
of generality in $N$ we have 
\[\begin{array}{ll}
r\forces_{\bR}\mbox{``}&\name{\eta}^\otimes\mbox{ is a
$(\bQ,\name{\eta})$--generic real over $M$ satisfying $q$ and }\name{p}^1\in
\bP\mbox{ and}\\
\ &\name{p}^1\mbox{ forces $(\forces_{\bP})$ that $\name{\eta}^\otimes$ is
$(\bQ,\name{\eta})$--generic over }M[\name{G}_{\bP}\cap \bP^M],\\
\ &\name{G}^\otimes_\bQ\mbox{ is a subset of $\bQ^M$ generic over $\V$ ''.}
  \end{array}\]
Let $G'_{\bR}\subseteq \bR$ be generic over $N[G_{\bR}]$, to which $r$ belongs,
so $N[G_{\bR}][G'_{\bR}]$ is a forcing extension of $N[G_{\bR}]$ so both are
generic extensions of $N$ by a small forcing.

Now $\name{G}^\otimes_\bQ$ is essentially a complete embedding of
$\bQ\rest(\ge q)$ into $\bR$ (by basic forcing theory, see the footnote to
\ref{0.9}(1)(d); and we can use the value for 0 of the function $\bigcup\{f:
f\in G_{\bR}\}$ to choose $q'$, $q\le q'\in \bQ^N$). Hence, for some
$\bQ$--name $\name{\bR}^*$ we have $(\bQ\rest(\ge q))*\name{\bR}^*$ is $\bR$,
so $G_{\bR}=G^\otimes_{\bQ} * G_{\bR^*}$ for some $G_{\bR^*}\in \V[G_{\bR}]$,
where $\bR^*=\name{\bR}^*[G^\otimes_{\bQ}]$, $\name{\bR}^*$ a $\bQ$--name. So
we can represent $N[G_{\bR}][G'_{\bR}]$ also as $N^3\stackrel{\rm def}{=}N[
G^\otimes_{\bQ}][G'_{\bR}][G_{\bR^*}]$; i.e.\ forcing first with $\bQ\rest
(\geq q)$, then with $\bR$, lastly with $\name{\bR}^*[G^\otimes_{\bQ}]$. Now
let $N^2\stackrel{\rm def}{=} N[ G^\otimes_{\bQ}][G'_{\bR}]$, so $N^2$ is a
generic extension of $N$ and $N^3$ is a generic extension of $N^2$ (both by
``small'' forcing), and in $N^3$ we have $p^1$ and $\eta^\otimes$ and
$G^\otimes_\bQ$. But $G^\otimes_{\bQ}\times G'_{\bR}$ is a generic subset of
$(\bQ^N\rest\geq q)\times \bR$ over $N$, so essentially a generic (over $N$)
subset of $\bQ^N\times \bR$ to which $(q,r)$ belongs, hence (by $(\otimes)$
above) $N^2\models\psi(\name{\eta}[G^\otimes_{\bQ}])$. Therefore there is no
$p'\in N^2$ such that\footnote{of course, we can use a weaker demand on
$G^\otimes_\bQ$}:   
\begin{enumerate}
\item[$(\boxtimes)$] \quad $N^2\models[p'\in \bP$, $p\le p'$, $p'\forces_{\bP
}$``$\name{\eta}[G^\otimes_{\bQ}]$ is $(\bQ,\name{\eta})$--generic over
$M[\name{G}_{\bP}\cap \bP^M]$'']. 
\end{enumerate}
In $N^2$ we can find a countable $M'\prec ({\cal H}(\chi_2)^{N^2},\in)$ to
which $\bP$, $\bQ$, $\name{\eta}^\otimes$, $\name{p}^1$, $\name{\bR}^*$, and
$G^\otimes_{\bQ}$, $G'_{\bR}$, $M$ belong (so $N^2\models$`` $M'$ is countable
and is a $\bP$--candidate and a $\bQ$-candidate'') and $M''=M'\rest N\in N$. 
In $N^2$ we can find $G'_{\name{\bR}^*[G^\otimes_{\bQ}]}\subseteq\name{\bR}^*[
G^\otimes_{\bQ}]\cap M'$ which is generic over $M'$. In $M^3=M''[G^\otimes_\bQ]
[G'_{\name{\bR}^*[G^\otimes_{\bQ}]}]$, again by the forcing theorem, there is
$p'$ as required in $(\boxtimes)$ above. So $M^3$ is a $\bP$--candidate inside
$N^2$, hence there is $p'_1$ such that $N^2\models$``$p'\le p'_1$ and $p'_1$ is
$(M^3,\bP)$--generic''. By the amount of absoluteness we require (moving up
from $M^3$ to $N^2$) this $p'_1$ can serve in $N^2$ for $(\boxtimes)$,
contradiction to the previous assertion. \QED$_{\ref{6.5}}$  

\begin{proposition}
\label{6.8}
Assume (a),(b),(c) and (d) of \ref{6.5} and
\begin{enumerate}
\item[(e)] $\ZFC^-_{**}$ is weakly normal,
\item[(f)] $\bQ$ is c.c.c. and simple (for simplicity) and correct.
\end{enumerate}
{\em Then} $\bP$ is strongly $I_{(\bQ,\name{\eta})}$--preserving.
\end{proposition}

\Proof First note that $I^\ex_{(\bQ,\name{\eta})}=I^\dx_{(\bQ,\name{\eta})}$.
Now, if the conclusion fails as witnessed by a set $X$, then, by \ref{6.7},
the statements (A), (B), (C) of \ref{6.7} fail. Hence, by $\neg$(A), $X\in
(I^{\ex}_{(\bQ,\name{\eta})})^+$ and $p\in\bP$ and a $\bP$--name $\name{y}$
such that  
\[\begin{array}{ll}
p\forces_{\bP}\mbox{``}&\name{y}\in {\cal H}_{<\aleph_1}(\bQ)\mbox{ and for
no $Q$-candidate $M$ such that $\name{y}\in M$}\\
\ &\mbox{there is }\nu\in X\mbox{ which is $(\bQ,\name{\eta})$--generic over
}N[\name{G}_{\bP}]\mbox{ ''.}
  \end{array}\]
As we can increase $p$, without loss of generality $\name{y}$ is a
hereditarily countable $\bP$--name. As $\ZFC^-_{**}$ is weakly normal we can
find a model $N$ of $\ZFC^-_{**}$ which is a $\bP$--candidate and a
$\bQ$--candidate and to which $p,\name{y}$ belongs. Let $\eta^*\in X\subseteq
\baire$ (in $\V$) be $(\bQ,\name{\eta})$--generic over $N$ (exists by the
negation of (B) of \ref{6.7}). By $(\beta)$ of \ref{6.5} there is $q\in\bP$
such that $p\le q$, $q$ is $\langle N,\bP\rangle$-generic and $q\forces_{
\bP}$``$\eta^*$ is $(\bQ,\name{\eta})$--generic over $N[\name{G}_{\bP}]$'', a
contradiction. \QED 

\begin{proposition}
\label{6.9}
Assume (a), (b), (c)  of \ref{6.5}. Let $\bP,\bQ$ be normal and forcing with
$\bP$ preserves ``${\cal I}\subseteq\bP$ is predense''. Then
\begin{enumerate}
\item[$(\alpha)'$] $\bP$ is strongly $I_{(\bQ,\name{\eta})}$--preserving,
\item[$(\beta)$]   for $\chi$ large enough, if $N\prec ({\cal H}(\chi),\in)$
is countable (and ${\frak C},{\frak B}^{\bQ},{\frak B}^{\bP},\name{\eta}\in
N)$ and $N\models$``$p\in\bP$'' and $\eta^*$ is $(\bQ,\name{\eta})$--generic
over $N$ then for some $q$ we have
\begin{enumerate}
\item[(i)]   $p\le q$, $q\in\bP$,
\item[(ii)]  $q$ is $\langle N,\bP\rangle$--generic; i.e.\
$q\forces_{\bP}$``$\name{G}_{\bP}\cap\bP^N$ is generic over $N$'', 
\item[(iii)] $q\forces_{\bP}$``$\eta^*$ is $(\bQ,\name{\eta})$--generic over
$N[\bP^N\cap G_{\bP}]$''. 
\end{enumerate}
\end{enumerate}
\end{proposition}

\Proof Let $\chi_1$ be a large enough strong limit, and $\chi_2=\beth_\omega(
\chi_1)$, $\chi=\beth_\omega(\chi_2)$, and repeat the proof of \ref{6.5} using
$N\prec ({\cal H}(\chi_3),\in)$ to which ${\frak C},{\frak B}^{\bQ},{\frak
B}^{\bP}$ and $\bP,\bQ,\theta,\chi_1,\chi_2$ belong. \QED$_{\ref{6.9}}$

\begin{proposition}
\label{6.10}
Assume (a),(b),(c) of \ref{6.5} and
\begin{enumerate}
\item[(d)$'$] $\ZFC^-_{**}$ is a version of set theory including, for some
$\chi_1 <\chi_2$ 
\begin{enumerate}
\item[(i)]   $({\cal H}(\chi_2),\in)$ is a (well defined) model of $\ZFC^-_*$,
\item[(ii)]  (a) and (b) (with ${\frak B}^{\bP},{\frak B}^{\bQ},\name{\eta}$
as individual constants) and (c),
\item[(iii)] $\bQ,\bP\in {\cal H}(\chi_1)$ and $({\cal H}(\chi_2),\in)$ is a
$\bP$--candidate and a $\bQ$--candidate with ${\frak B}^{\bP}$ interpreted as 
$({\frak B}^{\bP})^N\rest {\cal H}(\chi_2)$ and similarly for $\bQ$, so
\[\mbox{`` }{\frak B}^{\bP},{\frak B}^{\bQ}\in {\cal H}(\chi_2)\mbox{ ''},\] 
\item[(iv)]  forcing with $\bP$ preserves being a $\bQ$--candidate,
\item[(v)]   $\bQ$ satisfies the c.c.c.\ and being incompatible in $\bQ$ is
upward absolute from $\bQ$--candidates. 
\end{enumerate}
\end{enumerate}
{\em Then} for models of $\ZFC^-_{**}$, forcing with $\bP$ preserves ``${\cal
I}$ is a predense subset of $\bQ$'' (i.e.\ (d)(v) of\ref{6.5}).
\end{proposition}

\noindent{\bf Remark:}\qquad 1)\quad When $\ZFC^-_{**}$ is normal, this
applies to \ref{9.6}.\\
2)\quad Compare with \ref{5.5A}(iii). Here we have redundant assumptions as we
have a use for \ref{6.5} in mind.

\Proof Let $N$ be a model of $\ZFC^-_{**}$ and let $N\models$``${\cal I}$ is a
predense subset of $\bQ$''. As $N\models$``$\bQ$ satisfies the c.c.c.'' we can
find in $N$ a set ${\cal J}\subseteq\{q:N\models(\exists p \in {\cal I})(p\le
q\in \bQ)\}$ such that 
\[N\models\mbox{`` }{\cal J}\mbox{ is countable, say $\{p_n:n <\omega\}$, and
${\cal J}$ is predense in $\bQ$ ''.}\]
Toward contradiction assume $p^*\in \bQ^N$ and $\name{r}\in N$ are such that
\[N\models\mbox{`` }p^*\forces_{\bP}\name{r}\in\bQ\mbox{ is incompatible with
every $p_n$ ''.}\]
We can replace $N$ by ${\cal H}(\chi_2)^N$. Let $M\in N$ be such that 
\[\begin{array}{ll}
N\models&\mbox{`` }M\prec ({\cal H}(\chi_2)^N,\in)\mbox{ is countable and}\\
\ &\bP,\bQ,{\frak B}^{\bP},{\frak B}^{\bQ},\langle p_n:n <\omega\rangle,p^*,
\name{r}\in M\mbox{ ''.}
\end{array}\]
In $N$ we can find $G\subseteq \bP^N$ generic over $M$, so $M$ inherits from
${\cal H}(\chi_2)$ the property $M[G]$ is a $\bQ$--candidate and also $M[G]
\models$``$\name{r}[G],p_n$ are incompatible in $\bQ$''. So $\name{r}[G]\in
\bQ^N$ contradicts the choice of ${\cal J}$. \QED

\begin{conclusion}
For \ref{6.5}$(\beta)$ to hold, we can omit clause (v) of (e) there if we add:
\begin{enumerate}
\item[(g)] $\bQ$ satisfies the c.c.c. in $\bQ$-candidates and being
incompatible in $\bQ$ is upward absolute from $\bQ$--candidates. 
\end{enumerate}
\end{conclusion}

We can conclude (phrased for simplicity for strongly c.c.c. nep).

\begin{conclusion}
\label{6.11}
Assume that
\begin{enumerate}
\item[(a)] $\bQ$ is strongly c.c.c. explicitly nep (see Definition
\ref{5.14}) and simple and correct, 
\item[(b)] $\name{\eta}\in \baire$ generic for $\bQ$, a hereditarily countable
$\bQ$--name. 
\end{enumerate}
If $\bP_0$ is nep, $I_{(\bQ,\name{\eta})}$--preserving and
$\forces_{\bP_0}$``$\name{\bP}_1$ is nep,
$I_{(\bQ,\name{\eta})}$--preserving'' 

then $\bP_0 * \name{\bP}_1$ is (nep and) $I_{(\bQ,\name{\eta})}$--preserving. 
\QED 
\end{conclusion}

The reader may ask: what about $\omega$ limits (etc)? We shall address these
problems in the continuation \cite{Sh:F264}.

\stepcounter{section}
\subsection*{\quad 8. Non-symmetry}
The following hypothesis \ref{7.1} will be assumed in this and the next
section, though for the end (including the main theorems
\ref{8.7}--\ref{8.10}) we assume snep (i.e.\ \ref{7.1A}). 

\begin{hypothesis}
\label{7.1}
$\bQ$ is correct c.c.c.\ simple, strongly c.c.c. nep, $\name{\eta}$ is a
hereditarily countable name of a generic real, i.e.\ $(\bQ,\name{\eta})\in\cK$
(see Definition \ref{5.2A}) for $\ZFC^-_*$ and $\ZFC^-_*$ (and the properties
above) are preserved by a forcing of cardinality $<\bar{\chi}$, $|\bQ|^{
\aleph_0}<\bar{\chi}$, for $\bQ$-candidates. 
\end{hypothesis}

\begin{hypothesis}
\label{7.1A}
Like \ref{7.1} with snep.
\end{hypothesis}
  
\begin{definition}
\label{7.2}
Let $\bQ,\name{\eta}$ be as in \ref{7.1} and let $\alpha$ be an ordinal.
\begin{enumerate}
\item Let $\bQ^{[\alpha]}$ be $\bP_\alpha$, where $\langle\bP_i,\nbQ_j:i\le
\alpha,j<\alpha\rangle$ is a FS iteration and $\nbQ_j= \bQ^{\V[\bP_j]}$.  
\item We let $\name{\eta}^{[\alpha]}$ be $\langle\name{\eta}_\ell:\ell<\alpha
\rangle$, where $\name{\eta}_\ell$ is $\name{\eta}$ ``copied to $\nbQ_\ell$''
(see \ref{7.3}(1) below).   
\item $(\bQ^{\langle\alpha\rangle},\name{\eta}^{\langle\alpha\rangle})$ is
defined similarly as an FS product. 
\item For a finite set $u\subseteq\alpha$ we define $F=F^{\alpha,u}_\bQ:\bQ
\longrightarrow \bQ^{[\alpha]}$ by $F(p)=\bar{p}$, where $\bar{p}=\langle
p_\ell:\ell<\alpha\rangle$, $p_\ell=p$ if $\ell\in u$ and $p_i=\emptyset_\bQ$
otherwise. 
\item The FS iteration $\bar{\bQ}=\langle\bP_i,\nbQ_j,\name{\eta}_j:i\le
\alpha,j<i\rangle$ of neps means $(\nbQ_j,\name{\eta}_j)\in\cK$.  

\noindent We write $\name{\eta}$ to mean $\name{\eta}_j=F^{\alpha,\{j\}}_{\bQ}
(\name{\eta})$.
\end{enumerate}
\end{definition}

\begin{proposition}
\label{7.3}
\begin{enumerate}
\item In Definition \ref{7.2}(4), for finite $u \subseteq\alpha$, $F=F^{\alpha,
u}_\bQ$ is a complete ($\lesdot$) embedding, as ``$p \le q$'', ``$p,q$
compatible'', ``$p,q$ incompatible'', ``$\langle p_n:n<\omega\rangle$ is
predense set above $q$'' are upward absolute from $\bQ$--candidates (holds as
$\bQ$ is strongly c.c.c.\ by \ref{7.1}). So $\name{\eta}_\ell$ is $F^{\alpha,
\{\alpha\}}_{\bQ}(\name{\eta})$ if $\alpha \in u$. 
\item $\bQ^{[\alpha]}$ satisfies the c.c.c. 
\item Same holds for $\bQ^{\langle \alpha \rangle}$. 
\item $(\bQ^{[\alpha]},\name{\eta}^{*[\alpha]})$ for $\alpha<\omega_1$ are as
in \ref{7.1}, too. 
\end{enumerate}
\end{proposition}

\Proof For example:

\noindent 3)\quad It is enough to prove it for finite $\alpha$, and this we
prove by induction on $\alpha$ for $\alpha=n+1$. For the c.c.c.\ use
``incompatibility is absolute'' for forcing by $\bQ^{\langle n\rangle}$, so we
can use the last phrase in \ref{7.1}.
\medskip

\noindent 4)\quad The main point here is the strong c.c.c., so let $N$ be a
$\bQ$-candidate (and $\alpha+1\subseteq N$) and 
\[N\models\mbox{`` }{\cal I}\subseteq\bQ^{[\alpha]}\mbox{ is predense ''}.\]
Let $G^{[\alpha]}\subseteq \bQ^{[\alpha]}$ be generic over $\V$ and for $\beta
\leq\alpha$, $G^{[\beta]}=G^{[\alpha]}\cap\bQ^{[\beta]}$. Show by induction on
$\beta$ that $G^{[\beta]}\cap (\bQ^{[\beta]})^N$ is a generic subset of
$(\bQ^{[\beta]})^N$ over $N\langle G^{[\beta]}\rangle$. \QED

\begin{definition}
\label{7.4}
\begin{enumerate}
\item We say that $\bQ$ is $[n]$--symmetric if:
\begin{quotation}
{\em if}\ $\langle\eta^*_\ell:\ell<n\rangle$ is generic for $\langle\bP_\ell,
\nbQ_\ell,\name{\eta}_\ell:\ell<n\rangle$ and $\sigma$ is a permutation of
$\{0,\ldots,n-1\}$\\ 
{\em then}\ $\langle\eta_{\sigma(\ell)}:\ell<n\rangle$ is generic for $\langle
\bP_\ell,\nbQ_\ell,\name{\eta}:\ell<n\rangle$.
\end{quotation}
\item  If $(\bQ',\name{\eta}')$, $(\bQ'',\name{\eta}'')$ are as in \ref{7.1},
we say that they commute if: 
\begin{quotation}
{\em if}\ $r'$ is $(\bQ',\name{\eta}')$--generic over $\V$ and $r''$ is
$(\bQ'',\name{\eta}'')$--generic over $\V[r']$\\ 
{\em then}\ $r'$ is $(\bQ',\name{\eta}')$--generic over $\V[r'']$
\end{quotation}
(note that $\eta''$ is $(\bQ'',\name{\eta}'')$--generic over $\V$ is always
true by \ref{5.3B}).
\item For $(\bQ',\name{\eta}')$, $(\bQ'',\name{\eta}'')$ we say that they
weakly commute if $(\bQ'\rest(\ge q'),\name{\eta}')$, $(\bQ''\rest(\ge q''),
\name{\eta}'')$ commute for some $q'\in\bQ'$ and $q''\in\bQ''$.
\end{enumerate}
\end{definition}

\begin{proposition}
\label{7.5}
\begin{enumerate}
\item ``Commute'' is a commutative relation.
\item For $n \ge 2$ we have:\\
$\bQ$ is $[n]$--symmetric \quad iff\\
$\bQ,\bQ^{[n-1]}$ commute and $\bQ$ is $[n-1]$--symmetric \quad iff\\
$\bQ$ is $[2]$-symmetric. 
\item If $\bP,\bQ^{[n]}$ commute, $m \le n$ then $\bP,\bQ^{[m]}$
commute. Similarly, if $\bP,\bQ$ commute and $\bQ'\lesdot\bQ$, the $\bP,\bQ'$
commute.
\item In part 3) we can replace $[-]$ by $\langle - \rangle$.
\item If $(\bQ',\name{\eta}')$, $(\bQ'',\name{\eta}'')$ weakly commute and
$\bQ',\bQ''$ are homogeneous, then they commute.
\end{enumerate}
\end{proposition}

\Proof  1)\quad Let $(\bQ',\name{\eta}')$, $(\bQ'',\name{\eta}'')$ be as in
\ref{7.1}. Then ``$(\bQ',\name{\eta}')$, $(\bQ'',\name{\eta}'')$ commute''
says $\bQ'*\name{\bQ}''=\bQ''*\name{\bQ}'$, which is symmetric.
\medskip

\noindent 2)\quad For the second ``iff'', use ``the permutations $\pi_\ell=
(\ell,\ell +1)$ for $\ell < n$ generate the group of permutations of
$\{0,\ldots,n-1\}$''.  \QED$_{\ref{7.5}}$ 

\begin{proposition}
\label{7.6}
\begin{enumerate}
\item If $\bQ^{[\omega]}$ and $\Cohen$ do not commute, {\em then} for some $n<
\omega$, $\bQ^{[n]}$ and $\Cohen$ do not commute.\\
(The inverse holds by \ref{7.5}(3), second phrase.) 
\item If $\bQ^{\langle\omega\rangle}$ and $\Cohen$ do not commute, {\em then}
for some $n<\omega$, $\bQ^{\langle n\rangle}$ and $\Cohen$ do not commute.
\end{enumerate}
\end{proposition}

\Proof  1)\quad Since $\Cohen$ and $\bQ^{[\omega]}$ do not commute, there is a
$\bQ^{[\omega]}$--name $\name{{\cal I}}$ of a dense open subset of $\Cohen$
(i.e.\ of $({}^{\omega >}2,\vartriangleleft)$) such that for some condition
$(p,\name{q})\in\Cohen*\name{\bQ}^{[\omega]}$ we have
\[(p,\name{q})\forces\mbox{`` }\name{\eta}^{\Cohen}\mbox{ has no initial
segment in }\name{{\cal I}}\mbox{ ''}.\]
Without loss of generality for some $n^*<\omega$ we have
$p\forces_{\Cohen}$``$\dom(\name{q})\subseteq\{0,\ldots,n^*-1\}$''. 
Let $\name{{\cal I}}'$ be the $\bQ^{[n^*]}$--name for the following set:
\[\{\eta\in {}^{\omega >}2:\mbox{for some } p \in\bQ^{[\omega]},\ p\rest n^*
\in \name{G}_{\bQ^{[n^*]}}\mbox{ and }p\forces_{\bQ^{[\omega]}}\mbox{``}\eta
\in \name{\cal I}\mbox{''}\}.\]
It should be clear that $\forces_{\bQ^{[n^*]}}$``$\name{\cal I}'$ is a dense
open subset of $({}^{\omega>}2,\vartriangleleft)$''. Now we ask the following
question.
\begin{quotation}
Does $(p,\name{q})\forces_{\Cohen*\nbQ^{[n^*]}}$`` $\name{\eta}^{\Cohen}\rest
n \notin \name{{\cal I}}$ for each $n<\omega$ ''?  
\end{quotation}
If yes, we have gotten the desired conclusion (i.e.\ $\Cohen$ and
$\nbQ^{[n^*]}$ do not commute). If not, for some $(p',\name{q}')$ such that 
$(p,\name{q})\le (p',\name{q}')\in\Cohen*\nbQ^{[n^*]}$ and for some
$n<\omega$ we have:
\[(p',\name{q}')\forces_{\Cohen*\nbQ^{[n^*]}}\mbox{`` }\name{\eta}^\Cohen\rest
n=\eta\in\name{{\cal I}}'\mbox{ ''}.\]
Without loss of generality, for some $p\in(\bQ^{[\omega]})^\V$ we have $(p',
q')\forces$``$p\rest n^*\in G_{\Cohen*\nbQ^{[n^*]}}$'' and $p\forces$``$\eta
\in \name{{\cal I}}$''. Then $(p',q'\cup p\rest [n^*,\omega))$ forces (in
$\Cohen*\nbQ^{[\omega]}$) that $\name{\eta}^\Cohen\rest n=\eta\in \name{{\cal
I}}$, a contradiction. 
\medskip

\noindent 2)\quad Similarly. \QED$_{\ref{7.6}}$

\begin{proposition}
\label{7.7}
\begin{enumerate}
\item If $\bQ^{[n]}$ and $\Cohen$ do not commute ($\bQ$ as before), then $\bQ$
and $\Cohen$ do not commute (both ``absolute").  
\item The following conditions are equivalent:
\begin{enumerate}
\item[(i)]    $\bQ$ commutes with Cohen,
\item[(ii)]   $\forces_{\bQ}$``$(\can)^\V$ is not meagre'',
\item[(iii)]  $(\forall A)[\V\models$``$A\subseteq\can$ non-meagre'' $\ \
\Rightarrow\ \ \forces_{\bQ}$``$A$ is non-meagre''$]$ 
\end{enumerate}
(all ``absolutely'', i.e.\ not only in the present universe but in its generic
extensions too). 
\item We can replace $\Cohen$ by others to which \ref{6.5} applies and are
homogeneous (see \ref{6.4A}). 
\end{enumerate}
\end{proposition}

\Proof  1)\quad Assume toward contradiction that $Q$ and $\Cohen$ commute
(absolutely). Let $\eta\in\can$ be a Cohen real over $\V$. Let $G_\ell
\subseteq\bQ^{\V[G_0,\ldots,G_{\ell-1},\eta]}$ be generic over $\V[G_0,\ldots,
G_{\ell-1},\eta]$ for $\ell<n$, and let $\eta_\ell=\name{\eta}[G_\ell]$. We
now prove by induction on $\ell$, that $\eta$ is a Cohen real over $\V[G_0,
\ldots,G_{\ell -1}]$. The induction step is by the assumption ``$\bQ$ and
$\Cohen$ commute''. The net result is that $\eta$ is a Cohen real over
$\V[\eta_0,\ldots,\eta_{n-1}]$, contradicting the assumption. 
\medskip

\noindent 2)\quad The second clause implies the third by \ref{6.5}.  The third
clause implies the second trivially.  

Let us argue that the implication $(i)\ \Rightarrow\ (ii)$ holds. Add
$\aleph_1$ Cohen reals $\{\eta_i:i<\omega_1\}$ and then force by $\bQ$. Let
$G_\bQ\subseteq \bQ^{\V[\langle \eta_i:i<\omega_1\rangle]}$ be generic over
$\V$, and $\name{\eta}=\name{\eta}_{\bQ}[G_{\bQ}]$. Then (i) implies that for
every $j<\omega_1$ we have: $\eta_j$ is Cohen over $\V[\langle\eta_i:i<
\omega_1,i\ne j\rangle,\eta]$. Hence in $\V[\langle\eta_i:i<\omega_1\rangle,
\eta]= \V[\langle\eta_i:i<\omega_1\rangle][G_{\bQ}]$, the set $\{\eta_i:i<
\omega_1\}$ is not meagre and consequently (ii) holds.  

Lastly, assume (iii) and let $\nu\in {}^\omega|\bQ|$ be generic for
$\Levy(\aleph_0,|\bQ|)$. Let $\eta$ be $(\bQ,\name{\eta})$--generic real over
$\V[\nu]$. By (ii), we can find in $\V[\nu]$ a real $\rho\in \can$ which is in
no meagre set from $\V[\eta]$ (note that there are countably many such meagre
sets from the point of view of $\V[\nu]$). Now we easily finish. 
\medskip

\noindent 3)\quad Same proof. \QED$_{\ref{7.7}}$

\stepcounter{section}
\subsection*{\quad 9. Poor Cohen commutes only with himself}

\begin{definition}
\label{8.1}
\begin{enumerate}
\item We say a $\bQ$--name $\name{x}$ of a subset of some countable $a^*\in
\V$ is [somewhere] essentially Cohen if $B_2(\bQ,\name{x})$ is [somewhere]
essentially countable; i.e.\ [above some $p$] has countable density.
\item We say $(\bQ,\name{\eta})\in \cK^{\neg c}$ (a non-Cohen pair) if:
\begin{enumerate}
\item[(a)] $(\bQ,\name{\eta})$ is as in \ref{7.1A},
\item[(b)] $(\bQ,\name{\eta})$ (see Definition \ref{5.3A}) is nowhere
essentially Cohen (i.e.\ above every condition). 
\end{enumerate}
\end{enumerate}
\end{definition}

\begin{hypothesis}
\label{8.2}
$\chi$ is regular large enough cardinal, and $(\bQ,\name{\eta})\in\cK^{\neg c}$
will be fixed as in \ref{8.1}, and $\ZFC^-_*$ is normal (see Definition
\ref{0.9}).  
\end{hypothesis}

\begin{definition}
\label{8.2A}
\begin{enumerate}
\item ${\cal D}={\cal D}_{\le \aleph_0}({\cal H}(\chi))$ is the filter of
clubs on $[{\cal H}(\chi)]^{\le \aleph_0}$. 
\item ${\cal C}_0=\{a:a\prec ({\cal H}(\chi),\in)$ is countable, and
$(\bQ,\name{\eta})\in a$ (i.e.\ their definitions) so is a
$\bQ$--candidate$\}$.
\end{enumerate}
\end{definition}

\begin{definition}
\label{8.2B}
We say that $q\in\bQ$ is strong on $a\in {\cal C}_0$ if:
\begin{enumerate}
\item[$(\circledast)_{a,q}$] the set $\{p\in a\cap\bQ$: $p,q$ are incompatible
in $\bQ\}$ is dense in the (quasi) order $\bQ\cap a$.
\end{enumerate}
\end{definition}

\begin{proposition}
\label{8.3}
\begin{enumerate}
\item For every $a\in {\cal C}_0$ there is $q\in\bQ$ which is strong on $a$.
\item Moreover, for every $p\in\bQ$ and $a\in {\cal C}_0$ there is $q$ strong
on $a$ such that $p\le^{\bQ} q$. 
\end{enumerate}
\end{proposition}

\Proof Clearly $\name{G}_{\bQ}\cap a$ is a $\bQ$--name of a countable subset
of an old set $\bQ\cap a$, so it can be considered as a real. Note that
\begin{enumerate}
\item[$(*)_1$] $\name{G}_{\bQ}\cap a$ is not somewhere essentially Cohen.
\end{enumerate}
Why? We can restrict ourselves to be above some fix $p\in\bQ$. From
$\name{G}_{\bQ}\cap a$ we can compute $\name{\eta}$ (as $\name{\eta}\in a$,
i.e.\ the relevant maximal antichains belong to $a$), so $\name{\eta}$ can be
considered a $B_2[\bQ,\name{G}_{\bQ}\cap a]$--name. But ``any (name of) a real
in an essentially Cohen forcing notion is essentially Cohen itself'', so
$\name{\eta}$ is essentially Cohen $\bQ$--name, contradicting Hypothesis
\ref{8.2}. 

Consequently, $\forces_{\bQ}$``$\name{G}_{\bQ}\cap a$ is not a generic subset
of $\bQ\rest a$ (over $\V$)'' and  hence $p\forces_Q$``$\name{G}_{\bQ}\cap a$
is not a generic subset of $\bQ\rest a$ (over $\V$)''. Thus there are $q$ and
${\cal I}$ such that: 
\begin{enumerate}
\item[(i)]   $p \le q \in \bQ$,
\item[(ii)]  ${\cal I}\subseteq\bQ\cap a$ is a dense open subset of $Q\rest
a$, 
\item[(iii)] $q\forces_{\bQ}$``$\name{G}_{\bQ}$ is disjoint to ${\cal I}$''.
\end{enumerate}
But this means that 
\begin{enumerate}
\item[$(*)_2$] $q$ is incompatible with every $r \in {\cal I}$.
\end{enumerate}
[Why? Otherwise $q\not\forces_{\bQ}$``$r\notin\name{G}_{\bQ}$''.]\\
So $\{r\in a\cap\bQ:q,r$ incompatible (in $\bQ$)$\}$ is a subset of $\bQ\cap
a$ including ${\cal I}$ hence it is dense in $\bQ\rest a$. \QED$_{\ref{8.3}}$

\begin{choice}
\label{8.3A}
We choose $\bar{p}=\langle p_a:a\in {\cal C}_0\rangle$ such that $p_a\in\bQ$ is
strong on $a$ (possible by \ref{8.3}). 
\end{choice}

\begin{definition}
\label{8.4}
\begin{enumerate}
\item For $R \subseteq \bQ$ let $A[R]\stackrel{\rm def}{=}\{a\in {\cal C}_0:
p_a\in R\}$.
\item ${\cal D}_{\bar{p}}={\cal D}_{\bQ,\bar{p}}\stackrel{\rm def}{=}\{R
\subseteq\bQ:A[R]\in{\cal D}\}$.\\
The family of ${\cal D}_{\bar{p}}$--positive sets will be denoted ${\cal
D}^+_{\bar{p}}$ (so for a set $S\subseteq \bQ$, $S\in {\cal D}^+_{\bar{p}}$
iff $R\cap S\neq\emptyset$ for each $R\in {\cal D}_{\bar{p}}$).
\item For $R\subseteq\bQ$ and $q\in\bQ$ let $R[q]\stackrel{\rm def}{=}\{p\in
R:p,q$ are incompatible in $Q\}$ (so $R[q]$ is in a sense the orthogonal 
complement of $q$ inside $R$).
\end{enumerate}
\end{definition}

\begin{fact}
\label{8.4A}
\begin{enumerate}
\item ${\cal D}_{\bar{p}}$ is an $\aleph_1$--complete filter on $\bQ$.
\item  For $R\subseteq\bQ$ we have $R\in {\cal D}^+_{\bar{p}}\ \Leftrightarrow 
\ A[R]\in {\cal D}^+$.
\end{enumerate}
\end{fact}

\begin{proposition}
\label{8.5}
If $R\in {\cal D}^+_{\bar{p}}$ then the set
\[R^\otimes\stackrel{\rm def}{=}\{q\in\bQ:R[q]\in {\cal D}^+_{\bar{p}}\}\]
is dense in $\bQ$.
\end{proposition}

\Proof Assume not, so for some $q^* \in Q$ we have
\begin{enumerate}
\item[$(*)_1$] there is no $q\in\bQ$ such that $q^*\le q\in\bQ\ \&\ R[q]\in
{\cal D}^+_{\bar{p}}$. 
\end{enumerate}
Thus
\[\begin{array}{lll}
q^*\le q\in\bQ&\Rightarrow &R[q]=\emptyset\mbox{ mod }{\cal D}_{\bar{p}}\ \
 \Rightarrow\ \ A[R[q]]=\emptyset\mbox{ mod } {\cal D}\\
\ &\Rightarrow&\mbox{for some club }{\cal C}_q \subseteq {\cal C}_0\mbox{ of }
 [{\cal H}(\chi)]^{\le \aleph_0} \mbox{ we have} \\ 
\ &\ &(\forall a \in {\cal C}_q)[p_a \notin R[q], \mbox{ i.e.\ $p_a,q$ are
compatible}]. 
  \end{array}\]
Let ${\cal C}^*=\{a\in {\cal C}_0:q^*\in a\mbox{ and }(\forall q)[q^*\le q\in
a \cap\bQ\ \Rightarrow\ a\in {\cal C}_q]\}$. As each ${\cal C}_q$ is a club of
$[{\cal H}(\chi)]^{\le \aleph_0}$ clearly ${\cal C}^*$ (as a diagonal
intersection) is a club of $[{\cal H}(\chi)]^{\le \aleph_0}$, i.e.\ ${\cal
C}^* \in {\cal D}$. Since $R\in {\cal D}^+_{\bar{p}}$ we have $A[R]\in {\cal
D}^+$, so together with the previous sentence we know that there is $a^*\in
A[R]\cap{\cal C}^*$. By the choice of $\bar{p}$ (see \ref{8.3A}, and
Definition \ref{8.2B}) as $q^*\in a^*\cap \bQ$ (see the choice of ${\cal
C}^*$) for some $q$ we have: 
\[q^* \le q\in a^*\quad\mbox{ and }\quad p_{a^*},q\mbox{ are incompatible}.\]
Now this contradicts ``$a^*\in C_q$''. \QED$_{\ref{8.6}}$

\begin{definition}
\label{8.6}
Assume $\chi_1=(2^\chi)^+$ (so ${\cal H}(\chi)\in {\cal H}(\chi_1))$ and $N$
is a countable elementary submodel of $({\cal H}(\chi_1),\in)$ to which
$\{\chi,\bQ,\bar{p}\}$ belong (so ${\cal D}_{\bar{p}} \in N$). Further, assume
that $\bQ$ is snep.
\begin{enumerate}
\item We let $\Cohen_N=\Cohen_{N,\bQ}$ be $({\cal D}^+_{\bQ,\bar{p}},
\supseteq)\rest N$ (so this is a countable atomless forcing notion and hence
equivalent to Cohen forcing).
\item If $G_N\in\mbox{ Gen}(N,{\cal D}^+_{\bQ,\bar{p}})\stackrel{\rm def}{=} 
\{G:G\subseteq\Cohen_N\mbox{ is generic for }(N,({\cal D}^+_{\bQ,\bar{p}},
\supseteq)\rest N)\}$ (possibly in a universe $\V'$ extending $\V$) {\em then}
let $\name{p}_N[G]$ be the sequence (i.e.\ in $\baire$ or just member of
${}^\omega\theta(\bQ)$) such that for each $\ell<\omega$ and $\gamma$ 
\[(\name{p}_N[G_N])(\ell)=\gamma\quad \Leftrightarrow\quad (\exists R\in G)(
\forall p\in R)[p(\ell)=\gamma].\]
\end{enumerate}
\end{definition}

\begin{proposition}
\label{8.7}
Assume \ref{7.1A} and, additionally, $\bQ$ is Souslin c.c.c.\ (i.e.\ the
incompatibility relation is $\Sigma^1_1$). If $\chi_1,N$ and $G\in\mbox{ Gen}
(N,{\cal D}^+_{\bQ})$ are as in \ref{8.6} (so $G$ is possibly in some generic
extension $\V_1$ of $\V$ but $\Cohen_N$ is from $\V$) {\em then}
\begin{enumerate}
\item[(a)]  $\name{p}_N[G]$ is an $\omega$--sequence (i.e.\ for each $\ell$
there is one and only one $\gamma$), 
\item[(b)]  $\name{p}_N[G]\in\bQ$,
\item[(c)]  $\name{p}_N[G]$ is strong for $N\rest {\cal H}(\chi)$ (which
belongs to ${\cal C}_0$). 
\end{enumerate}
\end{proposition}

\Proof For every $p\in\bQ$ there is $\nu_p\in\baire$ which witnesses $p\in\bQ$,
i.e.\ $p*\nu_p\in\lim(T^{\bQ}_0)$. So choose such a function $p\mapsto\nu_p$.
Now in $\V$, for $n<\omega$ the function $p_a\mapsto (p_a\rest n,\nu_{p_a}
\rest n)$ is a mapping from $\{p_a:a \in {\cal C}_0\}\in {\cal D}_{\bar{p}}$
with countable range. Since ${\cal D}_{\bar{p}}$ is $\aleph_1$--complete
\begin{enumerate}
\item[$(*)_1$] in $\V$, if $R\in {\cal D}^+_{\bar{p}}$ and $n<\omega$ then for
some $R'\subseteq R$ and $(\eta^n,\nu^n)$ we have 
\[R'\in {\cal D}^+_{\bar{p}}\quad\mbox{ and }\quad(\forall p\in R')[(p\rest n,
\nu_p \rest n)=(\eta^n,\nu^n)].\]
\end{enumerate}
This is inherited by $N$, hence $\name{p}_N[G]$ satisfies clauses (a), (b) (in
fact 
\[\name{\nu}[G]=\bigcup\{\nu^*:\mbox{ for some } n<\omega\mbox{ and } R\in G
\mbox{ we have } (\forall p\in R)[\nu_p\rest n=\nu^*]\}\]
is a witness for $\name{p}_N[G]\in\bQ$). Also for each $q\in\bQ\cap N$ the set
\[\begin{array}{ll}
{\cal J}_q=\bigl\{R\in{\cal D}^+_{\bar{p}}:&\mbox{for some }q'\in\bQ \mbox{
stronger than $q$ we have:}\\
\ &(\forall p \in R)[p,q' \mbox{ are incompatible (in $\bQ$)}]\bigr\}
  \end{array}\]
is a dense subset of $({\cal D}^+_{\bar{p}},\supseteq)$ (remember $p_a$ is
strong on $a$; use Fodor lemma). Clearly it belongs to $N$, so by the demand
on $G$ we know that $G\cap {\cal J}_q \ne\emptyset$. Choose $R_q\in G\cap
{\cal J}_q$ and let $q'\in \bQ\cap N$ witness it, so 
\[R_q\in {\cal D}^+_{\bar{p}}\cap N\quad\mbox{ and }\quad(\forall p\in R_q)
[p,q'\mbox{ are incompatible}].\]
Now ``incompatible in $\bQ$'' is a $\Sigma^1_1$--relation (belonging to $N$)
hence as above, $\name{p}_N[G],q'$ are incompatible. As $q$ was any member of
$\bQ\cap N$ we have finished proving clause (c). \QED$_{\ref{8.7}}$

\begin{proposition}
\label{8.8}
Assume \ref{7.1A} and let $\bQ$ be Souslin c.c.c. Then $\bQ^{[\omega]}$ (see
\ref{7.2}) and $\Cohen$ do not commute.
\end{proposition}

\Proof Assume that $\bQ^{[\omega]}$ and $\Cohen$ do commute. Let $\chi$ be
large enough, $N \prec ({\cal H}(\chi),\in)$ be countable such that $(\bQ,
\name{\eta})\in N$ (as in \ref{8.6}). Now we can interpret a Cohen real $\nu$
(over $\V$) as a subset of ${\cal D}^+_{\bar{p}}\cap N$ called $g_\nu$. Thus
it is $\Cohen_{N,\bQ}$--generic over $\V$ so $\name{p}_{N}[g_\nu]$ is well
defined, and it belongs to $\bQ^{\V[\nu]}$ (by \ref{8.7}). Moreover, in
$\V[\nu]$ we have: 
\[\{q\in\bQ^N:q,\name{p}_N(g_\eta)\mbox{ are incompatible $\}$ is dense in }
\bQ^N.\]
Let $\langle \eta_\ell:\ell<\omega\rangle$ be generic for $(\bQ^{[\omega]},
\name{\eta}^{[\omega]})$ and let $\nu$ be Cohen generic over $\V[\langle
\eta_\ell:\ell<\omega\rangle]$. For each $\ell$, clearly $\eta_\ell$ is $(\bQ,
\name{\eta})$--generic over $\V$, so let $\eta_\ell=\name{\eta}[G_\ell]$,
where $G_\ell\subseteq\bQ$ is generic over $\V$. Clearly $G_\ell\cap N$ is a
subset of $\bQ^N$ generic over $\V$ (by ``$\bQ$ is strongly c.c.c.''). So
$\langle G_\ell\cap N,g_\nu\rangle$ is a subset of $\bQ^N*({\cal D}^+_{
\bar{p}}\cap N,\supseteq)$ generic over $N$. By \ref{8.7}, for any $q\in\bQ^N$
and $R\in({\cal D}^+_{\bar{p}}\cap N)$, for some $R'\subseteq R$ and $q'$ we
have $R'\in ({\cal D}^+_{\bar{p}}\cap N)$, $N\models$``$q\leq q'\in\bQ$'' and
\[N\models(\forall a\in R')(p_a,q'\mbox{ are incompatible\/}).\]
So look at the set 
\[\{(q,R)\in\bQ^N\times({\cal D}^+_{\bar{p}}\cap N): (\forall a\in R')(p_a,
q\mbox{ are incompatible\/})\}\]
-- there is $(q,R)\in (G\cap N)\times g_\nu$ which belongs to it. Hence, as in
\ref{8.7}, for each $\ell$, $\name{p}_N[g_\nu]$ is incompatible with some
$q\in G_{\bQ}[\eta_\ell]$.  

By the assumption that the forcing notions commute we know that $\langle
\eta_\ell:\ell<\omega\rangle$ is generic for $(\bQ^{[\omega]},\name{\eta}^{
[\omega]})$ over $\V(\nu)$. Necessarily (by FS + genericity) for some $\ell$
we have $F^{\omega,\{\ell\}}_{\bQ}(\name{p}_N(g_\eta))\in G_{\bQ}[\langle
\eta_\ell:\ell<\omega\rangle]$; a contradiction. \QED$_{\ref{8.8}}$

\begin{conclusion}
\label{8.9}
Assume \ref{7.1A} and let $\bQ$ be Souslin c.c.c. Then $(\bQ,\name{\eta})$
does not commute with $\Cohen$ (even above any $q\in\bQ$).
\end{conclusion}

\Proof If we restrict ourselves above $q_0\in\bQ$, the Hypothesis \ref{8.2}
still holds so we can ignore this. By \ref{8.8} we have $(\bQ^{[\omega]},
\name{\eta}^{[\omega]})$ does not commute with $\Cohen$. So by \ref{7.6} we
have that, for some $n$, $(\bQ^{[n]},\name{\eta}^{[n]})$ does not commute with
$\Cohen$ and by \ref{7.7} we finish. \QED$_{\ref{8.9}}$

\begin{proposition}
\label{8.9A}
If $\bQ$ is Souslin c.c.c.\ then for suitable $\ZFC^-_*$, $\bQ$ satisfies
\ref{7.1A}. 
\end{proposition}

\Proof Let $\rho\in\can$ be the real parameter in the definition of $\bQ$. Let
$\ZFC^-_*$ say:
\begin{enumerate}
\item[(a)] ZC (i.e.\ the axioms of Zermelo satisfied by $({\cal H}(
\beth_\omega),\in)$),
\item[(b)] $\bQ$ (defined from $\rho$ which is an individual constant)
satisfies the c.c.c. 
\item[(c)] for each $n<\omega$, generic extensions for forcing notions of
cardinality $\le\beth_\omega$ preserve (b) (and, of course (a)). 
\end{enumerate}
Now the desired properties are easy. \QED$_{\ref{8.9A}}$

\begin{conclusion}
\label{8.10}
If $\bQ$ is a Souslin c.c.c.\ forcing notion which is not $\baire$--bounding
(say $p\forces$`` there is an unbounded $\name{\eta}\in\baire$ ''), but adds
an essentially non-Cohen real {\em then} $\bQ$ does not commute with itself.
\end{conclusion}

\Proof By \cite{Sh:480}, $\bQ$ adds a Cohen real; now by the assumptions, for
some $\bQ$--name $\name{\eta}$, $(\bQ,\name{\eta})\in\cK^{\neg c}$. By
\ref{8.9} we know that $\bQ$ and Cohen do not commute, so by \ref{7.5}(3) we
are done. \QED$_{\ref{8.10}}$ 

\begin{conclusion}
\label{8.11}
If $\bQ$ is a Souslin c.c.c.\ forcing notion adding a non-Cohen real, {\em
then} the forcing by $\bQ$ makes the old reals meagre. 
\end{conclusion}

\stepcounter{section}
\subsection*{\quad 10. Some c.c.c.\ nep forcing notions are not nice}
We may wonder can we replace the assumption ``$\bQ$ is Souslin c.c.c.'' by
weaker one in \S8 and in \cite{Sh:480}. We review limitations and then see how
much we can weaken it. 

\begin{proposition}
\label{9.1}
Assume that $\eta^* \in\can$  and $\aleph_1=\aleph^{{\bf L}[\eta^*]}_1$. Then
there is a definition of a forcing notion $\bQ$ (i.e.\ $\bar{\varphi}$) such
that 
\begin{enumerate}
\item[(a)] the definition is $\Sigma^1_1$ (with parameter $\eta^*$), so $p
\in\bQ$, $p \le^{\bQ} q$, ``$p,q$ incompatible'', ``$\{p_n:n<\omega\}\subseteq
a$ is a maximal antichain of $\bQ$'' are preserved by forcing extensions, 
\item[(b)] $\bQ$ is c.c.c.\ (even in a forcing extension; even
$\sigma$--centered), 
\item[(c)] there is $\bQ$--name $\name{\eta}$ of a generic for $\bQ$,
\item[(d)] $\name{\eta}$ is not essentially Cohen (preserved by extensions not
collapsing $\aleph_1$), in fact has cardinality $\aleph_1$, 
\item[(e)] $\bQ$ commutes with $\Cohen$,
\item[(f)] $\bQ$ is nep (though not Souslin c.c.c.).
\end{enumerate}
\end{proposition}

\Proof A condition $p$ in $\bQ$ is a quadruple $\langle E_p,X_p,u_p,w_p
\rangle$ consisting of:\quad a 2-place relation $E_p$ on $\omega$ and subset
$X_p$ of $\omega$ and a finite subset $u_p$ of $X_p$ and a finite subset $w_p$
of $\omega$ such that:
\begin{quotation}
$N_p\stackrel{\rm def}{=}(\omega,E_p)$ is a model of $\ZFC^- + \V={\bf L}$ 
(let $<^{N_p}_*$ be the canonical ordering of $N_p$, we do not require well
foundedness) such that: 
\end{quotation}

\[\begin{array}{l}
(N_p,X_p)\models\mbox{`` }(\alpha)\ \mbox{every } x\in X_p\mbox{ is an
  infinite subset of }\omega,\\ 
\quad(\beta)\ \mbox{if } x\ne y\mbox{ are from } X_p \mbox{ then }x \cap y
  \mbox{ is finite},\\ 
\quad(\gamma)\ \mbox{if } x\in X\mbox{ then there is no } y\mbox{ satisfying}
\\ 
\ y<^{N_p}_* x\ \&\ (\forall z \in X_p)(z <^{N_p}_* x\Rightarrow z\cap y
\mbox{ finite})\ \&\ y \mbox{ an infinite subset of }\omega,\\
\quad(\delta)\ \bigwedge\limits_{n<\omega} (\forall z_1\ldots z_n\in X_p)\bigl(
\bigwedge\limits_{\ell =1}^n z_\ell <^{N_p}_* x \Rightarrow (\exists^\infty m
< \omega)(m \notin x \cup\bigcup\limits_{\ell =1}^n z_\ell \bigr)\mbox{ ''.}
  \end{array}\]
The order is defined by:\qquad $p \le q$ if and only if one of the following
occurs: 
\begin{enumerate}
\item[(A)] $p=q$,
\item[(B)] there are $Y\subseteq\omega$ and $a\in N_q$ and $f\in {}^{
\textstyle Y}\omega$ such
that: 
\begin{enumerate}
\item[(i)]    $[x\in Y\ \&\  N_p \models y\in x]\quad \Rightarrow\quad y \in
Y$, 
\item[(ii)]   $[N_p\models$``$\rk(x)=y$'', $y\in Y]\quad\Rightarrow\quad x\in
Y$,  
\item[(iii)]  $N_p\rest Y$ is a model of $(\ZFC^-+\V={\bf L})$,
\item[(iv)]   the set $\{x:N_p \models$``$x$ an ordinal'', $x \notin Y\}$ has
no first element, 
\item[(v)]    $N_q\models$``$a$ is a transitive set'', 
\item[(vi)]   $f$ is an isomorphism from $N_p\rest Y$ onto $N_q\rest\{b:N_q 
\models b\in a\}$,
\item[(vii)]  $f$ maps $X_p$ onto $X_q\rest\rng(f)$,
\item[(viii)] $f$ maps $u_p \cap Y$ into $u_q \cap\rng(f)$,
\item[(ix)]   $w_p\subseteq w_q$,
\item[(x)]    if $n\in w_q\backslash w_p$ and $x\in f(u_p)$ then $N_q
\models$``the $n$-th natural number does not belong to $x$''.
\end{enumerate}
\end{enumerate}
The reader can now check (note that $\name{w}=\bigcup\{w^p:p\in\name{G}_{\bQ}
\}$ is forced to be an infinite subset of $\omega$ almost disjoint to every $A
\in X^*$, $X^*$ a reasonably defined MAD family in ${\bf L}$); see more
details in the proof of \ref{9.3}. $\QED_{\ref{9.1}}$ 

\begin{proposition}
\label{9.2}
Assume $\V={\bf L}$. There is $\bQ=\bQ_0 * \nbQ_1$ such that:
\begin{enumerate}
\item[(a)] $\bQ_0$ is nep c.c.c.\ not adding a dominating real,
\item[(b)] $\forces_{\bQ_0}$``$\nbQ_1$ is nep c.c.c. (even Souslin c.c.c.) not
adding a dominating real'',
\item[(c)] $\bQ$ adds a dominating real,
\item[(d)] in fact, $bQ_0$ is the Cohen forcing (so in any $\V_1$ it is
c.c.c.\ strongly c.c.c., correct, very simple nep (and snep), and it is really
absolute, i.e.\ it is the same in $\V_1$ and $\V$, and its definition uses no
parameters), 
\item[(e)] moreover, $\bQ_1$ is defined in ${\bf L}$, really absolute, and in
any $\V_1$ it is c.c.c., strongly c.c.c.\ nep (and even snep). In $\V_1$,
$\bQ_1$ adds a dominating real iff $(\baire)^{\bf L}$ is a dominating family
in $\V_1$.
\end{enumerate}
\end{proposition}

\Proof Let $\bQ_0$ be $\Cohen$. We shall define $\bQ_1$ in a similar manner as
$\bQ$ in the proof of \ref{9.1}. 

A condition in $\nbQ_1$ is a triple $\langle E_p,u_p,w_p\rangle$ such that
$E_p$ is a 2-place relation on $\omega$, $u_p$ is a finite subset of $\omega$
and $w_p$ is a finite function from a subset of $\omega$ to $\omega$ and:
\begin{quotation}
$N_p\stackrel{\rm def}{=}(\omega,E_p)$ is a model of $\ZFC^- + \V={\bf L}$
(let $<^{N_p}_*$ be the canonical ordering of $N_p$, we do not require well
foundedness); so in formulas we use $\in$. 
\end{quotation}
[What is the intended meaning of a condition $p$? Let 
\[M_p=N_p\rest\{x: ({\rm Tc}(x)^{N_p}, E_p\rest {\rm Tc}(x)^{N_p})\mbox{ is
well founded\/}\},\]
where ${\rm Tc}(x)$ is the transitive closure of $x$. Let $M_p'$ be the
Mostowski collapse of $M_p$, $h_p:M_p\longrightarrow M_p'$ be the isomorphism. 
Now, $p$ gives us information on the function $\name{w}=\bigcup\{w_p:p\in
\name{G}\}$ from $\omega$ to $\omega$, it says: $\name{w}$ extends the function
$w_p$ and if $x\in M_p\cap u_p$ is a function from $\omega$ to $\omega$ then
for every natural number $n\notin\dom(w_p)$ we have $x(n)\leq\name{w}(n)$. 
Note that $h_p(x)$ is a function from $\omega$ to $\omega$ iff $M_p\models$``
$x$ is a function from $\omega$ to $\omega$ '' iff $N_p\models$`` $x$ is a
function from $\omega$ to $\omega$ ''.]\\
The order is defined by:\qquad $p \le q$ if and only if one of the following
occurs: 
\begin{enumerate}
\item[(A)] $p=q$,
\item[(B)] there are $Y \subseteq \omega$ and $a\in N_q$ and $f\in {}^{
\textstyle Y}\omega$ such that 
\begin{enumerate}
\item[(i)]     $[x\in Y \ \&\ N_p\models y\in x]\quad \Rightarrow\quad y\in
Y$,  
\item[(ii)]    $[N_p\models$``$\rk(x)= y$'' $\ \&\ y\in Y]\quad\Rightarrow
\quad x\in Y$,
\item[(iii)]   $N_p\rest Y$ is a model of $(\ZFC^- + \V={\bf L})$,
\item[(iv)]    the set $\{x:N_p \models$ ``$x$ an ordinal'', $x \notin Y\}$
has no first element (by $E_p$),
\item[(v)]     $N_q\models$``$a$ is a transitive set'',
\item[(vi)]    $f$ is an isomorphism from $N_p\rest Y$ onto $N_q\rest \{b:
N_q \models b\in a\}$, 
\item[(vii)]   $f$ maps $u_p\cap Y$ into $u_q\cap\rng(f)$,
\item[(viii)]  $w_p\subseteq w_q$,
\item[(ix)]    if $n\in\dom(w^q)\backslash\dom(w^p)$ and $x\in u_p$, $N_p
\models$``$x$ is a function from the natural numbers to the natural numbers''
and $x^*=f(x)$ {\em then} $N_q \models$``if $y$ is the $n$-th natural number
then $w^q(y) > x(y)$''.  
\end{enumerate}
\end{enumerate}
Clearly $\bQ$ is equivalent to $\bQ'=(\mbox{the Hechler forcing})^{\bf L}$,
just let us define, for $p\in\bQ_1$, $g(p)=(w^p,F^p)$ where $F^p=\{h_p(x):
x\in M_p\}$. Now, $g$ is onto $\bQ'$ and
\[\begin{array}{ll}
\bQ_1\models p\leq q\quad\Rightarrow\ &\bQ'\models g(p)\leq g(q)\quad
\Rightarrow\\
\ &\neg(\exists p')(p\leq^\bQ p'\ \&\ p',q\mbox{ are
incompatible in }\bQ).
  \end{array}\]
The rest is left to the reader.  \QED$_{\ref{9.2}}$ 

\begin{proposition}
\label{9.3}
\begin{enumerate}
\item Assume that:
\begin{enumerate}
\item[(a)]  $\bar{\varphi}=(\varphi_0(x),\varphi_1(x,y))$ defines, in any
model of $\ZFC^-_*$, a forcing notion $\bQ_{\bar{\varphi}}$ with parameters
from ${\bf L}_{\omega_1}$, 
\item[(b)]  for every $\beta<\omega_1$ such that ${\bf L}_\beta\models
\ZFC^-_*$, for every $x,y\in {\bf L}_\beta$ we have: 
\[[x\in\bQ^{{\bf L}_\beta}_{\bar{\varphi}}\ \Leftrightarrow\ x\in \bQ^{{\bf
L}_{\omega_1}}_{\bar{\varphi}}]\quad\mbox{and}\quad [x< y\mbox{ in }\bQ^{{\bf
L}_\beta}_{\bar{\varphi}}\ \Leftrightarrow\ x < y\mbox{ in }\bQ^{{\bf
L}_{\omega_1}}_{\bar{\varphi}}],\] 
\item[(c)]  for unboundedly many $\alpha<\omega_1$ we have ${\bf L}_\alpha
\models \ZFC^-_*$,
\item[(d)]  any two compatible members of $\bQ^{{\bf L}_{\omega_1}}_{\bar{
\varphi}}$ have a lub,
\item[(e)]  like (c) for compatibility and for existence of lub.
\end{enumerate}
{\em Then}\ there is an $\aleph_0$--snep forcing notion $\bQ$ equivalent to
$\bQ^{{\bf L}_{\omega_1}}_{\bar{\varphi}}$: the $\Sigma^1_1$ (i.e.\ Souslin)
relations have just the real parameters of $\bar{\varphi}$.
\item We can use a real parameter $\rho$ and replace ${\bf L}_\alpha$ by ${\bf
L}_\alpha[\rho]$. 
\end{enumerate}
\end{proposition}

\Proof  It is similar to the proof of \ref{9.1}. Let $\bQ$ be the set of
quadruples $p=(E_p,n_p,\bar{\alpha}_p,\bar{a}_p)$ such that:
\begin{enumerate}
\item[$(\alpha)$] $E_p$ is a two-place relation on $\omega$,
\item[$(\beta)$]  $N_p\stackrel{\rm def}{=}(\omega,E_p)$ is a model of
$\ZFC^-_* +\V={\bf L}$, 
\item[$(\gamma)$] for some $n=n_p$ we have 
\[\bar{\alpha}_p=\langle \alpha_{p,\ell}:\ell<n\rangle,\qquad \bar{a}_p=
\langle a_{p,\ell}:\ell<n\rangle,\]
\item[$(\delta)$] $N_p\models$``$\alpha_{p,\ell}$ is an ordinal, $a_{p,\ell}
\in {\bf L}_{\alpha_{p,\ell}}$, ${\bf L}_{\alpha_{p,\ell}}\models\ZFC^-_*$,
and for $k\leq\ell<n$ we have ${\bf L}_{\alpha_{p,\ell}}\models\varphi_0(a_{p,
k}),\alpha_{p,\ell}<\alpha_{p,\ell +1}$'',\quad and 
\item[$(\varepsilon)$] if $m\leq k\leq\ell<n$ then $N_p\models$`` ${\bf L}_{
\alpha_{p,\ell}}\models\varphi_1(a_p,\alpha_m,a_{p,\alpha_k})$ ''.
\end{enumerate}
The order is given by:\qquad $p_0 \leq_{\bQ} p_1$ if and only if ($p_0,p_1\in
\bQ$ and) for some $Y_0,Y_1 \subseteq \omega$ and $f$ we have:
\begin{enumerate}
\item[(i)]  for $\ell=0,1$:\quad $Y_\ell$ is an $E_p$--transitive subset of
$N_p$,
\[(\forall x\in N_{p_\ell})(x\in Y_\ell\equiv\rk^{N_{p_\ell}}(x)\in Y_\ell),\]
\item[(ii)]  $f$ is an isomorphism from $N_{p_0}\rest Y_0$ onto $N_{p_1}\rest
Y_1$,
\item[(iii)] in $\{x\in N_{p_0}:N_{p_0}$``$x\mbox{ is an ordinal}\}$ there is
no $E_p$--minimal element, 
\item[(iv)]  $f$ maps $\{\alpha_{p_0,\ell}:\ell<n^*\}\cap Y_0$ into
$\{\alpha_{p_{1,\ell}}:\ell<n_{p_1}\}\cap Y_1$, 
\item[(v)]   if $f(\alpha_{p_0,k})=\alpha_{p_1,m}$ then $N_{p_1}\models$``
${\bf L}_{\alpha_{p_1,m}}\models\varphi_1(f(a_{p_0,k}),a_{p_1,m})$ ''.
\end{enumerate}

\begin{claim}
\label{factA}
$\bQ$ is a quasi order.
\end{claim}

\noindent{\em Proof of the claim:}\qquad Check.

Now define $M_p,h_p,M_p'$ as in the proof of \ref{9.2}.

\begin{claim}
\label{factC}
The set 
\[\bQ'\stackrel{\rm def}{=}\{p\in\bQ:N_p\mbox{ is well founded, }n_p>0\}\]
is dense in $\bQ$. 
\end{claim}

\noindent{\em Proof of the claim:}\qquad Check.

Define $g:\bQ'\longrightarrow\bQ^{{\bf L}_{\omega_1}}_{\bar{\varphi}}$ by
$g(p)=h_p(a_{p,n_p-1})$. 

\begin{claim}
\label{factE}
$g$ is really a function from $\bQ'$ onto $\bQ^{{\bf L}_{\omega_1}}_{
\bar{\varphi}}$ and 
\[\begin{array}{l}
p_0\leq_{\bQ} p_1\quad\Rightarrow\quad \bQ^{{\bf L}_{\omega_1}}_{\bar{\varpi}}
\models g(p_0)\leq g(p_1)\quad\Rightarrow\\
\mbox{[if }p_1\leq_\bQ p_2\mbox{ then for some $p_3$ we have $p_2
\leq_\bQ p_3$ and $p_0\leq_\bQ p_3$].}
  \end{array}\]
\end{claim}

\noindent{\em Proof of the claim:}\qquad The first implication is immediate
(by clause (v) in the definition of $\leq_\bQ$. For the second implication
assume $\bQ^{{\bf L}_{\omega_1}}_{\bar{\varphi}}\models g(p_0)\leq g(p_1)$ and
let $p_1\leq_\bQ p_2$. For $\ell=0,1,2$ let 
\[n_\ell=\min\{n: n=n_p\mbox{ or }n<n_{p_\ell}\mbox{ and }\alpha_{p_\ell,n}
\notin M_p\}.\]
Let $p_3$ be defined as follows: $M_p={\bf L}_\gamma$, ${\bf L}_\gamma\models
\ZFC^-_*$, and $\gamma>M_{p_0}'\cap \omega_1,M_{p_1}'\cap \omega_1,M_{p_2}'
\cap \omega_1$. Let $g^*_\ell$ be the isomorphism from $M_{p_\ell}$ onto ${\bf
L}_{\gamma_\ell}$, $\gamma_\ell<\gamma$, and let $w=\{f_\ell(\alpha_{p_\ell,
m}): m<n_\ell, \ \ell<2\}$. List it as $\{\alpha_{p_3,k}: k<n_{p_3}\}$
(increasing enumeration) and let $\Upsilon=\{f_\ell(a_{p_\ell,m}):m<n_\ell,\
\ell<2\}$. Now, $f_2(a_{p_2,n_2-1})$ is a $\leq_{\bQ^{{\bf L}_{\omega_1}}_{
\bar{\varphi}}}$--upper bound of $\Upsilon$. Consequently, by clauses (d) and
(e) of the assumptions, we can define $a_{p_3,m}$ as required.
\QED$_{\ref{9.3}}$

\begin{proposition}
\label{9.4}
Assume that $\varphi=\varphi(x,y)$ is such that
\begin{enumerate}
\item[(i)]   $\ZFC^-_*\vdash$ for every infinite cardinal $x\in X\stackrel{\rm
def}{=}\{\alpha: \alpha=\omega$ or $\omega^\alpha=\alpha$ (ordinal
exponentiation) $\}$, there is a unique $A_x$, an unbounded subset of $x$ of
order type $x$ such that $\varphi(x,A_x)$, and $\psi(\cdot)$ defines a set
$S\subseteq X$ not reflecting,  
\item[(ii)]  $\ZFC^-_*\vdash$ if $\mu_1<\mu_2$ are from $X$ then $A_{\mu_1}
\nsubseteq A_{\mu_2}$,  
\item[(iii)] $\omega_1=\sup\{\alpha:{\bf L}_\alpha\models\ZFC^-_*\}$, and the
truth value of ``$\beta\in A_\gamma,\ \beta\in S$'' is the same in ${\bf
L}_\alpha$ for every $\alpha<\omega_1$ for which ${\bf L}_\alpha\models
\ZFC^-_*$,
\item[(iv)]  the set $S$, i.e.\ $\{\beta<\omega_1: (\exists\alpha)({\bf
L}_\alpha\models\ZFC^-_*\ \&\ \psi(\beta))\}$, is a stationary subset of
$\omega_1$

[a kind of ``$\aleph^{\V}_1$ is below first ineffable of ${\bf L}$'' and is
not weakly compact]. 
\end{enumerate}
{\em Then} for some $\bar{\varphi}$ as in the assumptions of \ref{9.3}, and
$\name{\eta}$ we have:
\begin{enumerate}
\item[(a)] $\bQ^{{\bf L}_{\omega_1}}_{\bar{\varphi}}$ is a c.c.c.\ forcing
notion, 
\item[(b)] $\name{\eta}\in\can$ is a generic real of $\bQ^{{\bf L}_{\omega_1}
}_{\bar{\varphi}}$, and is nowhere essentially Cohen, 
\item[(c)] $\bQ^{{\bf L}_{\omega_1}}_{\bar{\varphi}}$ commute with Cohen.
\end{enumerate}
\end{proposition}

\Proof Let $pr(\alpha,\beta)=(\alpha+\beta)(\alpha+\beta)+\alpha$, it is a
pairing function. By coding, without loss of generality (e.g.\ letting 
\[A'_\alpha=\{pr^+(n,pr_n(\beta_1,\ldots,\beta_n)):n<\omega,\{\beta_1,\ldots,
\beta_n\} \subseteq B_\alpha\},\]
where $pr_1(\beta)=\beta$, $pr_{n+1}(\beta_1,\ldots,\beta_{n+1})=pr(pr_n(
\beta_1,\ldots,\beta_n),\beta_{n+1})$)
\begin{enumerate}
\item[(ii)$'$] if $x,x_1,\ldots,x_n$ are distinct cardinals in ${\bf L}_{
\omega_1}$, then $A_x \nsubseteq\bigcup\limits^n_{\ell =1} A_{x_\ell}$.
\end{enumerate}
For $\delta\in X$ let $f_0(\delta)=\min(X\setminus (\delta+1))$ and let
$f^1_\delta$ be the first (in the canonical well ordering of ${\bf L}$)
one-to-one function from $f_0(\delta)$ onto $\delta$. Let $C_\delta$ be the
first club of $\delta$ disjoint to $S$. For $\alpha\in [\omega,\omega_1)$, let
$\delta_\alpha=\max(X\cap\delta)$ and let
\[B^*_\alpha=\{pr_3(\varepsilon,\zeta,\xi):\varepsilon\in C_{\delta_\alpha},\
\zeta=f^1_{\delta_\alpha}(\alpha),\ \xi\in A'_{\delta_\alpha}\mbox{ and }
\varepsilon>\zeta,\ \varepsilon>\xi\}.\]
Note that
\begin{enumerate}
\item[$(*)$] $B^*_\alpha$ is an unbounded subset of $\delta_\alpha$ such that 
\begin{enumerate}
\item[(a)] $\beta\in S\cap\alpha\quad\Rightarrow\quad\beta>\sup(B^*_\alpha\cap
\beta)$,
\item[(b)] if $\alpha_1,\ldots,\alpha_n\in [\omega,\omega_1)\setminus\{
\alpha\}$ then $B_\alpha\setminus \bigcup\limits_{\ell=1}^n B^*_{\alpha_\ell}$
is unbounded in $\delta_\alpha$.
\end{enumerate}
\end{enumerate}
[Why? For (a), suppose that $\beta\in S\cap\alpha$. Trivially, $\min(
B^*_\alpha)>\min(C_\alpha)$, so $\gamma=\sup(C_{\delta_\alpha}\cap\beta)$ is
well defined. Now, 
\[B^*_\alpha\cap\beta\subseteq \{{\rm pr}_3(\varepsilon,\zeta,\xi):
\varepsilon,\zeta,\xi\leq\gamma\}\subseteq (\gamma+\gamma+\gamma)^3<\beta\]
(the last inequality follows from the fact that $\beta\in X$). To show (b)
suppose that $\gamma_0<\delta_\alpha$ and choose $\xi\in B^*_\alpha\setminus
\bigcup\limits_{\ell=1}^n B^*_{\alpha_\ell}$. Let $\zeta=f^1_{\delta_\alpha}
(\alpha)$ and let $\varepsilon\in C_{\delta_\alpha}$ be large enough. So ${\rm
pr}_3(\varepsilon,\zeta,\xi)\in B^*_\alpha$ (by definition) and ${\rm pr}_3(
\varepsilon,\zeta,\xi)\notin B^*_{\alpha_\ell}$ (use the third coordinate) and
${\rm pr}_3(\varepsilon,\zeta,\xi)>\varepsilon>\gamma_0$.]

Let $I_\alpha$ be the ideal of subsets of $B^*_\alpha$ generated by
\[\{B^*_\alpha\cap B^*_\beta:\omega\leq\beta<\omega_1,\ \beta\neq\alpha\}\cup
\{B^*_\alpha\cap\beta:\beta<\delta_\alpha\}.\]
Let $\bQ$ be the set of finite functions $p$ from $\omega_1\setminus\omega$ to
$\{0,1,2\}$ ordered by:\\
$p \le q$ if and only if: 
\begin{quotation}
if $\alpha\in\dom(p)$, $\beta\in\dom(q)\cap A_\alpha\backslash\dom(p)$\\ 
then $q(\beta)=p(\alpha)$ and $\beta>\sup(\delta_\alpha\cap\dom(p))\ \vee\
q(\beta=2)$. 
\end{quotation}

\begin{claim}
\label{fact2A}
$\bQ$ is a partial order.
\end{claim}

\begin{claim}
\label{fact2B}
For each $\alpha\in [\omega,\omega_1)$ the set ${\cal I}_\alpha=\{p:\alpha\in
\dom(p)\}$ is dense in $\bQ$.
\end{claim}

\noindent{\em Proof of the claim:}\qquad Let $p\in\bQ$ and suppose that
$\alpha\notin \dom(p)$. Let $q=p\cup\{\langle\alpha,2\rangle\}$.
\medskip

Let $\name{f}$ be the $\bQ$--name defined by $\forces\name{f}=\bigcup
\name{G}_\bQ$. 

\begin{claim}
\label{fact2D}
For $\alpha\in [\omega,\alpha)$,
\[\begin{array}{ll}
\forces_{\bQ}&\mbox{`` for some $\ell< 3$, for any $m<3$ we have}\\
\ &\{\beta\in B^*_\alpha:\name{f}(\beta)=m\}\neq\emptyset\mod {\cal I}_\alpha\
\mbox{ iff }\ m\in\{0,\ell\}\mbox{ ''.}
  \end{array}\]
\end{claim}

\noindent{\em Proof of the claim:}\qquad Take $p\in\name{G}_\bQ$ such that
$\alpha\in\dom(p)$ and let 
\[B=\bigcup\{B^*_\alpha\cap B^*_\gamma:\gamma\in\dom(p)\setminus\{\alpha\}\},\]
so $B\in {\cal I}_\alpha$. Clearly, $p\forces$`` if $\beta\in B^*_\alpha
\setminus B$ then $\name{f}(\beta)\in\{2,p(\alpha)\}$ '', hence 
\[p\forces_\bQ\mbox{`` if }m\in\{0,1,2\}\setminus\{2,p(\alpha)\}\mbox{ then }
m\notin\rng(\name{f}\rest (B^*_\alpha\setminus B))\mbox{ ''.}\]
Now, if $B'\in {\cal I}_\alpha$, and $p\leq_\bQ q$ then there is $\gamma\in
A_\alpha\setminus B'\setminus\bigcup\{B^*_\gamma: \gamma\in\dom(p)\setminus\{
\alpha\}\}$ such that $q\cup\{\langle\gamma,2\rangle\}$ and
$q\cup\{\langle\gamma,p(\alpha)\rangle\}$ are in $\bQ$ above $q$. Reflecting
we are done. 

\begin{claim}
\label{fact2E}
One can define $\name{f}$ from $\name{f}\rest\omega\in {}^{\textstyle
\omega}3$. 
\end{claim}

\noindent{\em Proof of the claim:}\qquad Define $\name{f}\rest\alpha$ by
induction on $\alpha\in X$ using \ref{fact2D}.

\begin{claim}
\label{fact2F}
The forcing notion $\bQ$ is nowhere essentially Cohen.
\end{claim}

\noindent{\em Proof of the claim:}\qquad For every $\alpha^*<\omega_1$ and
for every large enough $\gamma<\omega_1$ and for $\ell\in \{0,1,2\}$, the
condition $q_\ell=\{\langle\gamma,\ell\rangle\}$ is compatible with every
$q\in\bQ$ such that $\dom(p)\subseteq\alpha^*$.

\begin{claim}
\label{fact2G}
The $\bQ$--name $\name{f}$ (for a real) is nowhere essentially Cohen.
\end{claim}

\noindent{\em Proof of the claim:}\qquad By \ref{fact2E}, \ref{fact2F}. 

\begin{claim}
\label{fact2H}
The forcing notion $\bQ$ satisfies the demands in \ref{9.4}.
\end{claim}

\noindent{\em Proof of the claim:}\qquad Check.

\begin{claim}
\label{fact2I}
The forcing notion $\bQ$ satisfies the c.c.c.
\end{claim}

\noindent{\em Proof of the claim:}\qquad Use ``$S\subseteq\omega_1$ is
stationary''. \QED$_{\ref{9.4}}$

\begin{remark}
\begin{enumerate}
\item Of course, such forcing can make $\aleph_1$ to be $\aleph^{{\bf
L}[\eta]}_1$. But it seems that we can have such forcing which preserves the
${\bf L}_{\omega_1}$--cardinals (and even their being ``large'' in suitable
senses). For this it should be like ``coding the universe by a real'' of
Jensen Beller Welch \cite{BJW}, and see Shelah Stanley \cite{ShSt:340}.  
\item Instead of coding $\aleph_1$--Cohen we can iterate adding dominating
reals or whatever.
\end{enumerate}
\end{remark}

\begin{definition}
\label{9.6}
\begin{enumerate}
\item We say that forcing notions $\bQ_0,\bQ_1$ are equivalent if their
completions to Boolean algebras ($\BA(\bQ_0), \BA(\bQ_1)$) are isomorphic. 
\item Forcing notions $\bQ_0,\bQ_1$ are locally equivalent if
\begin{enumerate}
\item[(i)] for each $p_0\in\bQ_0$ there are $q_0,q_1$ such that
\[p_0\le q_0\in\bQ_0\ \&\ q_1\in\bQ_1\ \&\ \BA(\bQ_0\rest(\ge q_0))\cong\BA(
\bQ_1\rest(\ge q_1)),\] 
\item[(ii)] for every $p_1\in\bQ_1$ there are $q_0,q_1$ such that 
\[q_0\in\bQ_0\ \&\ p_1\le q_1\in\bQ_1\ \&\ \BA(\bQ_0\rest(\ge q_0))\cong\BA(
\bQ_1\rest (\ge q_1)).\]
\end{enumerate}
\end{enumerate}
\end{definition}

Now we may phrase the conclusions of \ref{9.3}, \ref{9.4}.

\begin{proposition}
\label{9.5}
\begin{enumerate}
\item Assume $\bar{\varphi}_1=\langle \varphi^1_0,\varphi^1_1\rangle$ and
$\bar{\varphi}^2_2=\langle \varphi^2_0,\varphi^2_1\rangle$ are as in
\ref{9.3}. {\em Then}\ we can find $\bar{\varphi}^3_3$ as there, only with the
parameters of $\bar{\varphi},\bar{\varphi}_2$  and such that:
\begin{enumerate}
\item[(a)] if in ${\bf L}_{\omega_1}$ there is a last cardinal $\mu$ (i.e.\ 
$\aleph^{\V}_1$ is a successor cardinal in ${\bf L}$), then $\bQ^{{\bf
L}_{\omega_1}}_{\bar{\varphi}_3}$ is locally equivalent to 
\[\bigcup\{\bQ^{{\bf L}_\alpha}_{\bar{\varphi}_1}:\mu<\alpha,\ {\bf L}_\alpha
\models\mu\mbox{ is the last cardinal}\},\]
\item[(b)] if in ${\bf L}_{\omega_1}$ there is no last cardinal (i.e.\
$\aleph^\V_1$ is a limit cardinal in ${\bf L}$), then $\bQ^{{\bf
L}_{\omega_1}}_{\bar{\varphi_3}}$ is locally equivalent to 
\[\bigcup\{\bQ^{{\bf L}_{\omega,\alpha}}_{\bar{\varphi}^1_2}:{\bf
L}_{\omega_1} \models\alpha\mbox{ a cardinal}\}.\]
\end{enumerate}
\item  In \ref{9.3}, \ref{9.4}(1) we can replace ${\bf L}_{\omega_1}$ by ${\bf
L}_{\omega_1}[\eta^*]$, $\eta^* \in\baire$. 
\item In \ref{9.3}, \ref{9.4}(1) we can replace ${\bf L}_{\omega_1}$ by${\bf
L}_{\omega_1}[A]$ where $A \subseteq \omega_1$ but have $\aleph_1$-snep
instead of $\aleph_0$-snep. 
\end{enumerate}
\end{proposition}

\Proof Let $\varphi_{3,0}(x)$ say
\begin{enumerate}
\item[(i)] $x=\langle\bar{\alpha}^x,\bar{\beta}^x,\bar{a}^x,\bar{b}^x\rangle$,
$\bar{\alpha}^x=\langle\alpha^x_\ell:\ell\le n^x\rangle$, $\langle\langle
\beta^x_{\ell,k}:k \le k^x_\ell\rangle:\ell\le n^x\rangle$, $\bar{a}^x=
\langle\langle a^x_{\ell,k}:k \le k^x_\ell \rangle:\ell<n^x\rangle$,
$\bar{a}^x=\langle a^x_\ell:\ell \le n^x \rangle$, $\bar{b}^x=\langle
b^x_\ell:\ell<n^x\rangle$, 
\item[(ii)] $\alpha^x_\ell<\beta^x_{\ell,0}<\beta^x_{\ell,1}\ldots$,  ${\bf
L}_{\beta^x_\ell}\models$``$\alpha^\ell_x$ the last cardinal'', 
\item[(iii)] ${\bf L}_{\beta^x_\ell}\models$``$\varphi_{1,0}(b^x_\ell)$'',
${\bf L}_{\alpha^x_{\ell +1}}\models$``$\varphi_{2,0}(a^x_\ell)$'' for $\ell <
n^x$, 
\item[(iv)] ${\bf L}_{\beta^x_{n^x}}\models$``$\alpha^x_\ell$ is a cardinal'',
\item[(v)]  ${\bf L}_{\alpha^x_{\ell +1}}\models$``$\varphi^2_1(a^x_\ell,
a^x_{\ell+1})$''. 
\end{enumerate}
Let $\beta(x)=x$. Let $\varphi_{3,1}(x,y)$ say:
\begin{enumerate}
\item[$(\alpha)$] $\beta^x_{n^x} \le \beta^y_{n_y}$,
\item[$(\beta)$]  $\{\alpha^x_\ell:\ell \le n^x$ and ${\bf L}_{\beta^y_y}
\models$``$\alpha^x_\ell$ is a cardinal''$\}$ is a subset of
$\{\alpha^y_\ell:\ell \le n^y\}$, 
\item[$(\gamma)$] if $\alpha^x_{\ell(*)}$ is maximal in $\{\alpha^x_\ell:\ell
< n^x, {\bf L}_{\beta^y_{n^y}}\models$``$\alpha^x_\ell$ is a cardinal''$\}$
then $\bar{\alpha}^x\rest \ell(*)=\bar{\alpha}^y \rest\ell(x)$, $\bar{\beta}^x
\rest\ell(*)=\bar{\beta}^y\rest\ell(*)$, $\bar{a}^x\rest\ell(*)=\bar{a}^y
\rest \ell(*)$, $\bar{b}^x\rest\ell(*)=\bar{b}^y\rest\ell(*)$,
\item[$(\delta)$] $\alpha^x_{\ell(*)}=\alpha^y_{\ell(*)}$,
\item[$(\varepsilon)$] $\beta^x_{\ell(*)}\le\beta^y_{\ell(*)}$ and ${\bf
L}_{\beta^y_{\ell(*)}}\models \varphi^1_{0,1}(b^x_{\ell(*)},b^y_{\ell(*)})$. 
\end{enumerate}
Now check. \QED$_{\ref{9.5}}$

\stepcounter{section}
\subsection*{\quad 11. Preservation of ``no dominating real''}
The main result of \S7: (for homogeneous c.c.c.\ $\bQ$) if a nep forcing $\bP$
preserves $(\baire)^\V\in (I^{\ex}_{(\bQ,\name{\eta})})^+$ {\em then}\ it
preserves $X\in (I^{\ex}_{(\bQ,\name{\eta})})^+$ (see \ref{6.5}) is a case of
the following 

\begin{thesis}
\label{6A.1}
Nep forcing notions do not discern sets $X\subseteq\baire$ built by
diagonalization, say between $X,Y\subseteq\baire$ which are generic enough.
\end{thesis}
But there are interesting cases not covered by \ref{6.5}, most prominent is: 

\begin{question}
\label{6A.2}
If a nep forcing notion preserves ``$F\subseteq\baire$ is unbounded'' for some
(unbounded) $F\subseteq\baire$ {\em then} does it preserve this for every
(unbounded) $F'\subseteq \baire$?
\end{question}

\begin{definition}
\label{6A.3}
\begin{enumerate}
\item For a Borel 2-place relation ${\cal R}$ on $\baire$ let $I_{\cal R}$ be
the $\aleph_1$--complete ideal on $\baire$ generated by the sets of the form
$A_\nu=\{\eta\in\baire:\neg(\eta R \nu)\}$.
\item  We say that $\nu$ is ${\cal R}$--generic over $N$ if $\eta\in N\cap
\baire\ \Rightarrow\ \eta R\nu$. 
\item  A forcing notion $\bP$ is weakly ${\cal R}$--preserving if for any
$\eta_0,\eta_1,\ldots,\eta_n,\ldots\in (\baire)^{\V^\bP}$ there is $\nu \in
(\baire)^\V$ such that $n<\omega\ \Rightarrow\ \eta_n {\cal R}\nu$ (i.e.\
$\forces_{\bP} (\baire)^\V \in I^+_{\cal R}$). 
\item We say that a forcing notion $\bP$ is ${\cal R}$--preserving if for any
Borel subset $B$ of $\baire$ from $\V$ which is in $I^+_{\cal R}$, for any
$\eta_0,\eta_1,\ldots,\eta_n,\ldots \in (\baire)^{\V^\bP}$ there is $\nu\in
B^\V$ such that $n<\omega\ \Rightarrow\ \eta_n {\cal R}\nu$ (i.e.\
$\forces_{\bP} B^\V\in I^+_{\cal R}$).
\item We say that a forcing notion $\bP$ is strongly ${\cal R}$--preserving if
for any $X\in I^+_{\cal R}$ (in $\V$) we have $\forces_{\bP}$``$X\in I^+_{\cal
R}$''. 
\item We say that a forcing notion $\bP$ is super ${\cal R}$--preserving as
witnessed by $({\frak B},\ZFC^-_{**})$ if 
\begin{enumerate}
\item[(a)] every $({\frak B},\ZFC^-_*)$--candidate is a $\bP$--candidate,
\item[(b)] for any $({\frak B},\ZFC^-_{**})$--candidate $N$ such that $p\in N$
and for any $\nu\in\baire$ which is ${\cal R}$--generic over $N$, {\em there
is} $q$ such that $p\le_{\bP} q\in \bP$, $q$ is $\langle N,\bP
\rangle$--generic and $q\forces$``$\nu$ is ${\cal R}$--generic over $N\langle
\name{G}_{\bP}\rangle = N\langle\name{G}_{\bP}\cap P^N \rangle$''.
\end{enumerate}
\end{enumerate}
\end{definition}

\begin{proposition}
\label{6A.4}
\begin{enumerate}
\item If $\bP$ is super ${\cal R}$--preserving, {\em then} $\bP$ is strongly
${\cal R}$--preserving (also for $\bP$ nep). 
\item If $\bP$ is strongly ${\cal R}$--preserving, {\em then} $\bP$ is ${\cal
R}$--preserving. 
\item If $\bP$ is ${\cal R}$--preserving, {\em then} it is weakly ${\cal
R}$--preserving. 
\end{enumerate}
\end{proposition}

\Proof Easy. \QED

\begin{proposition}
\label{6A.5}
A sufficient condition for ``$\bP$ is super ${\cal R}$--preserving as
witnessed by $({\frak B},\ZFC^-_{**})$'' is that for some nep forcing notion
$\bQ$ and $\hc$--$\bQ$-name $\name{\eta}^*$ and a Borel relation ${\cal R}_1$
we have 
\begin{enumerate}
\item[$(\alpha)$] every $({\frak B},\ZFC^-_{**})$--candidate is a
$\bQ$--candidate and $\name{\eta}^*\in N$ and also it is a $\bP$--candidate,
\item[$(\beta)$]  if $N$ is a $({\frak B},\ZFC^-_{**})$--candidate (so
countable) and $\nu$ is ${\cal R}$--generic over $N$ then for some  $G_{\bQ}
\subseteq \bQ^N$ generic over $N$, in $\V$ we have 
\[(\forall x)(x\;{\cal R}_1\;\name{\eta}^*[G_{\bQ}])\ \Rightarrow\ x\; {\cal
R}_1\; \nu)\]
and for every $G_{\bQ}\subseteq\bQ^N$, generic over $N$, and $p\in\bP^N$ there
is $q$ such that $p\le q\in\bP$ and $q$ is $(N,\bP)$--generic and
\[q\forces_{\bP}\mbox{`` }\name{\eta}^*[G_Q]\mbox{ is ${\cal R}_1$--generic
over }N \langle G_{\bP}\cap\bP^N \rangle\mbox{ ''}.\]
\end{enumerate}
\end{proposition}

\Proof Straight. \QED

\begin{remark}
\label{6A.6}
In \ref{6A.7} below we phrase a sufficient condition. Note that clause
$(\delta)$ can be naturally phrased as ``an appropriate sentence $\psi$
follows from $\ZFC^-_{**}$; this is slightly stronger as possibly $\psi$ holds
only for all $({\frak B},\ZFC^-_{**})$--candidates but not for some (e.g.\
non-well founded) models of $\ZFC^-_{**}$ (this does not matter).
\end{remark}

\begin{proposition}
\label{6A.7}
Assume that:
\begin{enumerate}
\item[$(\alpha)$] $\bP,\bQ$ are nep forcing notions, $\name{\eta}^*$ is a
$\hc$--$\bQ$--name, ${\cal R},{\cal R}_1$ are Borel relations,
\item[$(\beta)$]  every $({\frak B},\ZFC^-_{**})$--candidate $N$ is a
$\bQ$--candidate and $\bP$-candidate, $\name{\eta}^*\in N$, 
\item[($\gamma)$] for every $({\frak B},\ZFC^-_{**})$--candidate $N$ and $\nu\in
\baire$ which is ${\cal R}$--generic over $N$ and $r\in\bQ^N$ we can find $G_Q
\subseteq\bQ^N$ generic over $N$ such that 
\[r\in G_Q\quad\mbox{ and }\quad(\forall x)(x\; {\cal R}\; \nu\ \ \Rightarrow\
\ x\; {\cal R}_1\;\name{\eta}^*[G_Q]),\]
\item[$(\delta)$] if $N$ is a $({\frak B},\ZFC^-_{**})$--candidate then for
every $G_{\bQ}\subseteq\bQ^N$ generic over $N$ and $G_R\subseteq\Levy(
\aleph_0,2^{|\bP|} + 2^{|\bQ|})^N$ generic over $N[G_Q]$ we have 

$N[G_Q][G_R]\models$`` there are $G_\bQ$ and $q,p\le q\in\bP$ such that $q$ is
explicitly $(N,\bP)$--generic and $G_\bQ$ is a generic over $N$ subset of
$\bQ^N$ and $q\forces_{\bP}$``$\name{\eta}^*[G_{\bQ}]$ is ${\cal
R}_1$--generic over $N[G_{\bQ},\name{G}_{\bP}]$''. 
\end{enumerate}
{\em Then} ``$\bP$ is super ${\cal R}$--preserving as witnessed by
$\ZFC^-_{**}$'' . 
\end{proposition}

\Proof  Clause $(\alpha)$ of \ref{6A.5} holds by clause $(\beta)$ here. Next,
the first demand in clause \ref{6A.5}$(\beta)$ (``for some $G_\bQ\subseteq
\bQ^N$ generic over $N$'') follows from clause $(\gamma)$ of
\ref{6A.7}. Finally, suppose that $G_\bQ\subseteq\bQ^N$ is generic over $N$,
equivalently, $\nu^*$ is a $\bQ$--generic real. Let $G\subseteq\Levy(\aleph_0,
|2^{\bP}|^N)$ be generic over $N$, equivalently over $N[\nu^*]$. In $N$, by
clause (e) of the assumptions, in $N[G]$, there is a semi $\bP$--candidate
$M$, ${\cal P}(\bP^N)^N ={\cal P}(\bP^M)^M$. So in $N[G]$, $M[G]$ is a
$\bP$--candidate. So there is $q\in\bP^{N[G]}$ such that $N[G]\models$``$p\le
q$ and $q$ is $(M,\bP^M)$--generic''. As above possibly increasing $q$, 
\begin{enumerate}
\item[$(\circledast)$] $N[G]\models[q\forces_{\bP}$`` $\nu$ is ${\cal
R}$--generic over $M[\name{G}_{\bP}]$ '' and $\nu$ is Cohen over $M]$.
\end{enumerate}
So for some $\name{\nu}$, 
\[\begin{array}{ll}
N \models\mbox{``}&\name{\nu},\name{q}\mbox{ are $\Levy(\aleph_0,|2^{(\bP)}
|^N)$--names of a Cohen real and}\\
\ &\mbox{  a member of $\bP$, respectively, and some }r\in\Levy(\aleph_0,(
|2^{|\bP|}|)^N)\\ 
\ &\mbox{ forces the statement $(\circledast)$ above on }\name{q}, 
\name{\nu}\mbox{ ''}.
  \end{array}\]  
Now we can find $G'\subseteq\Levy(\aleph_0,(2^{|\bP|})^N)$ generic over $N$ to
which $r$ belongs and $\name{\nu}[G']=^*\nu$ (i.e.\ they are equal except for
finitely many coordinates). Let $q'\in\bP$ be $\ge\name{q}[G']$ and be
$(N\langle G'\rangle,\bP^{N\langle G'\rangle})$--generic, so we are
done. \QED$_{\ref{6A.7}}$ 

\begin{theorem}
\label{6A.8}
Assume that:
\begin{enumerate}
\item[(a)] ${\cal R}$ is:\quad $f {\cal R} g$ iff $g$ is non-decreasing and
$(\exists^\infty n)(f(n)>g(n))$;\\
${\cal R}_1$ is:\quad  $f {\cal R}_1 g$ iff $(\exists^\infty n)
(f(n)>\max\{g(m):m \le n\})$, 
\item[(b)] $\bQ=(\{\eta\in {}^{\omega >}\omega:\eta$ non-decreasing$\})$,
$\name{\eta}^*$ is the generic real $\bigcup\name{G}_{\bQ}$ (so really $\bQ$
is the Cohen forcing),
\item[(c)] $\bP$ is nep,
\item[(d)] every $({\frak B},\ZFC^-_{**})$--candidate $N$ is a
$\bP$--candidate (and easily it is a $\bQ$--candidate), 
\item[(e)] $\ZFC^-_{**}$ says: ``$\bP$ is nep, ${\cal P}(\bP)\in {\cal
H}(\chi)$, ${\cal H}(\chi)$ is a semi $\bP$--candidate and after forcing with
$\Levy(\aleph_0,2^{|\bP|} + 2^{|\bQ|})$ still is, and forcing with $\bP$ does
not add a dominating real''.
\end{enumerate}
{\em Then} the conditions $(\alpha)$--$(\delta)$ of \ref{6A.7} hold.
\end{theorem}

\Proof Let $N$ be a $\bP$--candidate and let $q\in\bQ$. Now, for any
$\nu\in\baire$ 
\begin{enumerate}
\item[$\bigotimes_1$] there are $g_1,g_2$ such that 
\begin{enumerate}
\item[(a)] $g_\ell$ is a subset of $\bQ^N$ generic over $N$ to which $q$
belongs; and let $\eta^*_\ell=\name{\eta}^*[g_\ell]\in\baire$,
\item[(b)] $m\in [\lh(q),\omega)\ \Rightarrow\ \eta_1(m)<\nu(m)\vee\eta_2(m) <
\nu(m)\vee\nu(m)=0$.
\end{enumerate}
\end{enumerate}
[Why?  Quite easy, letting $\langle {\cal I}_k:k<\omega\rangle$ list the dense
open subsets of $\bQ^N$ in $N$, we choose inductively $m_k$ and $(\eta^1_k,
\eta^2_k)$ such that $\eta^1_k\in {}^{(m_k)} \omega$, $\eta^2_k\in {}^{(m_k)}
\omega$, $\eta^1_0=\eta^2_0=q$, $\eta^1_k\vartriangleleft\eta^1_{k+1}$,
$\eta^2_k \vartriangleleft \eta^2_{k+1}$, $\eta^1_{2k+1}\in {\cal I}_k$,
$\eta^2_{2k+2}\in {\cal I}_k$ and the demand in (b) is satisfied and
$\eta_1\stackrel{\rm def}{=}\bigcup\limits_{k<\omega}\eta^1_k$, $\eta_2
\stackrel{\rm def}{=}\bigcup\limits_{k<\omega}\eta^2_k$ induce $g_1,g_2$
respectively]. 

Next,
\begin{enumerate}
\item[$\bigotimes_2$] if $\nu$ is ${\cal R}$--generic over $N'$ (any
$\bP$--candidate), {\em then} $\eta^*_2$ is ${\cal R}$--generic over $N'$ or
$\eta^*_1$ is ${\cal R}$--generic over $N'$.
\end{enumerate}
[Why? Assume this fails, as $\eta^*_1$ is not ${\cal R}$--generic over $N'$
then some $f_1 \in (\baire)^M$ dominates $\eta^*_1$, and as $\eta^*_2$ is not
${\cal R}$--generic over $N'$ some $f_2 \in (\baire)^M$ dominates $\eta^*_2$,
so $f^*=\max\{f_1+1,f_2+1\}\in N'$ dominates $\nu$ (i.e.\ $f^*(m)=\max\{f_1(m)
+ 1,f_2(m) + 1\}$).] 
\begin{enumerate}
\item[$\bigotimes_3$] if $\nu \in \baire$ is non-decreasing ${\cal
R}$--generic over $N$ (so it is not dominated by $N$ and is non-decreasing),
$q\in\bQ$, {\em then} there is $G_{\bQ} \subseteq \bQ^N$ generic over $N$,
$q\in G_{\bQ}$ and $\name{\eta}^*[G_Q]\le^*\nu$.
\end{enumerate}
[Why? Let $\langle {\cal I}_n:n<\omega \rangle$ list the dense open subset of
$\bQ$ in $N$. We choose by induction on $n$, $q_n\in {}^{k_n}\omega\subseteq
\bQ$ such that $q_0 = 1$, $q_n\le q_{n+1}$, and $q_n \rest[k_0,k_n)\le\nu\rest
[k_0,k_n)$ and $q_{n+1}\in {\cal I}_n$. For $n=0$ trivial, for $n+1$ choose in
$N$ by induction on $\ell$, $m_{n,\ell},\rho_{n,\ell}$ such that $m_{n,0}=k_n$
and 
\[m_{n,\ell}<m_{n,\ell +1},\quad\rho_{m,\ell}\in {}^{[m_{n,\ell},m_{n,\ell
+1})}\omega,\quad q_n \cup 0_{[k_n,m_{n,\ell})} \cup \rho_{m,\ell} \in {\cal
I}_n.\] 
This is easy and $\langle \rho(m_{n,\ell},\rho_{m,\ell}):\ell<\omega\rangle
\in N$. Now define $\rho^*\in\baire$ by: 
\[\rho^* \rest [m_{n,\ell},m_{n,\ell +1})\mbox{ is constantly }\max\bigl(
\bigcup_{i\le \ell +1}\rng(\rho_i)\cup\rng(q_n)\bigr).\]
So $(\exists^\infty j < \omega)(\rho^*(j)<\nu(j))$ hence for some $\ell$ and
some $m \in [m_{n,\ell},m_{n,\ell +1})$ we have $\rho^*(m)<\nu(m)$. So
\[(\forall m')(m_{n,\ell +1}\le m'<m_{n,\ell +2}\ \Rightarrow\ \rho_{n,\ell
+1}(m')<\nu(m)),\]
but $\nu(m)<\min\{\nu(j):m_{n,\ell+1}\le j<m_{n,\ell +2}\}$, so we are done.]

Now we have to check the conditions in \ref{6A.7}, so obviously clauses
$(\alpha),(\beta)$ hold. Also clause $(\gamma)$ there holds by $\otimes_3$. So
let us prove clause $(\delta)$. Let $N$ be a $({\frak B},\ZFC^-_{
**})$--candidate and $p\in\bP^N$. Let $q$ be $\langle N,\bP \rangle$--generic,
$p\le q$ (by $\otimes_1 + \otimes_2$). Let $G\subseteq\Levy(\aleph_0,|2^{\bP}
|^N)$ be generic over $N$ (equivalently over $N[\nu^*]$. In $N$, by clause (e)
of the assumptions, in $N[G]$, there is a semi $\bP$--candidate $M$, ${\cal
P}(\bP^N)^N ={\cal P}(\bP^M)^M$. Then in $N[G]$, $M[G]$ is a
$\bP$--candidate. So there is $q\in\bP^{N[G]}$ such that $N[G]\models$``$p\le
q$ and $q$ is $(M,\bP^M)$--generic''. As above possibly increasing $q$, 
\[N[G]\models[q\forces_{\bP}\mbox{`` }\nu\mbox{ is ${\cal R}$--generic over }
M[\name{G}_{\bP}]\mbox{''}\mbox{ and $\nu$ is Cohen over }M].\]
\QED$_{\ref{6A.8}}$

\begin{remark}
\label{6A.11}
Clearly this proof is similar to \S7, so we can replace ``$\Cohen$'' by more
general $\bQ$. More exactly, the point is that in \S7 the demand was
$q\forces_{\bP}$``$\name{\eta}$ is $\bQ$--generic over $N$''. Here we replace
it by other demands. 
\end{remark}

\begin{conclusion}
\label{6A.9}
For any Souslin proper forcing notion $\bP$, if $\bP$ add no dominating real,
{\em then} forcing with $\bP$ adds no member of $\baire$ dominating some $F
\subseteq\baire$ from $\V$ not dominated there.
\end{conclusion}

\Proof  By \ref{6A.8}, \ref{6A.6}, \ref{6A.7}. \QED

\begin{conclusion}
\label{6A.10}
\begin{enumerate}
\item Suppose that
\begin{enumerate}
\item[(a)] $\bP$ is a forcing notion adding no dominating real,
\item[(b)] $\nbQ$ is a $\bP$--name for a Souslin proper forcing notion not
adding a dominating real. 
\end{enumerate}
{\em Then} $P * \nbQ$ adds no dominating real.
\item $P * \nbQ$ adds no real dominating an old undominated family {\em if}\
both $\bP$ and $\nbQ$ satisfy this and are Souslin proper.
\end{enumerate}
\end{conclusion}

\Proof By \ref{6A.9}. \QED

\stepcounter{section}
\subsection*{\quad 12. Open problems} 

\begin{problem}
\label{20.1}
\begin{enumerate}
\item Can we in \cite{Sh:480} weaken the assumptions (from Souslin c.c.c.)
to ``$\bQ$ is nep and c.c.c.''? 
\item Similarly in the symmetry theorem.
\item Similarly other problems here have such versions too.
\end{enumerate}
\end{problem}

\begin{problem}
\label{20.3x}
\begin{enumerate}
\item {\em (von Neumann)} Is it consistent that every c.c.c.\
$\baire$--bounding atomless forcing notion is a measure algebra? We may now
rephrase: is the non-existence consistent?
\item {\em (Velickovic)} Is it consistent that every c.c.c.\ forcing notion
adding new reals adds a real $\name{f}\in\baire$ such that 

if $S\in \prod\limits_{n<\omega} [\omega]^{\textstyle 2^n}\cap \V$ then
$(\exists^\infty n\in\omega)(\name{f}(n)\notin S(n))$.

[Note that \cite{Sh:480} answers a relative of \ref{20.3x}(2): there is no
such Souslin c.c.c.\ forcing notion.]
\end{enumerate}
\end{problem}

A relative of the von Neumann problem is a problem which Fremlin \cite{Fe94}
stresses and has many equivalent versions (see \cite{Fe94} on its
history). Half way between them and our context is the following.

\begin{problem}
\label{20.2}  Assume $\bQ$ is a Souslin c.c.c.\ $\baire$--bounding forcing
notion. Is it random forcing? 
\end{problem}

\begin{problem}
\label{20.3}
\begin{enumerate}
\item Is it consistent that every c.c.c.\ forcing notion adding an
unbounded real adds a Cohen real?  (See B{\l}aszczyk Shelah \cite{Sh:F151} for
a proof of the $\sigma$-centered version).
\item  If $\bP$ satisfies \cite[1.5]{Sh:480}, does it imply $\bP$ adds a Cohen
real?
\end{enumerate}
\end{problem}

\begin{problem}
\label{20.4}
Are there any symmetric (or $(<\omega)$-symmetric) c.c.c.\ Souslin forcing
notions in addition to Cohen forcing and random forcing?

[``Yes'' here implies ``no'' to \ref{20.2} so not of present interest.]
\end{problem}
 
\begin{problem}
[Gitik Shelah \cite{GiSh:357}, \cite{GiSh:412}]
\label{20.5}
\begin{enumerate}
\item Assume $I$ is an $\aleph_1$--complete ideal on $\kappa$ such that ${\cal
P}/I$ is atomless.  Can $I^+$ (as a forcing notion) be a c.c.c.\ Souslin forcing
generated by a real.
\item Replace Souslin by ``definable in an $({\cal H}_{<\sigma}(\theta),\in,
{\frak B})$, ${\frak B}$ has universe $\kappa$ or ${\cal
H}_{<\sigma}(\kappa)$, and $I$ is $(\theta+\kappa)^+$--complete (see
\cite{GiSh:357}).
\item Generalize the results of the form ``if ${\cal P}(\kappa)/I$ is the
measure algebra with Maharam dimension $\mu$ (or is the adding of $\mu$ Cohen
reals) then $\lambda$ is large enough'', see \cite{GiSh:412}, \cite{GiSh:582}
for those results.  
\item Combine (2) and (3).
\end{enumerate}
\end{problem}

\begin{problem}
[Judah]
\label{20.6}
Can a Souslin c.c.c.\ forcing notion add a minimal real?  (Note: this is of
interest only if the answer in \ref{20.2} is NO and/or the answer to
\ref{20.10} is NO.).
\end{problem}

\begin{problem}
\label{20.7}
Give examples of a Souslin forcing notion which is only temporarily c.c.c.\
and/or proper (${\bf L}$) (see \S10).
\end{problem}

\begin{problem}
\label{20.8}
Do iterations (CS,FS) of Souslin c.c.c.\ forcing notions not adding a
dominating real have this property?  Is each almost $\baire$--bounding?  
[Maybe \ref{6.5} answers need better: replace $\eta^*$ is generic real for
$(N,\bQ,\name{eta})$ by less].  See \S11 + \S10.
\end{problem}

\begin{problem}
\label{20.9}
\begin{enumerate}
\item Is there a pair $(\bQ,\name{r})$ such that:
\begin{enumerate}
\item[(a)] $\forces_{\bQ}$``$\name{r}\in\can$ is new'',
\item[(b)] if $\bP$ is a Souslin c.c.c.\ forcing notion with no $\bP$--name
$\name{r}'$ of a real such that the forcing notion ${\cal B}_{\bP}(\name{r}')$
is $\baire$--bounding but $\bP$ adds a nowhere essentially Cohen real\\ 
{\em then} forcing with $\bP$ adds a $(\bQ,\name{r})$ real,i.e.\ for some
$\bP$--name $\name{r}''$ for a real we have $\forces_\bP$``for some $G''
\subseteq\bQ^\V$ generic over $\V$, $\name{r}''[G_\bP]=\name{r}[G'']$''.
\end{enumerate}
\item As above $\bP$ is $\sigma$--centered.
\item If $\bP$ is a Souslin c.c.c.\ forcing notion adding new reals but not
adding a real $\name{r}'$ with ${\cal B}_\bP(\name{r}')$ being
$\baire$--bounding,\\ 
{\em then} forcing with $\bP$ adds a new real $\name{r}''$ such that ${\cal
B}_\bP(r'')$ is $\sigma$--centered. 
\end{enumerate}
\end{problem}

\begin{problem}
\label{20.10}
\begin{enumerate}
\item Let $\bQ$ be a Souslin c.c.c.\ forcing notion and $\forces_\bQ$``$
\name{r}\in \can$''. Is ${\cal B}(\name{r})$ also a Souslin c.c.c.\ forcing
notion? 
\item Similarly for nep c.c.c.
\end{enumerate}
\end{problem}

\begin{problem}
\label{20.11}
Assume $\bQ$ is a Souslin c.c.c.\ forcing notion which is snep and even $``x
\in\bQ$'', ``$x \le^{\bQ} y$'', ``$\{p_n:n<\omega\}$ predense above $q'$'' are 
$\Sigma^1_1$ relations.  Does $\bQ$ add Cohen or random real?
\end{problem}

\begin{problem}
\label{20.12}
Develop the theory of ``definable forcing notions'' when we allow an
ultrafilter on $\omega$ as a parameter.
\end{problem}

\begin{problem}
\label{20.13}
Does nep$\neq$snep? (the case $\theta=\kappa=\aleph_0$, of course).
\end{problem}

\begin{problem}
\label{20.14}
Try to generalize our present context to $\lambda$--complete forcing notions
(Baumgartner's Axiom; \cite{Sh:186}, \cite{Sh:655}).
\end{problem}

\begin{problem}
\label{20.15}
When $\bQ^\V\lesdot\bQ^{\V^\bP}$?
\end{problem}

\begin{problem}
\label{20.19}
Does ${\rm Ax}_{\omega_1}[(\aleph_1,\aleph_1)\mbox{--nep}]$ imply
$2^{\aleph_0}=\aleph_2$? 

Or does ${\rm Ax}_{\omega_1}[\mbox{nep}]$ imply $2^{\aleph_0}=\aleph_2$?

\noindent [The parallel question for Souslin proper was formulated in xxx]
\end{problem}
\shlhetal

\end{document}